\documentclass[12pt, english, reqno, oneside]{smfart}

\usepackage[utf8]{inputenc}
\usepackage{float}
\usepackage{enumitem}
\usepackage{colonequals}
\usepackage{textcomp}
\usepackage{amsmath,amssymb,amsthm,bm,mathtools}
\usepackage[bottom=32mm, top=32mm , left=30mm , right=30mm]{geometry}
\usepackage{rotating}
\usepackage{makecell}
\usepackage{graphicx}
\usepackage{subcaption}
\usepackage{url}
\usepackage[colorlinks=true,
            linkcolor=blue,
            filecolor=blue,
            citecolor=red,
            urlcolor=blue, 
            colorlinks=true]{hyperref}
\usepackage[capitalize,noabbrev]{cleveref}
\usepackage{csquotes}
\usepackage{crossreftools}
\usepackage{multirow}
\usepackage{tikz}
\usepackage{tikz-cd}
\usepackage{comment}
\usepackage{hyphenat}
\usetikzlibrary{arrows}
\usetikzlibrary{decorations.pathmorphing}
\usetikzlibrary{decorations.markings}
\usetikzlibrary{arrows.meta,bending}
\usetikzlibrary{decorations.pathreplacing,calligraphy}
\usetikzlibrary{shapes.geometric}
\usetikzlibrary{patterns}

\newtheorem{maintheorem}{Theorem}[section]
\renewcommand{\themaintheorem}{\Alph{maintheorem}} 

\newtheorem{theorem}{Theorem}[section]
\newtheorem{proposition}[theorem]{Proposition}
\newtheorem{lemma}[theorem]{Lemma}
\newtheorem{corollary}[theorem]{Corollary}
\newtheorem{conjecture}[theorem]{Conjecture}

\theoremstyle{definition} 
\newtheorem{definition}[theorem]{Definition}

\newtheorem{remark}[theorem]{Remark}
\newtheorem{convention}[theorem]{Convention}
\newtheorem{example}[theorem]{Example}
\newtheorem{assumption}[theorem]{Assumption}
\newtheorem{question}[theorem]{Question}

\definecolor{qqqqff}{rgb}{0,0,1}
\definecolor{ggreen}{RGB}{0, 94, 25}
\definecolor{claudiaColor}{RGB}{102, 0, 204}

\newcommand{\yassine}[1]{{\color{blue} \sf $\clubsuit$ Yassine: [#1]}}
\newcommand{\paul}[1]{{\color{orange} \sf $\clubsuit$ Paul: [#1]}}
\newcommand{\Felix}[1]{{\color{magenta} \sf Felix: [#1]}}

\newcommand{\dR}{{\operatorname{dR}}}
\newcommand{\CC}{{\operatorname{CC}}}
\newcommand{\an}{{\operatorname{an}}}
\newcommand{\CG}{{\operatorname{CG}}}
\newcommand{\PP}{{\operatorname{PP}}}
\newcommand{\GL}{{\operatorname{GL}}}
\newcommand{\PGL}{{\operatorname{PGL}}}
\newcommand{\inv}{{\operatorname{inv}}}
\newcommand{\Hom}{{\operatorname{Hom}}}
\newcommand{\trop}{{\operatorname{trop}}}
\newcommand{\Jac}{{\operatorname{Jac}}}
\newcommand{\Berk}{{\operatorname{Berk}}}
\newcommand{\Spec}{{\operatorname{Spec}}}
\newcommand{\Orb}{{\operatorname{Orb}}}
\newcommand{\Supp}{{\operatorname{Supp}}}
\newcommand{\red}{{\operatorname{red}}}
\newcommand{\res}{{\operatorname{res}}}
\newcommand{\Res}{{\operatorname{Res}}}
\newcommand{\Symb}{{\operatorname{Symb}}}
\newcommand{\Log}{{\operatorname{Log}}}
\newcommand{\Lie}{{\operatorname{Lie}}}
\newcommand{\gl}{{\operatorname{gl}}}
\newcommand{\Div}{{\operatorname{Div}}}
\renewcommand{\div}{{\operatorname{div}}}
\newcommand{\Pic}{{\operatorname{Pic}}}
\newcommand{\ch}{{\operatorname{char}}}
\newcommand{\SL}{{\operatorname{SL}}}
\newcommand{\vol}{{\operatorname{vol}}}
\renewcommand{\mod}{{\operatorname{mod}}}
\newcommand{\dil}{\mathrm{dil}}
\newcommand{\free}{\rm free}

\newcommand{\longhookrightarrow}{\lhook\joinrel\longrightarrow}

\newcommand{\bbGamma}{{\mathpalette\makebbGamma\relax}}
\newcommand{\makebbGamma}[2]{%
  \raisebox{\depth}{\scalebox{1}[-1]{$\mathsurround=0pt#1\mathbb{L}$}}%
}

\DeclareMathOperator{\id}{id}					
\newcommand{\op}{\mathrm{op}}					
\DeclareMathOperator{\Aut}{Aut}					
\DeclareMathOperator{\Stab}{Stab}					
\DeclareMathOperator{\val}{val}					
\DeclareMathOperator{\Star}{Star}				

\newcommand{\lra}{\longrightarrow}
\newcommand  \scdot{\,\cdot\,}
\newcommand*\angles[1]{\langle #1 \rangle}
\newcommand{\catname}[1]{\upshape{\textbf{#1}}}

\newcommand*\dsum{\oplus}
\newcommand*\Dsum{\bigoplus}

\newcommand*\N{\mathbb{N}}
\newcommand*\Z{\mathbb{Z}}
\newcommand*\Q{\mathbb{Q}}
\newcommand*\R{\mathbb{R}}
\newcommand*\C{\mathbb{C}}
\newcommand*\F{\mathbb{F}}

\newcommand*\G{\mathbb{G}}
\newcommand*\A{\mathbb{A}}
\renewcommand*\P{\mathbb{P}}
\newcommand*\T{\mathbb{T}}

\newcommand*\Zp{\Z_p}
\newcommand*\Qp{\Q_p}
\newcommand*\Cp{\C_p}

\newcommand*\cA{\mathcal{A}}
\newcommand*\cB{\mathcal{B}}
\newcommand*\cC{\mathcal{C}}
\newcommand*\cD{\mathcal{D}}
\newcommand*\cF{\mathcal{F}}
\newcommand*\cH{\mathcal{H}}
\newcommand*\cL{\mathcal{L}}
\newcommand*\cO{\mathcal{O}}
\newcommand*\cP{\mathcal{P}}
\newcommand*\cM{\mathcal{M}}
\newcommand*\cN{\mathcal{N}}
\newcommand*\cR{\mathcal{R}}

\newcommand{\MS}{\Delta_{g,p}}

\newcommand{\lin}{{\mathrm{lin}}}
\newcommand{\nlin}{{\mathrm{nlin}}}
\newcommand{\aff}{{\mathrm{aff}}}
\newcommand{\st}{{\mathrm{st}}}
\newcommand{\sep}{\mathrm{sep}}
\newcommand{\et}{\mathrm{et}}

\newcommand*\mfrak{\mathfrak{m}}

\newcommand*\wt{\widetilde}
\renewcommand{\tilde}{\widetilde}		
\newcommand*\wb{\overline}
\newcommand{\bigmid}{\mathrel{\big\vert}}

\newcommand{\vertex}[3][-1]{\fill (#2, #3) circle (2pt);%
	\ifthenelse{\equal{#1}{-1}}{}{\vertexgenus{#2}{#3}{#1}}}

\newcommand{\vertexgenus}[3]{\draw (#1, #2 + 0.3) node {$#3$};}

\newcommand{\fatvertex}[3][-1]{\fill (#2, #3) circle (3pt);%
	\ifthenelse{\equal{#1}{-1}}{}{\vertexgenus{#2}{#3}{#1}}}

\newcommand{\sqvertex}[3][-1]{\node[xscale=0.4,yscale=0.4,diamond,draw,fill=black] (d) at (#2, #3) {};%
    \ifthenelse{\equal{#1}{-1}}{}{\vertexgenus{#2}{#3}{#1}}}

\newcommand{\selfloop}[3][0.25]{\draw [line width= #1 mm] (#2, #3) .. controls (#2 + 1, #3 + 1) and (#2 + 1, #3 - 1) .. (#2, #3);}
\newcommand{\selfloopleft}[3][0.25]{\draw [line width= #1 mm] (#2, #3) .. controls (#2 - 1, #3 + 1) and (#2 - 1, #3 - 1) .. (#2, #3);}

\newcommand{\dilatedselfloop}[4][0.25]{\draw [line width= #1 mm] (#2, #3) .. controls (#2 + 1, #3 + 1) and (#2 + 1, #3 - 1) .. (#2, #3);}

\newcommand{\orient}[1]{#1 node[currarrow,
   		pos=0.5, 
   		xscale=-1,
   		sloped,
   		scale=2] {};}

\newcommand{\thinarrowup}[1]{
    \draw #1 -- ++ (-0.15, -0.15);
    \draw #1 -- ++ (0.15, -0.15);
}
\newcommand{\thinarrowright}[1]{
    \draw #1 -- ++ (-0.15, -0.15);
    \draw #1 -- ++ (-0.15, 0.15);
}

\newcommand{\fatarrowup}[1]{
    \draw[-triangle 90] #1 -- +(0, 0.15);
}
\newcommand{\fatarrowright}[1]{
    \draw[-triangle 90] #1 -- +(0.15, 0);
}

\newcommand{\verticaldots}[2]{\fill (#1, #2) circle (0.6pt);%
	\fill (#1, #2 - 0.07) circle (0.6pt);%
	\fill (#1, #2 + 0.07) circle (0.6pt);}

\newcommand{\butterfly}[2]{\draw (#1, #2) .. controls (#1 + 1, #2 + 1.2) and (#1 + 1, #2 + 0.2) .. (#1, #2);%
	\verticaldots{#1 + 0.6}{#2}
	\draw (#1, #2) .. controls (#1 + 1, #2 - 1.2) and (#1 + 1, #2 - 0.2) .. (#1, #2);}

\newcommand{\smallspiral}[2]{\smallsp{#1}{#2}{1}}
\newcommand{\bigspiral}[2]{\bigsp{#1}{#2}{1}}
\newcommand{\smallspiralleft}[2]{\smallsp{#1}{#2}{-1}}
\newcommand{\bigspiralleft}[2]{\bigsp{#1}{#2}{-1}}

\newcommand{\smallsp}[3]{%
	\vertex{#1}{#2}
	\vertex{#1}{#2 + 0.4}
	\vertex{#1}{#2 + 0.8}
	\begin{scope}
		\clip (#1, #2 - 0.2) rectangle (#1 + #3, #2 + 1);
		\foreach \shift in {-0.4, 0, 0.4, 0.8}{
			\draw (#1, #2 + \shift) 
				.. controls (#1 + #3, #2 - 0.5 + \shift) and (#1 + #3, #2 + 0.4 + 0.5 + \shift) .. 
				(#1, #2 + 0.4 + \shift);
		}
	\end{scope}%
}

\newcommand{\bigsp}[3]{%
	\vertex{#1}{#2}
	\vertex{#1}{#2 + 0.4}
	\vertex{#1}{#2 + 0.8}
	\begin{scope}
		\clip (#1, #2 - 0.2) rectangle (#1 + #3, #2 + 1);
		\foreach \shift in {-0.8, 0, 0.8}{
			\draw (#1, #2 + \shift) 
			.. controls (#1 + #3, #2 - 0.5 + \shift) and (#1 + #3, #2 + 0.8 + 0.5 + \shift) .. 
			(#1, #2 + 0.8 + \shift);
		}
	\end{scope}
	\begin{scope}
		\clip (#1, #2) rectangle (#1 + #3, #2 + 0.8);
		\foreach \shift in {-0.4, 0.4}{
			\draw (#1, #2 + \shift) 
			.. controls (#1 + #3, #2 - 0.5 + \shift) and (#1 + #3, #2 + 0.8 + 0.5 + \shift) .. 
			(#1, #2 + 0.8 + \shift);
		}
	\end{scope}%
}

\makeatletter
\long\def\my@drawfill#1#2;{%
	\@skipfalse
	\fill[#1,draw=none] #2;
	\@skiptrue
	\draw[#1,fill=none] #2;
}

\newif\if@skip

\newcommand{\skipit}[1]{%
	\if@skip
	\else
	#1
	\fi
}

\newcommand{\drawfill}[1][]{%
	\my@drawfill{#1}}
\makeatother

\title[Topology of $\Delta_{g, \Z/p}$]{On the topology of the moduli of tropical unramified $p$-covers}

\author[Y.~El~Maazouz]{\small Yassine~El~Maazouz}
\address{(RWTH~Aachen~University)}
\email{yassine.el-maazouz@rwth-aachen.de}

\author[P.~A.~Helminck]{\small Paul~Alexander~Helminck}
\address{(University~of~Tsukuba)}
\email{p.a.helminck@math.tsukuba.ac.jp}

\author[F.~R\"ohrle]{\small Felix~R\"ohrle}
\address{(University~of~T\"ubingen)}
\email{roehrle@math.uni-tuebingen.de}

\author[P.~Souza]{\small Pedro~Souza}
\address{(Goethe~University~Frankfurt)}
\email{souza@math.uni-frankfurt.de}

\author[C.~H.~Yun]{\small Claudia~He~Yun}
\address{(University~of~Michigan)}
\email{clyun@umich.edu}


\begin{document}

\begin{abstract}
    We study the topology of the moduli space of unramified $\mathbb{Z}/p$-covers of tropical curves of genus $g \geq 2$, where $p$ is a prime number. We use recent techniques by Chan--Galatius--Payne to identify contractible subcomplexes of the moduli space. We then use this contractibility result to show that this moduli space is simply connected. In the case of genus $2$, we determine the homotopy type of this moduli space for all primes $p$. This work is motivated by prospective applications to the top-weight cohomology of the space of prime cyclic \'etale covers of smooth algebraic curves.
\end{abstract}

\maketitle

\setcounter{tocdepth}{1}

\tableofcontents

\section{Introduction}

The moduli space $\cM_g$ of smooth algebraic curves of genus $g$ and the Deligne-Mumford compactification $\overline{\cM}_g$ are arguably some of the most intricate objects in modern algebraic geometry and have been studied from several perspectives \cite{Grot60, DM69, HarrisMumford82, Mumford83, EisenbudHarris87, ACP15}. Closely related to $\cM_g$ is the moduli space $\cM_{g, \Z/p}$ parametrizing isomorphism classes of \'etale $\Z/p$-Galois covers $\tilde C \to C$ with $C$ a smooth curve of genus $g$, see for example \cite[Sections~5 and~6]{BR11} and \cite{ACV2003}. On the level of coarse moduli spaces, this is the same as the moduli space of smooth curves $C$ of genus $g$ together with a non-trivial $p$-torsion point of their Jacobian.
The space $\cM_{g, \Z/p}$ is itself a smooth \'etale cover of $\cM_g$ of degree $p^{2g} - 1$, see \cite[Lemma~5.7]{DM69} or \cite[Chapter~XVI]{ACG_II} for a more recent treatment.

One of the more recent insights into the topology of $\cM_g$ was achieved in \cite{CGP1}, where the authors introduce new tropical methods to investigate the cohomology of $\cM_g$ and in particular find new non-zero classes in the unstable part of its cohomology. The key observation in this story was the comparison theorem \cite[Theorem~1.2]{CGP1} which identifies the top-weight cohomology of $\cM_g$ with the reduced rational homology of the dual boundary complex of its compactification $\cM_g \subset \overline{\cM}_g$. The latter in turn, by virtue of \cite[Theorem~1.2.1]{ACP15}, can be identified with the moduli space $\Delta_g$ parametrizing stable tropical curves of genus $g$ of volume one. 
This motivates a more careful study of the topology of $\Delta_g$. 
To this end, \cite{CGP2} develops a general framework to study loci in so-called \emph{symmetric $\Delta$-complexes}, of which the tropical moduli space $\Delta_g$ is an example. In particular, the authors provide criteria to prove contractibility of subcomplexes which are defined by \emph{vertex properties}.

\medskip

In this article, we study the topology of a tropical moduli space $\Delta_{g, \Z/p}$, which is the tropical counterpart to $\cM_{g,\Z/p}$ as $\Delta_g$ is to $\cM_g$. Throughout this article, a tropical curve will mean a finite metric graph $\Gamma$ together with a vertex genus function $g \colon V(\Gamma) \to \Z_{\geq 0}$. Following \cite{Galeotti, BCK23} we define \emph{tropical $\Z/p$-covers} as harmonic morphisms $\pi\colon \tilde \Gamma \to \Gamma$ of tropical curves, which are unramified in the sense that they satisfies the local Riemann-Hurwitz condition everywhere, and which are additionally equipped with a $\Z/p$-action on $\tilde \Gamma$ and a \emph{dilation flow} on $\Gamma$ (see \cref{def:Gcover}). This definition is natural, because these are precisely the combinatorial objects that arise when degenerating \'etale $\Z/p$-covers of smooth curves, see~\cref{subsec:relatedWork}.

The moduli space $M^\trop_{g, \Z/p}$ parametrizing tropical $\Z/p$-covers is naturally a generalized cone complex in the sense of \cite{ACP15}. Its link $\Delta_{g,\Z/p}$ around the cone point, which is a symmetric $\Delta$-complex, is the central object of study in this article. Its study is motivated by prospective applications to the top-weight cohomology of the moduli space $\mathcal{M}_{g,\Z/p}$.

\subsection{Results}
For ease of notation, we shall from now on write $\Delta_{g,p}$ to denote $\Delta_{g,\Z/p}$. We define the following nested loci in $\Delta_{g,p}$:

\begin{enumerate}
	\item The \emph{equivariant weight locus} $\MS^{w}$ is defined as the locus (see also \cref{lem:weight_locus_criterion}):
	\[ \Big\{ \pi \colon \tilde \Gamma \to \Gamma \mathrel{\Big |} \text{there exists } v \in V(\Gamma) \text{ such that } g(v) \geq 0\text{ and }\sum_{\tilde v \in \pi^{-1}(v)} g(\tilde v) \geq p \Big\}. \]
    
	\item The \emph{equivariant loop or weight locus} $\MS^{lw}$ is defined as the closure of the locus of covers  $\pi \colon \tilde \Gamma \to \Gamma$ such that $\pi \in \MS^{w}$ or there is a loop-edge $e$ in $\Gamma$ such that $\pi^{-1}(e)$ consists of $p$ loop edges. 
 
	\item The \emph{equivariant bridge locus} $\MS^{br}$ is the closure of the union of $\MS^{lw}$ and
	\begin{equation*}
		\left\{ 
		\pi \colon \tilde \Gamma \to \Gamma \mathrel{\bigg |} 
		\begin{minipage}{0.65\textwidth}
			there is a bridge-edge $e \in \Gamma$ such that $|\pi^{-1}(e)| = p$ and every $\tilde e \in \pi^{-1}(e)$ is a bridge as well
		\end{minipage}		
		\right\}.
	\end{equation*}
	
	\item The \emph{sparsely connected locus} $\MS^{scon}$ is the closure of the union of $\MS^{br}$ and
	\begin{equation*}
		\left\{
		\pi \colon \tilde \Gamma \to \Gamma \mathrel{\bigg |}
		\begin{minipage}{0.65\textwidth}
			there is an edge $e \in \Gamma$ and a connected component $\Gamma_0 \subseteq \Gamma \setminus \{e\}$ such that $\pi^{-1}(\Gamma_0)$ is disconnected
		\end{minipage}
		\right\}.
	\end{equation*}

	\item The \emph{equivariant parallel edge locus} $\MS^{par}$ is the closure of the union of $\MS^{scon}$ and
	\begin{equation*}
		\left\{
		\pi \colon \tilde \Gamma \to \Gamma \mathrel{\Bigg |} \ 
		\begin{minipage}{0.75\textwidth}
			there exists a pair of parallel edges $e_1, e_2 \in \Gamma$ such that
			$|\pi^{-1}(e_1)| = |\pi^{-1}(e_2)| = p$ and every lift of the simple closed cycle defined by $e_1$ and $e_2$ to $\tilde \Gamma$ is again a simple closed cycle
		\end{minipage}
		\ \right\}.
	\end{equation*}
\end{enumerate}
These loci are nested as follows:
\[ \MS^{w} \subseteq \MS^{lw} \subseteq \MS^{br} \subseteq \MS^{scon} \subseteq \MS^{par}.  \]
We illustrate the moduli spaces $\Delta_{2,2}$ and $\Delta_{2,3}$, as well as the loci defined above in Figures~\ref{fig:Delta_2_2}, \ref{fig:Delta_2_3}, and \ref{fig:loci}, respectively. Our first result provides a very beneficial simplification for further study of the topology of $\MS$.

\begin{maintheorem} \label{thm:contractible_loci}
    For all $g \geq 2$ and $p$ a prime number, the loci $\MS^{w}$, $\MS^{lw}$, $\MS^{br}$ are contractible. Moreover, when $p > 2$, the loci $\MS^{scon}$ and $\MS^{par}$ are contractible as well.
\end{maintheorem}

We prove Theorem~\ref{thm:contractible_loci} in Section~\ref{sec:contractible_loci}. 
Except for the sparsely connected locus, all the loci we list above are generalizations of subcomplexes with similar properties in $\Delta_g$. The contractibility of these loci follows from similar arguments to \cite{CGP2}; we give new ideas for the sparsely connected locus. The subcomplex $\MS^{par}$ is the largest locus in $\Delta_{g,p}$ that one can prove to be contractible with the vertex property techniques from \cite{CGP2}, see \cref{rem:moreLoci} for a more detailed discussion.

\medskip

As an application of \cref{thm:contractible_loci} we show the following:

\begin{maintheorem} \label{thm:simply_connected}
	The moduli space $\MS$ is simply connected for all $g \geq 2$ and $p$ prime.
\end{maintheorem}

Theorem~\ref{thm:simply_connected} is proved in Section~\ref{sec:simply_connected}. The key idea is to use \cite[Theorem~3.1]{AllcockCoreyPayne}, which states that the fundamental group of the geometric realization of $\MS$ is generated by cycles supported on the 1-skeleton $\MS^{(1)}$. This reduces the number of combinatorial types of $\Z/p$-covers to be considered. 

\medskip

In \cref{sec:genus_2}, we turn to the special case $g = 2$ and explicitly compute the homotopy type of $\Delta_{2, p}$ for any prime number $p$.

\begin{maintheorem} \label{thm:MainThmGenus2}
    The number of maximal cells of $\Delta_{2,2}$ is $7$ and $\Delta_{2,2}$ is contractible. The number of maximal cells of $\Delta_{2,3}$ is $9$ and $\Delta_{2,3}$ is contractible as well. For $p \geq 5$ a prime number, the number of maximal cells in $\Delta_{2,p}$ is 
    \[
        \frac{4p^2 +9p  - 13}{6}
    \]
    and $\Delta_{2,p}$ has the homotopy type of a wedge of 
        \[
            \frac{(p-1)(p-5)}{6} + \frac{(p-3)^2}{4}
        \]
        spheres of dimension $2$.
\end{maintheorem}

Finally we remark that in this paper we only consider unramified tropical $\Z/p$-covers. Naturally one could ask about an extension of our results to the more general case of $\Z/p$-covers with prescribed ramification profiles over marked points of the target. We believe that our results and techniques carry over to this setting, and that there should be an analogue of the contractibility result for the \emph{repeated marking locus} \cite[Theorem~1.1.3]{CGP2}. However, at this stage, we refrain from elaborating more on this.

\subsection{Motivation and historical background}

One difficulty in studying $\mathcal{M}_{g}$ and $\overline{\mathcal{M}}_{g}$ arises from the fact that the parametrized objects can have nontrivial automorphisms. To untangle these automorphisms, one can introduce auxiliary moduli spaces $\mathcal{M}_{g}(N)$ as in \cite{MumfordPicard1965} and \cite{DM69}, which parametrize smooth curves $C$ together with a basis of the $N$-torsion of the Jacobian of $C$. For $N\geq{3}$, these spaces are always schemes, and this allows one to study the stacks $\mathcal{M}_{g}$ as finite quotients of $\mathcal{M}_{g}(N)$. In the special case when $g=1$, the moduli spaces $\mathcal{M}_{g}(N)$ are better known as \emph{modular curves}\footnote{Technically speaking, one has to add a marked point as well, so that we are considering $\mathcal{M}_{1,1}(N)$ rather than $\mathcal{M}_{1}(N)$.}, see \cite{DR1973, KatzMazur1985}. Here, the coarse moduli spaces of the $\mathcal{M}_{g}(N)$ are often denoted by $Y(N)$, and their natural compactifications using generalized elliptic curves are denoted by $X(N)$. The modular curve that describes generalized elliptic curves with a specified $p$-torsion point is the quotient of $X(p)$ by the group of unipotent upper-triangular matrices in $\mathrm{PSL}_{2}(\Z/p)$, and it is usually denoted by $X_{1}(p)$. We similarly view the moduli space $\mathcal{M}_{g,\Z/p}$ under consideration in this paper as a quotient of $\mathcal{M}_{g}(p)$ by a subgroup of the group $\mathrm{Sp}_{2g}(\Z/p)$ of symplectic matrices over $\Z/p$.

\subsection{Applications}

We now describe prospective results on the (top-weight) cohomology of $\cM_{g,\Z/p}$. This is based on the following general fact. Let $\mathcal{X}$ be a smooth and separated Deligne-Mumford stack of dimension $d$ over $\mathbb{C}$. By \cite{Deligne75}, there is a so-called \emph{weight filtration} 
\begin{equation*}
    W_{0}H^{k}(\mathcal{X};\mathbb{Q})\subseteq \cdots \subseteq W_{2d}H^{k}(\mathcal{X};\mathbb{Q})=H^{k}(\mathcal{X};\mathbb{Q})
\end{equation*}
for every $k\geq{0}$. One then considers the graded pieces
\begin{equation*}
    \mathrm{Gr}^{W}_{j}H^{k}(\mathcal{X};\mathbb{Q}) \coloneqq W_{j}H^{k}(\mathcal{X};\mathbb{Q})/W_{j-1}H^{k}(\mathcal{X};\mathbb{Q})
\end{equation*}
where $\mathrm{Gr}^W_{2d}H^k(\mathcal{X}; \mathbb{Q})$ is called \emph{top-weight cohomology} of $\mathcal{X}$ in degree $k$. Now let $\mathcal{X} \subseteq \overline{\mathcal{X}}$ be a normal crossings compactification, $D = \overline{\mathcal{X}}\setminus\mathcal{X}$ its boundary, and denote the dual boundary complex by $\Delta(D)$ (see \cite[Definition 5.2]{CGP1}). By \cite[Theorem 5.8]{CGP1}, we have an isomorphism
\begin{equation} \label{eq:topweight_cohomology}
    \mathrm{Gr}^{W}_{2d}H^{2d-k}(\mathcal{X};\mathbb{Q})\simeq H_{k-1} \big(\Delta(D);\mathbb{Q} \big).
\end{equation}

In order to apply this to $\mathcal{X}=\mathcal{M}_{g,\Z/p}$, we recall that there is a natural compactification of $\cM_{g,\Z/p}$, the moduli space of admissible $\Z/p$-covers $\overline{\cM}_{g,\Z/p}$, whose construction goes back to ideas in \cite[Section~4]{HarrisMumford82} and is a special case of \cite{AbramovichVistoli2002}. Write $D=\overline{\cM}_{g,\Z/p}\backslash \cM_{g,\Z/p}$ for the boundary of $\cM_{g,\Z/p}$ in this compactification, which is a normal crossings divisor with dual intersection complex $\Delta(D)$. The dimension of $\cM_{g, \Z/p}$ is $d = 3g-3$, and one can state the following analogue of \cite[Theorem~1.2.1]{ACP15}, whose proof is subject to ongoing work by the fourth author. The ideas behind this proof are based on the description of the boundary of $\mathcal{M}_{g,\Z/p}$ in the admissible covers compactification, see \cite[Section 7.4]{BR11}.

\begin{conjecture} \label{conj:PedrosTheorem}
 There is an isomorphism of symmetric $\Delta$-complexes
    \begin{equation*}
        \Delta(D)\longrightarrow \MS.
    \end{equation*}
\end{conjecture}

Conjecture~\ref{conj:PedrosTheorem} motivates the study of the topology of 
$\MS$, and combining this with Equation~\eqref{eq:topweight_cohomology} we obtain the following:

\begin{corollary}\label{cor:topweight_cohomology}
    Conditional on \cref{conj:PedrosTheorem}, we have for all $p\geq 5$:
    \[
        \mathrm{dim}_{\mathbb{Q}} H^{3}(\mathcal{M}_{2,\Z/p};\mathbb{Q})\geq \frac{(p-1)(p-5)}{6} + \frac{(p-3)^2}{4}.
    \]
\end{corollary}

\subsection{Related work}\label{subsec:relatedWork}

Tropical methods have been employed to study a number of moduli spaces: the moduli space of principally polarized abelian varieties was studied in \cite{BBCMGC22}, and more recently \cite{BCK23} study the moduli space of hyperelliptic curves. 
\cite{CMR16} defines a tropical analogue of the Hurwitz space $\mathcal{H}_{g, d}(\mu)$ parametrizing degree $d$ covers of genus $g$ curves with prescribed ramification profiles $\mu = (\mu_1, \ldots, \mu_n)$ at marked points. This is based on the notion of harmonic morphism stemming from analytic considerations, see \cite{ABBR15}. 
However it was shown in \cite{CMR16} that the natural morphism relating the dual boundary complex of the admissible cover compactification of the algebraic Hurwitz space on one hand, and tropical Hurwitz space on the other, is neither injective nor surjective. The latter is because of the \emph{Hurwitz existence problem}, see \cite[Section 1]{PP06}. Therefore, the topology of $H^{\trop}_{g, d}(\mu)$ cannot be directly used to compute the top-weight cohomology of $\mathcal{H}_{g, d}(\mu)$.

We now discuss two tropical spaces in the literature that mimick part of the moduli space of $G$-covers. In \cite{LUZ19} a moduli space of tropical $G$-covers is defined for abelian groups $G$, but this space does not take the dilation flows into account. This space thus inherits the shortcomings of the space of tropical Hurwitz covers. Similarly, in \cite{APS_roots} the space $\mathrm{Root}^{p}_{g}$ of tropical curves with genus $g$ together with a $p$-torsion point in the tropical Jacobian is studied. Note that by \cite[Theorem 5.16]{Helminck2023}, we have an identification between $p$-torsion points in the tropical Jacobian and algebraic unramified $\Z/p$-covers that induce a tropical cover with dilation over at least one edge. We thus see that the other missing covers in this case are the ones that correspond to covers with no dilation.

The tropical moduli space that we consider in this paper essentially combines these two points of view by parametrizing tropical $\Z/p$-covers with dilation flows. From the non-archimedean point of view this space is natural because the skeleton of the non-archimedean analytification of an \'etale $\Z/p$-cover of algebraic curves is a tropical $\Z/p$-cover in our sense. Here the dilation data can be interpreted as additional gluing data on the annuli corresponding to the edges, see \cite[Section 7.1]{ABBR15} and \cite[Section 4.1]{Helminck2023}. In terms of algebraic geometry, a one-parameter degeneration of a smooth $G$-curve is a nodal curve, and the stabilizer subgroups of the nodes act on their tangent spaces. The dilation flow is the combinatorial shadow of this data, see \cite[Section 7.4]{BR11}. Conversely, we note that one we can use the material in \cite{Helminck2023} to see that the stable $\Z/p$-covers considered in this paper are always realizable as algebraic $\Z/p$-covers: The conditions on the dilation flow allow one to define an automorphism of order $p$ of enhanced coverings of metrized complexes as in \cite[Section 4.1]{Helminck2023}, and the equivalence of Galois categories in that paper implies that this corresponds to a $\Z/p$-cover of algebraic curves. A version of this argument also appeared in the proof of \cite[Theorem 4.4]{LUZ19}. In summary this shows that our tropical moduli space is expected to be the correct space for top-weight cohomology applications.

\begin{figure}[ht]
    \begin{center}
    \begin{tikzcd}
        & M^{\trop}_{g, \Z/p} \arrow[->>, dl] \arrow[->>, dr] & & & \\
        \operatorname{Root}_g^p \arrow[->>, dr] & & \big(H_{g, \Z/p}^\trop(\emptyset)\big)^{\operatorname{real}} \arrow[->>, dl] \arrow[hook, r] & H_{g, \Z/p}^\trop(\emptyset) \arrow[r, "(*)"] & H^\trop_{g, d}(\emptyset) \arrow[->>, dlll, bend left = 10] \\
        & M_g^{\trop} & & & 
    \end{tikzcd}
\end{center}
    \caption{Summary of the different moduli spaces in the literature related to $M^{\trop}_{g,\Z/p}$. Here we denote by \enquote{real} the realizability locus and $(*)$ is the map which sends a $\Z/p$-cover to its underlying degree $d = |\Z/ p| = p$ Hurwitz cover.}
    \label{fig:diagramDiffModuli}
\end{figure}
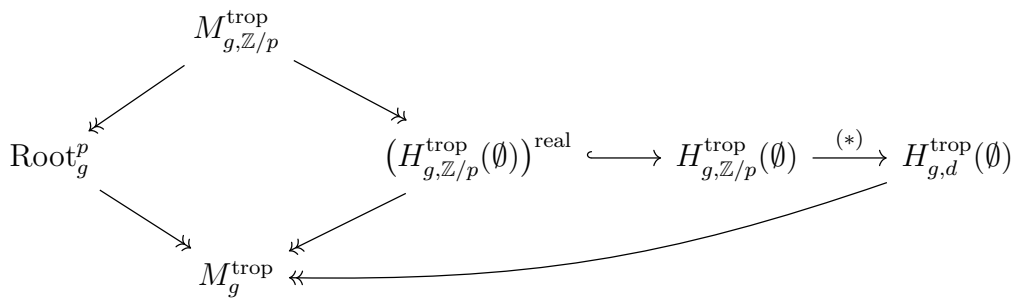

The relation between all the moduli spaces mentioned above is summarized in \cref{fig:diagramDiffModuli}.

\begin{remark}
    To find a tropical space that reflects the top-weight cohomology of $\mathcal{H}_{g,d}(\mu)$, one can consider the branch locus map $br \colon \overline{\mathcal{H}}_{g,d}(\mu)\to \overline{\mathcal{M}}_{g,b}$. This map is \'{e}tale over the smooth locus and, away from the primes dividing the order of the Galois group, it is tamely ramified over the boundary. Since the boundary of $\mathcal{M}_{g,b}$ is a normal crossings divisor, it directly follows, from Abhyankar's Lemma and the purity of the branch locus as in \cite[Chapter XIII, Appendice 1]{SGA1}, that the inverse image is a normal crossings divisor after normalizing. The normalization of this moduli space can be interpreted as the space of twisted bundles, see \cite{ACV2003}. When the covers are Galois, this space is already normal \cite[Theorem 4.3.2]{ACV2003}, so that the notion of a twisted cover is unnecessary. All in all we see that we only need to characterize the inverse images of these divisors and their self-intersections to obtain the correct boundary complex.

    For Galois covers, we can use the results in \cite[Section 7.4]{BR11}, which gives a description of the different components in the inverse image in terms of tropical $G$-covers enhanced with compatible dilation data, which are called \emph{graphes modulaires de Hurwitz}.
\end{remark}

\begin{remark}\label{rem:luzSpace}
    The work in the present article started with an analysis of the space $H^{\trop}_{g,\Z/p}(\emptyset)$ in the sense of \cite{LUZ19} and the results presented here were first found to hold for that space, as well as the realizable locus within it (the only exception being the computational result in Theorem~\ref{thm:MainThmGenus2}, where the correct formula for the homotopy type of $H^{\trop}_{2,\Z/p}(\emptyset)$ is a wedge of $(p-1)(p-5)/12$ spheres). Since the proofs for $H^{\trop}_{2, \Z/p}(\emptyset)$ are verbatim as presented here, we refrain from elaborating any further.
\end{remark}

\section{The moduli space of tropical \texorpdfstring{$G$}{}-covers } \label{sec:Background_ModuliSpace}

Let $g \geq 2$ and fix a finite abelian group $G$.
In this section we introduce the moduli space of tropical $G$-covers of genus $g$ and the objects it parametrizes. The space $\MS$, which is the main player in the rest of the article, is obtained in the special case $G = \Z/p$. We first review the notions of graphs and their morphisms.

\subsection{Weighted graphs}

A graph is a tuple ${\Gamma} =(X, \iota, r)$ consisting of:
\begin{itemize}
    \item a finite set $X$,
    \item an idempotent \emph{root map} $r \colon X \to X$; that is $r\circ r = r$, and
    \item an involution $\iota \colon X \to X$; that is $\iota \circ \iota = \id$,
\end{itemize}
such that $\iota \circ r = r$ and the fixed points of $\iota$ are exactly the elements in the image of $r$.
The set $V(\Gamma) = r(X)$ is the set of \emph{vertices} of $\Gamma$ and $H(\Gamma) = X \setminus V(\Gamma)$ denotes the set of \emph{half-edges} of $\Gamma$. The edge set $E(\Gamma)$ of $\Gamma$ is defined by identifying half-edges that are exchanged by the involution $\iota$, i.e. 
\[
    E(\Gamma) = \Big\{  \{h,h'\} \colon h,h' \in H(\Gamma) \text{ and } \iota(h) = h'\Big\}.
\]
A \emph{loop} is an edge $e = \{h,h'\} \in E(\Gamma)$ such that $r(h) = r(h')$. The \emph{tangent space $T_v \Gamma$} of $\Gamma$ at a vertex $v \in V(\Gamma)$ is defined as $T_v \Gamma \coloneqq r^{-1}(\{v\}) \setminus \{v\}$ and the cardinality of $T_v \Gamma$ is denoted by $\val(v)$ and is called the \emph{valence} of $v$. A \emph{subgraph} of $\Gamma = (X, \iota, r)$ is a graph $\Gamma' = (X', \iota', r')$ such that 
\[
    X' \subseteq X, \quad \iota' = \iota_{|X'} \quad \text{ and } \quad  r' = r_{|X'}.
\]
Unless explicitly mentioned, all graphs considered in this paper will be finite and connected.

A \emph{weighted graph} $( {\Gamma}, g)$ is a graph $ {\Gamma}$ together with a function $g \colon V( {\Gamma}) \to \Z_{\geq 0}$. 
We say that the weighted graph $(\Gamma,g)$ is \emph{stable} at a vertex $v \in V(\Gamma)$ if we have:
\begin{equation*} \label{eq:stability}
	2 g(v) - 2 + \val(v) > 0.
\end{equation*}
A weighted graph is called \emph{stable} if it is stable at every vertex. This is equivalent to saying that vertices of weight $0$ must have valence at least $3$ and $\Gamma$ is not an isolated vertex of weight 1.
In the context of tropical geometry, the weight $g(v)$ of a vertex $v$ is also referred to as the \emph{genus} of $v$ and usually $g$ is omitted from the notation. That is we will denote the weighted graph $(\Gamma, g)$ simply by $\Gamma$ again. Moreover, the \emph{genus} of a weighted graph $\Gamma$ is given by:
    \[
        \big|E(\Gamma)\big| - \big|V(\Gamma)\big| + 1+\sum_{v\in V(\Gamma)}g(v),
    \]
and by abuse of notation this quantity is denoted by $g$ as well. Since the vertex genus and the genus of the graph are clearly related, we believe that confusion is unlikely.

\subsection{Maps of graphs}

A map of graphs from ${\Gamma_1} = (X_1, \iota_1, r_1)$ to ${\Gamma_2} = (X_2, \iota_2, r_2)$ is a map $f \colon X_1 \to X_2$ which is compatible with the root and involution maps in the sense that $f\circ r_1 = r_2 \circ f$ and $f\circ \iota_1 = \iota_2 \circ f$. We say that $f$ is \emph{finite} if the inverse image of any vertex $v_2 \in V(\Gamma_2)$ consists only of vertices in $V(\Gamma_1)$. We say that $f$ is an \emph{isomorphism} if it has an inverse map.

Stable graphs of fixed genus $g$ form a category $\bbGamma_g$ whose morphism are defined as follows.

\begin{definition}\label{def:MorphismsOfGenusgGraphs}
    Let $g \geq 1$ be an integer. A $\bbGamma_g$-\emph{morphism} between two genus $g$ weighted stable graphs $(\Gamma_1, g_1)$ and $(\Gamma_2, g_2)$ is a (not necessarily finite) map of graphs $f \colon \Gamma_1 \to \Gamma_2$ such that: 
    \begin{enumerate}
        \item For each $h\in H(\Gamma_2)$, the set $f^{-1}(h)$ consists of only one element in $H(\Gamma_1)$.
        
        \item For each $v \in V(\Gamma_2)$, the inverse image $f^{-1}(v)$, which is a subgraph of $ {\Gamma_1}$, is connected and has genus equal to $g_2(v)$.    
    \end{enumerate}
    We say that the morphism $f$ is an isomorphism if it has an inverse morphism. The set of automorphisms of a weighted graph $(\Gamma,g)$ is denoted by $\Aut(\Gamma)$.
\end{definition}

The purpose of the category $\bbGamma_g$ is to define the moduli space of genus $g$ tropical curves in a concise way, see \cref{subsec:SymmDeltaComplexes}. 

\begin{example}
    Let $(\Gamma, g)$ be a weighted stable graph and $e = \{h,h'\} \in E(\Gamma)$ an edge of $\Gamma$. We define $(\Gamma/e, g')$ as the weighted graph obtained from $\Gamma$ by contracting the edge $e$ and identifying its two end points $v_1,v_2 \in V(\Gamma)$ into a single vertex $v$. The vertex $v$ in $\Gamma/e$ is then assigned the weight $g'(v) = g(v_1) + g(v_2)$ if $e$ is not a loop and $g'(v) = g(v_1) + 1$ otherwise and all other vertices in $\Gamma/e$ keep their original weight as in $\Gamma$. Note that $(\Gamma/e, g')$ has the same genus as $(\Gamma, g)$ and is stable as well. The natural map $(\Gamma, g) \to (\Gamma/e, g')$ is a morphism of genus $g$ weighted stable graphs. We call $\Gamma$ an \emph{uncontraction} of $\Gamma/e$.
\end{example}

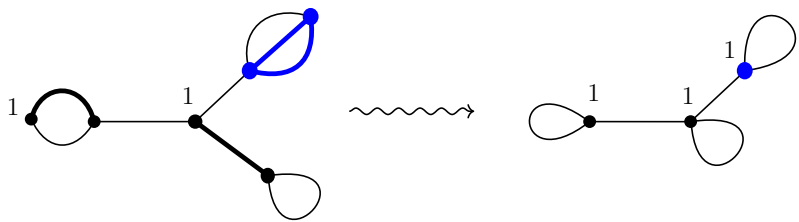
\begin{figure}[ht]
    \begin{center}
    \scalebox{0.8}{
        \tikzset{every picture/.style={line width=0.75pt}}
        
        \begin{tikzpicture}[x=0.75pt,y=0.75pt,yscale=-1,xscale=1]
        \draw    (154.92,187.54) -- (91.92,187.54) ;
        \draw [color={rgb, 255:red, 0; green, 0; blue, 0 }  ,draw opacity=1 ][line width=2.25]    (52.66,185.99) .. controls (58.66,161.99) and (85.66,161.99) .. (91.92,187.54) ;
        \draw    (52.66,185.99) .. controls (62.66,207.99) and (80.89,206.99) .. (91.92,187.54) ;
        \draw [color={rgb, 255:red, 0; green, 0; blue, 254 }  ,draw opacity=1 ][line width=2.25]    (226.92,121.54) -- (188.92,155.54) ;
        \draw    (188.92,155.54) .. controls (180.92,131.54) and (198.92,112.54) .. (226.92,121.54) ;
        \draw [color={rgb, 255:red, 0; green, 0; blue, 254 }  ,draw opacity=1 ][line width=2.25]    (188.92,155.54) .. controls (215.92,162.54) and (230.92,150.54) .. (226.92,121.54) ;
        \draw [color={rgb, 255:red, 0; green, 0; blue, 0 }  ,draw opacity=1 ][line width=2.25]    (200.43,221.54) -- (154.92,187.54) ;
        \draw  [fill={rgb, 255:red, 0; green, 0; blue, 0 }  ,fill opacity=1 ] (196.23,221.54) .. controls (196.23,219.08) and (198,217.09) .. (200.17,217.09) .. controls (202.35,217.09) and (204.11,219.08) .. (204.11,221.54) .. controls (204.11,223.99) and (202.35,225.98) .. (200.17,225.98) .. controls (198,225.98) and (196.23,223.99) .. (196.23,221.54) -- cycle ;
        \draw  [color={rgb, 255:red, 0; green, 0; blue, 255 }  ,draw opacity=1 ][fill={rgb, 255:red, 0; green, 0; blue, 254 }  ,fill opacity=1 ][line width=1.5]  (222.98,121.54) .. controls (222.98,119.08) and (224.74,117.09) .. (226.92,117.09) .. controls (229.09,117.09) and (230.86,119.08) .. (230.86,121.54) .. controls (230.86,123.99) and (229.09,125.98) .. (226.92,125.98) .. controls (224.74,125.98) and (222.98,123.99) .. (222.98,121.54) -- cycle ;
        \draw    (154.92,187.54) -- (188.92,155.54) ;
        \draw  [color={rgb, 255:red, 0; green, 0; blue, 255 }  ,draw opacity=1 ][fill={rgb, 255:red, 0; green, 0; blue, 254 }  ,fill opacity=1 ][line width=1.5]  (184.98,155.54) .. controls (184.98,153.08) and (186.74,151.09) .. (188.92,151.09) .. controls (191.09,151.09) and (192.86,153.08) .. (192.86,155.54) .. controls (192.86,157.99) and (191.09,159.98) .. (188.92,159.98) .. controls (186.74,159.98) and (184.98,157.99) .. (184.98,155.54) -- cycle ;
        
        \path[draw, ->, decorate, decoration ={snake, amplitude = 1.5}] (251,181.04) -- (329,181.04);    
        
        \draw    (460.42,187.58) -- (400.92,187.54) ;
        \draw  [fill={rgb, 255:red, 0; green, 0; blue, 0 }  ,fill opacity=1 ] (404.42,187.49) .. controls (404.44,189.43) and (402.89,191.01) .. (400.96,191.04) .. controls (399.03,191.06) and (397.44,189.51) .. (397.42,187.58) .. controls (397.39,185.65) and (398.94,184.06) .. (400.87,184.04) .. controls (402.81,184.01) and (404.39,185.56) .. (404.42,187.49) -- cycle ;
        \draw  [fill={rgb, 255:red, 0; green, 0; blue, 0 }  ,fill opacity=1 ] (467.42,187.49) .. controls (467.44,189.43) and (465.89,191.01) .. (463.96,191.04) .. controls (462.03,191.06) and (460.44,189.51) .. (460.42,187.58) .. controls (460.39,185.65) and (461.94,184.06) .. (463.87,184.04) .. controls (465.81,184.01) and (467.39,185.56) .. (467.42,187.49) -- cycle ;
        \draw    (463.92,187.54) -- (497.92,155.54) ;
        \draw    (400.92,187.54) .. controls (351.92,226.2) and (349.92,149.2) .. (400.66,187.54) ;
        \draw    (497.92,155.54) .. controls (569.92,145.2) and (494.92,83) .. (497.66,155.54) ;
        \draw    (200.43,221.54) .. controls (205.17,288.2) and (271.17,211.2) .. (200.17,221.54) ;
        \draw    (464.17,187.54) .. controls (468.92,254.2) and (534.92,177.2) .. (463.92,187.54) ;
        \draw  [fill={rgb, 255:red, 0; green, 0; blue, 0 }  ,fill opacity=1 ] (95.42,187.49) .. controls (95.44,189.43) and (93.89,191.01) .. (91.96,191.04) .. controls (90.03,191.06) and (88.44,189.51) .. (88.42,187.58) .. controls (88.39,185.65) and (89.94,184.06) .. (91.87,184.04) .. controls (93.81,184.01) and (95.39,185.56) .. (95.42,187.49) -- cycle ;
        \draw  [fill={rgb, 255:red, 0; green, 0; blue, 0 }  ,fill opacity=1 ] (56.16,185.95) .. controls (56.18,187.88) and (54.63,189.47) .. (52.7,189.49) .. controls (50.77,189.52) and (49.18,187.97) .. (49.16,186.04) .. controls (49.13,184.1) and (50.68,182.52) .. (52.61,182.49) .. controls (54.55,182.47) and (56.13,184.02) .. (56.16,185.95) -- cycle ;
        \draw  [fill={rgb, 255:red, 0; green, 0; blue, 0 }  ,fill opacity=1 ][line width=1.5]  (158.42,187.49) .. controls (158.44,189.43) and (156.89,191.01) .. (154.96,191.04) .. controls (153.03,191.06) and (151.44,189.51) .. (151.42,187.58) .. controls (151.39,185.65) and (152.94,184.06) .. (154.87,184.04) .. controls (156.81,184.01) and (158.39,185.56) .. (158.42,187.49) -- cycle ;
        \draw  [color={rgb, 255:red, 0; green, 0; blue, 255 }  ,draw opacity=1 ][fill={rgb, 255:red, 0; green, 0; blue, 254 }  ,fill opacity=1 ][line width=1.5]  (493.72,155.54) .. controls (493.72,153.08) and (495.49,151.09) .. (497.66,151.09) .. controls (499.84,151.09) and (501.6,153.08) .. (501.6,155.54) .. controls (501.6,157.99) and (499.84,159.98) .. (497.66,159.98) .. controls (495.49,159.98) and (493.72,157.99) .. (493.72,155.54) -- cycle ;
        
        \draw (145,164.4) node [anchor=north west][inner sep=0.75pt]    {$1$};
        \draw (36,171.4) node [anchor=north west][inner sep=0.75pt]    {$1$};
        \draw (457,164.94) node [anchor=north west][inner sep=0.75pt]    {$1$};
        \draw (398,162.94) node [anchor=north west][inner sep=0.75pt]    {$1$};
        \draw (483,135.94) node [anchor=north west][inner sep=0.75pt]    {$1$};
        \end{tikzpicture}
     }
    \end{center}
    \caption{A morphism of stable weighted graphs of genus $6$. The bolded edges are contracted. The blue vertex in the target graph has weight $1$ which is equal to the genus of the contracted subgraph colored in blue in the source graph.}
    \label{fig:ExampleOfEdgeContractions}
\end{figure}

\subsection{Abelian Galois covers of graphs}

In this subsection, we define abelian Galois covers which are the central objects of study of the present paper. The following definition is closely related to \cite[Definition~3.1]{BCK23}.

\begin{definition}\label{def:Gcover}
     Let $G$ be a finite abelian group.
     A \emph{(stable) $G$-cover} of graphs consists of a finite map of weighted graphs $\pi \colon (\tilde{\Gamma}, \tilde{g}) \to (\Gamma, g)$ together with the following data:
    \begin{enumerate}
        \item a $G$-action on $\tilde{\Gamma}$, and
        \item a map $\varphi \colon H(\Gamma^{\rm dil}) \to G$ where $\Gamma^{\rm dil}$ is the (possibly disconnected) subgraph of $\Gamma$ consisting of all elements $x \in X(\Gamma)$ with $|\pi^{-1}(x)| < |G|$,
    \end{enumerate}
    such that the following hold:
    \begin{enumerate}[label = (\alph*)]

        \item The weighted graphs $(\tilde \Gamma, \tilde g)$ and $(\Gamma, g)$ are stable.
        
        \item The map $\pi$ is canonically isomorphic to the quotient map $\tilde{\Gamma} \to \tilde{\Gamma}/G$ and the stabilizer $\Stab_G(\tilde{h})$ is a cyclic subgroup of $G$ for any $\tilde{h} \in H(\tilde{\Gamma})$.
        
        \item Local Riemann-Hurwitz condition: for any $v \in V(\Gamma)$ and $\tilde{v} \in \pi^{-1}(v)$ we have
        \begin{equation}\label{eq:local-RH}
            2 \tilde{g}(\tilde{v}) - 2 = d_{\pi}(\tilde{v}) \big(2 g(v) - 2 \big) + \sum_{\tilde{h} \in T_{\tilde v}(\tilde \Gamma)} \big(d_{\pi}(\tilde{h})  - 1 \big),
        \end{equation}
        with the degree map $d_{\pi} \colon X(\tilde{\Gamma}) \to \Z_{\geq 1}$ defined by $d_{\pi}(\tilde{x}) = |\Stab_G(\tilde{x})|$ for $\tilde x \in X(\tilde{\Gamma})$.
        
        \item  The map $\varphi$ satisfies:

                \begin{enumerate}[label = (\roman*)]
                
                    \item The fiber $\pi^{-1}(h)$ is isomorphic to $G/\langle \varphi(h) \rangle$ as left $G$-sets for any $h \in H(\Gamma)$.

                    \item For any $h \in H(\Gamma)$ we have $\varphi(\iota(h)) = - \varphi(h)$.

                    \item For any $v \in V(\Gamma)$ we have $\sum_{h \in T_v(\Gamma)} \varphi(h) = 0$.
                \end{enumerate}

    \end{enumerate}
\end{definition}

Every $G$-cover is harmonic in the sense of \cite{ABBR15}. More precisely, this means that for any $v \in V(\Gamma)$, $h \in T_v(\Gamma)$, and $\tilde{v} \in \pi^{-1}(v)$:
        \[
            d_\pi(\tilde{v}) = \sum_{\tilde{h} \in \pi^{-1}(h) \cap T_{\tilde v}(\tilde{\Gamma})}  d_\pi(\tilde{h}).
        \]
Note that harmonicity implies that for any vertex $v \in V(\Gamma)$ and any edge $e \in E(\Gamma)$, the following sums are independent of the choice of $v$ and $e$ and are equal to $|G|$ which is also the \emph{degree} of $\pi$:
\[
    \deg(\pi) \coloneqq \sum_{\tilde v \in \pi^{-1}(v)} d_{\pi}(\tilde{v}) = \sum_{\tilde e \in \pi^{-1}(e)} d_{\pi}(\tilde{e}) = |G|.
\]
We also note that, due to the local Riemann-Hurwitz condition \eqref{eq:local-RH}, there is a unique vertex weight function $\tilde{g}$ on $\tilde{\Gamma}$ that makes $\pi$ a $G$-cover.
Moreover, the respective genera $\tilde{g}$ and $g$ of $\tilde{\Gamma}$ and $\Gamma$ satisfy the \emph{global Riemann-Hurwitz} relation
\[
    \tilde{g} = |G|(g - 1) + 1.
\]
The map $\varphi$ is called \emph{dilation flow}\footnote{In \cite[Definition 3.1]{BCK23} dilation flows were called \emph{monodromy marking}.} and the (possibly disconnected) graph $\Gamma^{\rm dil}$ is the \emph{dilation subgraph} of $\Gamma$. We say that elements in $X(\Gamma)$ which are in $X(\Gamma^{\rm dil})$ are \emph{dilated}, every other element is called \emph{free}. Moreover, we say that a $G$-cover is \emph{dilated} if $\Gamma^\dil \neq \emptyset $ and \emph{free} otherwise. The integers $d_\pi(\tilde{v})$, $d_\pi(\tilde{e})$ and $d_{\pi}(\tilde{h})$ are sometimes called \emph{dilation factors}.

\medskip
We now discuss maps of $G$-covers.
\begin{definition} \label{def:mapOfCovers}
    Let $\pi_1 \colon \tilde{\Gamma}_1 \to \Gamma_1$ and $\pi_2 \colon \tilde{\Gamma}_2 \to \Gamma_2$ be two $G$-covers with dilation flows $\varphi_1, \varphi_2$ respectively. A \emph{map of $G$-covers} $(\tilde{f},f)$ consists of two maps of graphs $\tilde{f} \colon \tilde{\Gamma}_1 \to \tilde{\Gamma}_2$ and $f\colon \Gamma_1 \to \Gamma_2$ such 
    \[
        \begin{tikzcd}
         \tilde{\Gamma}_1 \arrow[r,"\tilde{f}"] \arrow[d, "\pi_1"'] & \tilde{\Gamma}_2 \arrow[d, "\pi_2"]\\
         \Gamma_1   \arrow[r,"f"]                    & \Gamma_2
        \end{tikzcd}
    \]
    commutes and: 
    \begin{enumerate}
        \item The map $\tilde{f}$ is $G$-equivariant.
        \item The map $f$ is compatible with the dilation flows i.e. $\varphi_2 \circ f = \varphi_1$.
    \end{enumerate}
\end{definition}

\begin{remark}
    There are several possible ways to classify $G$-covers, depending on whether or not we want to think of the target graph as ``rigid" and on whether we take automorphisms of the cover into account or not. One also encounters the same considerations when classifying covers of algebraic curves, see for example \cite{Fried77,Mochizuki95,Debes97}. In this sense, we made a choice in \cref{def:mapOfCovers} which reflects the space we want to study in the end.
\end{remark}

Similar to $\bbGamma_{g}$, we define a category $\bbGamma_{g,G}$ whose objects are $G$-covers of genus $g$ curves and, following \cite{BR11}, the morphisms of $\bbGamma_{g,G}$ are defined as follows:

\begin{definition}\label{def:MorphismsAreSquares}
    Let $\pi_1 \colon \tilde{\Gamma}_1 \to \Gamma_1$ and $\pi_2 \colon \tilde{\Gamma}_2 \to \Gamma_2$ be two $G$-covers of genus $g$ curves. A \emph{$\bbGamma_{g,G}$-morphism of $G$-covers} is a map of $G$-covers $(\tilde{f}, f)$ where $\tilde f$ and $f$ are both morphisms of graphs in the sense of \cref{def:MorphismsOfGenusgGraphs}, i.e. they are compositions of isomorphisms and edge contractions.
\end{definition}

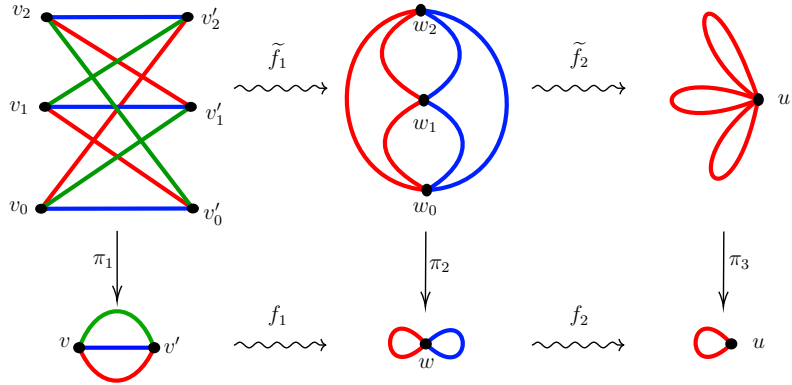
\begin{figure}
    \centering
    \scalebox{0.7}{
        \tikzset{every picture/.style={line width=0.75pt}} 
        \begin{tikzpicture}[x=0.75pt,y=0.75pt,yscale=-1,xscale=1]
        
        \draw [color={rgb, 255:red, 8; green, 167; blue, 0 }  ,draw opacity=1 ][line width=2.25]    (61.61,326.32) .. controls (74.87,291.12) and (102.62,291.66) .. (115.06,326.32) ;
        \draw [color={rgb, 255:red, 255; green, 0; blue, 0 }  ,draw opacity=1 ][line width=2.25]    (61.61,326.32) .. controls (74.87,357.43) and (100.71,358.35) .. (115.06,326.32) ;
        \draw [color={rgb, 255:red, 0; green, 13; blue, 255 }  ,draw opacity=1 ][line width=2.25]    (61.61,326.32) -- (115.06,326.32) ;
        
        \draw [color={rgb, 255:red, 255; green, 0; blue, 0 }  ,draw opacity=1 ][line width=2.25]    (38.5,88.83) -- (141.17,153.23) ;
        \draw [color={rgb, 255:red, 255; green, 0; blue, 0 }  ,draw opacity=1 ][line width=2.25]    (37.25,153.22) -- (142.64,226.28) ;
        \draw [color={rgb, 255:red, 255; green, 0; blue, 0 }  ,draw opacity=1 ][line width=2.25]    (34.51,226.46) -- (139.7,88.05) ;
        \draw [color={rgb, 255:red, 0; green, 32; blue, 255 }  ,draw opacity=1 ][line width=2.25]    (38.5,88.83) -- (138.87,88.61) ;
        \draw [color={rgb, 255:red, 0; green, 32; blue, 255 }  ,draw opacity=1 ][line width=2.25]    (34.51,226.46) -- (142.64,226.28) ;
        \draw [color={rgb, 255:red, 0; green, 32; blue, 255 }  ,draw opacity=1 ][line width=2.25]    (37.25,153.22) -- (141.17,153.23) ;
        \draw [color={rgb, 255:red, 1; green, 152; blue, 4 }  ,draw opacity=1 ][line width=2.25]    (37.25,153.22) -- (139.7,88.05) ;
        \draw [color={rgb, 255:red, 1; green, 152; blue, 4 }  ,draw opacity=1 ][line width=2.25]    (38.5,88.83) -- (142.64,226.28) ;
        \draw [color={rgb, 255:red, 1; green, 152; blue, 4 }  ,draw opacity=1 ][line width=2.25]    (34.51,226.46) -- (141.17,153.23) ;
        \draw  [fill={rgb, 255:red, 0; green, 0; blue, 0 }  ,fill opacity=1 ] (57.8,326.32) .. controls (57.8,324.38) and (59.5,322.8) .. (61.61,322.8) .. controls (63.72,322.8) and (65.42,324.38) .. (65.42,326.32) .. controls (65.42,328.26) and (63.72,329.83) .. (61.61,329.83) .. controls (59.5,329.83) and (57.8,328.26) .. (57.8,326.32) -- cycle ;
        \draw  [fill={rgb, 255:red, 0; green, 0; blue, 0 }  ,fill opacity=1 ] (111.25,326.32) .. controls (111.25,324.38) and (112.96,322.81) .. (115.06,322.81) .. controls (117.17,322.81) and (118.87,324.38) .. (118.87,326.32) .. controls (118.87,328.26) and (117.17,329.84) .. (115.06,329.84) .. controls (112.96,329.84) and (111.25,328.26) .. (111.25,326.32) -- cycle ;
        \draw  [fill={rgb, 255:red, 0; green, 0; blue, 0 }  ,fill opacity=1 ] (138.85,226.28) .. controls (138.85,224.59) and (140.55,223.22) .. (142.64,223.22) .. controls (144.73,223.22) and (146.43,224.59) .. (146.43,226.28) .. controls (146.43,227.97) and (144.73,229.35) .. (142.64,229.35) .. controls (140.55,229.35) and (138.85,227.97) .. (138.85,226.28) -- cycle ;
        \draw  [fill={rgb, 255:red, 0; green, 0; blue, 0 }  ,fill opacity=1 ] (137.38,153.23) .. controls (137.38,151.54) and (139.08,150.17) .. (141.17,150.17) .. controls (143.27,150.17) and (144.96,151.54) .. (144.96,153.23) .. controls (144.96,154.92) and (143.27,156.3) .. (141.17,156.3) .. controls (139.08,156.3) and (137.38,154.92) .. (137.38,153.23) -- cycle ;
        \draw  [fill={rgb, 255:red, 0; green, 0; blue, 0 }  ,fill opacity=1 ] (135.08,88.61) .. controls (135.08,86.92) and (136.78,85.55) .. (138.87,85.55) .. controls (140.96,85.55) and (142.66,86.92) .. (142.66,88.61) .. controls (142.66,90.3) and (140.96,91.67) .. (138.87,91.67) .. controls (136.78,91.67) and (135.08,90.3) .. (135.08,88.61) -- cycle ;
        \draw  [fill={rgb, 255:red, 0; green, 0; blue, 0 }  ,fill opacity=1 ] (30.72,226.46) .. controls (30.72,224.76) and (32.42,223.39) .. (34.51,223.39) .. controls (36.6,223.39) and (38.3,224.76) .. (38.3,226.46) .. controls (38.3,228.15) and (36.6,229.52) .. (34.51,229.52) .. controls (32.42,229.52) and (30.72,228.15) .. (30.72,226.46) -- cycle ;
        \draw  [fill={rgb, 255:red, 0; green, 0; blue, 0 }  ,fill opacity=1 ] (33.46,153.22) .. controls (33.46,151.52) and (35.16,150.15) .. (37.25,150.15) .. controls (39.35,150.15) and (41.04,151.52) .. (41.04,153.22) .. controls (41.04,154.91) and (39.35,156.28) .. (37.25,156.28) .. controls (35.16,156.28) and (33.46,154.91) .. (33.46,153.22) -- cycle ;
        \draw  [fill={rgb, 255:red, 0; green, 0; blue, 0 }  ,fill opacity=1 ] (34.71,88.83) .. controls (34.71,87.13) and (36.41,85.76) .. (38.5,85.76) .. controls (40.59,85.76) and (42.29,87.13) .. (42.29,88.83) .. controls (42.29,90.52) and (40.59,91.89) .. (38.5,91.89) .. controls (36.41,91.89) and (34.71,90.52) .. (34.71,88.83) -- cycle ;
        \draw    (88.68,242.38) -- (88.68,290.73) ;
        \draw [shift={(88.68,292.73)}, rotate = 270] [color={rgb, 255:red, 0; green, 0; blue, 0 }  ][line width=0.75]    (10.93,-3.29) .. controls (6.95,-1.4) and (3.31,-0.3) .. (0,0) .. controls (3.31,0.3) and (6.95,1.4) .. (10.93,3.29)   ;
        \draw    (308.68,242.84) -- (308.39,295.05) ;
        \draw [shift={(308.38,297.05)}, rotate = 270.32] [color={rgb, 255:red, 0; green, 0; blue, 0 }  ][line width=0.75]    (10.93,-3.29) .. controls (6.95,-1.4) and (3.31,-0.3) .. (0,0) .. controls (3.31,0.3) and (6.95,1.4) .. (10.93,3.29)   ;
        
        \path[draw, ->, decorate, decoration ={snake, amplitude = 1.5}] (171.87,139.86) -- (238.62,139.86);

        \path[draw, ->, decorate, decoration ={snake, amplitude = 1.5}] (171.87,324.39) -- (237.93,324.39);    
        
        \draw [color={rgb, 255:red, 255; green, 0; blue, 0 }  ,draw opacity=1 ][line width=2.25]    (308.81,323.59) .. controls (273.62,357.19) and (275.68,286.8) .. (308.64,323.59) ;
        \draw [color={rgb, 255:red, 0; green, 32; blue, 255 }  ,draw opacity=1 ][line width=2.25]    (308.58,323.59) .. controls (344.2,289.53) and (343.51,357.19) .. (308.64,323.59) ;
        \draw [color={rgb, 255:red, 255; green, 0; blue, 0 }  ,draw opacity=1 ][line width=2.25]    (305.12,83.83) .. controls (267.43,102.96) and (268.12,126.88) .. (307.25,148.45) ;
        \draw [color={rgb, 255:red, 255; green, 0; blue, 0 }  ,draw opacity=1 ][line width=2.25]    (307.25,148.45) .. controls (269.56,167.58) and (270.25,191.5) .. (309.37,213.07) ;
        \draw [color={rgb, 255:red, 255; green, 0; blue, 0 }  ,draw opacity=1 ][line width=2.25]    (305.12,83.83) .. controls (233.75,90.18) and (234.43,207.73) .. (309.37,213.07) ;
        \draw [color={rgb, 255:red, 0; green, 32; blue, 255 }  ,draw opacity=1 ][line width=2.25]    (305.12,83.83) .. controls (385,90.18) and (387.06,209.09) .. (309.37,213.07) ;
        \draw [color={rgb, 255:red, 0; green, 32; blue, 255 }  ,draw opacity=1 ][line width=2.25]    (305.12,83.83) .. controls (343.75,102.96) and (339.62,128.93) .. (307.25,148.45) ;
        \draw [color={rgb, 255:red, 0; green, 32; blue, 255 }  ,draw opacity=1 ][line width=2.25]    (307.25,148.45) .. controls (345.87,167.58) and (341.75,193.55) .. (309.37,213.07) ;
        \draw  [fill={rgb, 255:red, 0; green, 0; blue, 0 }  ,fill opacity=1 ] (301.84,83.83) .. controls (301.84,81.57) and (303.31,79.75) .. (305.12,79.75) .. controls (306.93,79.75) and (308.41,81.57) .. (308.41,83.83) .. controls (308.41,86.08) and (306.93,87.9) .. (305.12,87.9) .. controls (303.31,87.9) and (301.84,86.08) .. (301.84,83.83) -- cycle ;
        \draw  [fill={rgb, 255:red, 0; green, 0; blue, 0 }  ,fill opacity=1 ] (303.96,148.45) .. controls (303.96,146.2) and (305.43,144.37) .. (307.25,144.37) .. controls (309.06,144.37) and (310.53,146.2) .. (310.53,148.45) .. controls (310.53,150.7) and (309.06,152.53) .. (307.25,152.53) .. controls (305.43,152.53) and (303.96,150.7) .. (303.96,148.45) -- cycle ;
        \draw  [fill={rgb, 255:red, 0; green, 0; blue, 0 }  ,fill opacity=1 ] (306.09,213.07) .. controls (306.09,210.82) and (307.56,208.99) .. (309.37,208.99) .. controls (311.19,208.99) and (312.66,210.82) .. (312.66,213.07) .. controls (312.66,215.32) and (311.19,217.15) .. (309.37,217.15) .. controls (307.56,217.15) and (306.09,215.32) .. (306.09,213.07) -- cycle ;
        \draw  [fill={rgb, 255:red, 0; green, 0; blue, 0 }  ,fill opacity=1 ] (304.82,323.59) .. controls (304.82,321.65) and (306.53,320.07) .. (308.64,320.07) .. controls (310.74,320.07) and (312.45,321.65) .. (312.45,323.59) .. controls (312.45,325.53) and (310.74,327.1) .. (308.64,327.1) .. controls (306.53,327.1) and (304.82,325.53) .. (304.82,323.59) -- cycle ;
        \draw    (521.08,242.82) -- (520.79,295.03) ;
        \draw [shift={(520.77,297.03)}, rotate = 270.32] [color={rgb, 255:red, 0; green, 0; blue, 0 }  ][line width=0.75]    (10.93,-3.29) .. controls (6.95,-1.4) and (3.31,-0.3) .. (0,0) .. controls (3.31,0.3) and (6.95,1.4) .. (10.93,3.29)   ;
        \path[draw, ->, decorate, decoration ={snake, amplitude = 1.5}]  (384.26,139.85) -- (451.01,139.85);
        \path[draw, ->, decorate, decoration ={snake, amplitude = 1.5}] (384.26,324.37) -- (450.33,324.37);    
        \draw [color={rgb, 255:red, 255; green, 0; blue, 0 }  ,draw opacity=1 ][line width=2.25]    (526.58,323.57) .. controls (491.39,357.17) and (493.45,286.78) .. (526.4,323.57) ;
        \draw  [fill={rgb, 255:red, 0; green, 0; blue, 0 }  ,fill opacity=1 ] (522.59,323.57) .. controls (522.59,321.63) and (524.3,320.06) .. (526.4,320.06) .. controls (528.51,320.06) and (530.21,321.63) .. (530.21,323.57) .. controls (530.21,325.51) and (528.51,327.08) .. (526.4,327.08) .. controls (524.3,327.08) and (522.59,325.51) .. (522.59,323.57) -- cycle ;
        \draw [color={rgb, 255:red, 255; green, 0; blue, 0 }  ,draw opacity=1 ][line width=2.25]    (545.97,148.02) .. controls (463.71,180.87) and (462.81,116.91) .. (545.79,148.02) ;
        \draw [color={rgb, 255:red, 255; green, 0; blue, 0 }  ,draw opacity=1 ][line width=2.25]    (545.79,148.02) .. controls (482.87,103.89) and (507.04,34.52) .. (545.62,148.02) ;
        \draw [color={rgb, 255:red, 255; green, 0; blue, 0 }  ,draw opacity=1 ][line width=2.25]    (545.62,148.02) .. controls (479.82,200.69) and (527.27,243.94) .. (545.44,148.02) ;
        \draw  [fill={rgb, 255:red, 0; green, 0; blue, 0 }  ,fill opacity=1 ] (542.51,148.02) .. controls (542.51,145.77) and (543.98,143.95) .. (545.79,143.95) .. controls (547.61,143.95) and (549.08,145.77) .. (549.08,148.02) .. controls (549.08,150.28) and (547.61,152.1) .. (545.79,152.1) .. controls (543.98,152.1) and (542.51,150.28) .. (542.51,148.02) -- cycle ;
        \draw (45.87,317.7) node [anchor=north west][inner sep=0.75pt]    {$v$};
        \draw (120.34,317.01) node [anchor=north west][inner sep=0.75pt]    {$v'$};
        \draw (148.72,220.7) node [anchor=north west][inner sep=0.75pt]    {$v_{0} '$};
        \draw (147,80.4) node [anchor=north west][inner sep=0.75pt]    {$v_{2} '$};
        \draw (10.17,219.7) node [anchor=north west][inner sep=0.75pt]    {$v_{0}$};
        \draw (10.53,147.67) node [anchor=north west][inner sep=0.75pt]    {$v_{1}$};
        \draw (12.26,78.58) node [anchor=north west][inner sep=0.75pt]    {$v_{2}$};
        \draw (150.08,149.03) node [anchor=north west][inner sep=0.75pt]    {$v_{1} '$};
        \draw (69.9,257.92) node [anchor=north west][inner sep=0.75pt]    {$\pi _{1}$};
        \draw (302,332.4) node [anchor=north west][inner sep=0.75pt]    {$w$};
        \draw (298,220.4) node [anchor=north west][inner sep=0.75pt]    {$w_{0}$};
        \draw (298,90.4) node [anchor=north west][inner sep=0.75pt]    {$w_{2}$};
        \draw (298,160.4) node [anchor=north west][inner sep=0.75pt]    {$w_{1}$};
        \draw (309.54,261.34) node [anchor=north west][inner sep=0.75pt]    {$\pi _{2}$};
        \draw (194.68,105.96) node [anchor=north west][inner sep=0.75pt]    {$\tilde{f}_{1}$};
        \draw (194.68,294.79) node [anchor=north west][inner sep=0.75pt]    {$f_{1}$};
        \draw (540,318) node [anchor=north west][inner sep=0.75pt]    {$u$};
        \draw (523.14,258.76) node [anchor=north west][inner sep=0.75pt]    {$\pi _{3}$};
        \draw (409.55,105.06) node [anchor=north west][inner sep=0.75pt]    {$\tilde{f}_{2}$};
        \draw (409.55,294.79) node [anchor=north west][inner sep=0.75pt]    {$f_{2}$};
        \draw (558,143) node [anchor=north west][inner sep=0.75pt]    {$u$};
        
        \end{tikzpicture}
        
    }
    \caption{A sequence of morphisms of $\Z/3$-covers corresponding to first contracting the green edge in the target graph, and then the blue edge. The dilation graph is empty in the first two $3$-covers and it is a single vertex in the third picture, so there are no dilation flows.}
    \label{fig:G-contraction}
\end{figure}

\begin{example}
    An edge contraction on the target graph induces a morphism of $G$-covers as follows. Let $\pi \colon \tilde{\Gamma} \to \Gamma$ be a $G$-cover and $e$ an edge in $\Gamma$. Contracting the edge $e$ gives a morphism of graphs $f \colon \Gamma \to \Gamma/e$. Contracting the edges in $\tilde\Gamma$ mapping to $e$ gives a morphism or graphs $\tilde{f} \colon \tilde \Gamma \to \tilde \Gamma /\pi^{-1}(e)$. It is elementary to verify that this yields a $G$-cover $\tilde \Gamma / \pi^{-1}(e) \to \Gamma/e$, which we denote by $\pi/e$. In addition, $(\tilde f, f)$ is a morphism of $G$-covers. We call $(\tilde f, f)$ an edge contraction and say that the $G$-cover $\pi$ is an \emph{uncontraction} of the $G$-cover $\pi/e$. See \cref{fig:G-contraction}.
\end{example}

\subsection{Tropical \texorpdfstring{$G$}{}-covers}

     An abstract \emph{tropical curve} $C = (\Gamma, g, \ell)$ of genus $g$ is a weighted stable graph $(\Gamma, g)$ of genus $g$ and a length function $\ell \colon E(\Gamma) \to \R_{>0}$. The weighted graph $(\Gamma, g)$ is called the \emph{combinatorial type} of the tropical curve $C$. The volume $\vol(C)$ of $C$ is
    \[
        \vol(C) = \sum_{e \in E(\Gamma)} \ell(e).
    \]

    \begin{definition} \label{def:tropicalGcover}
        A \emph{tropical $G$-cover} is a $G$-cover of graphs $\pi \colon \tilde \Gamma \to \Gamma$, together with a tropical curve structure -- that is,  edge-lengths -- on the target $\Gamma$ and compatible edge-lengths on $\tilde \Gamma$ i.e.
        \[
            \ell \big(\pi(\tilde e) \big) = d_\pi(\tilde e) \ \tilde \ell(\tilde e) \qquad \text{for all } \tilde e \in E(\tilde\Gamma).
        \]
    \end{definition}

    \noindent Definition~\ref{def:tropicalGcover} justifies the terminology \enquote{dilation factor}.
    
    We emphasize that even without the dilation flow data, tropical $G$-covers carry more information than tropical Hurwitz covers in the sense of \cite{CMR16}. This is illustrated in the following example.

\begin{example}\label{ex:spiralCovers}
    \cref{fig:exampleOfGCovering} depicts two different $\Z/5$-covers of the same weighted graph of genus $2$. The covers map every vertex in the source curves to the unique vertex in the target and the edges are all sent to the one loop in the target. The group $\Z/5$ acts on both $\tilde \Gamma_1$ and $\tilde \Gamma_2$ by permuting the vertices cyclically as the permutation $(0,1,2,3,4)$. Note that the two source graphs are isomorphic as graphs (and as topological covers) but the $\Z/5$-covers $\pi_1$ and $\pi_2$ are \emph{not} isomorphic. This will be discussed further in \cref{rem:OnSpirals}.
    
    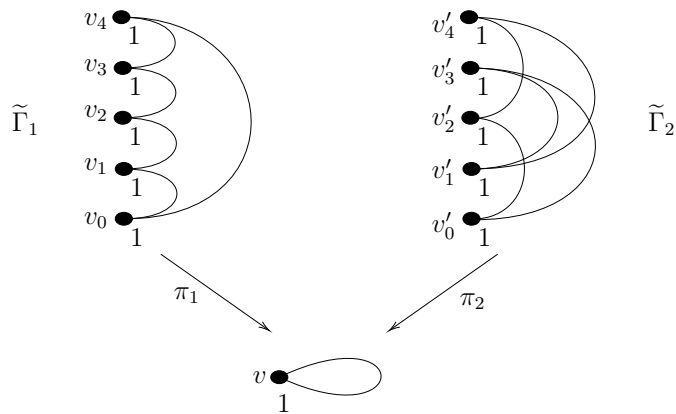
\begin{figure}[ht]
        \centering
            \scalebox{0.5}{      
            \tikzset{every picture/.style={line width=0.75pt}} 
            \begin{tikzpicture}[x=0.75pt,y=0.75pt,yscale=-1,xscale=1]
             
            \draw  [fill={rgb, 255:red, 0; green, 0; blue, 0 }  ,fill opacity=1 ] (58.46,79.4) .. controls (58.46,75.83) and (62.28,72.93) .. (66.99,72.93) .. controls (71.7,72.93) and (75.51,75.83) .. (75.51,79.4) .. controls (75.51,82.98) and (71.7,85.87) .. (66.99,85.87) .. controls (62.28,85.87) and (58.46,82.98) .. (58.46,79.4) -- cycle ;
             
            \draw  [fill={rgb, 255:red, 0; green, 0; blue, 0 }  ,fill opacity=1 ] (59.78,180.59) .. controls (59.78,177.02) and (63.6,174.12) .. (68.31,174.12) .. controls (73.02,174.12) and (76.83,177.02) .. (76.83,180.59) .. controls (76.83,184.16) and (73.02,187.06) .. (68.31,187.06) .. controls (63.6,187.06) and (59.78,184.16) .. (59.78,180.59) -- cycle ;
             
            \draw  [fill={rgb, 255:red, 0; green, 0; blue, 0 }  ,fill opacity=1 ] (59.78,130.5) .. controls (59.78,126.92) and (63.6,124.03) .. (68.31,124.03) .. controls (73.02,124.03) and (76.83,126.92) .. (76.83,130.5) .. controls (76.83,134.07) and (73.02,136.97) .. (68.31,136.97) .. controls (63.6,136.97) and (59.78,134.07) .. (59.78,130.5) -- cycle ;
             
            \draw  [fill={rgb, 255:red, 0; green, 0; blue, 0 }  ,fill opacity=1 ] (60.87,282.2) .. controls (60.87,278.63) and (64.68,275.73) .. (69.39,275.73) .. controls (74.1,275.73) and (77.92,278.63) .. (77.92,282.2) .. controls (77.92,285.77) and (74.1,288.67) .. (69.39,288.67) .. controls (64.68,288.67) and (60.87,285.77) .. (60.87,282.2) -- cycle ;
             
            \draw  [fill={rgb, 255:red, 0; green, 0; blue, 0 }  ,fill opacity=1 ] (60.87,232.11) .. controls (60.87,228.54) and (64.68,225.64) .. (69.39,225.64) .. controls (74.1,225.64) and (77.92,228.54) .. (77.92,232.11) .. controls (77.92,235.68) and (74.1,238.58) .. (69.39,238.58) .. controls (64.68,238.58) and (60.87,235.68) .. (60.87,232.11) -- cycle ;
             
            \draw  [fill={rgb, 255:red, 0; green, 0; blue, 0 }  ,fill opacity=1 ] (405.39,79.4) .. controls (405.39,75.83) and (409.21,72.93) .. (413.91,72.93) .. controls (418.62,72.93) and (422.44,75.83) .. (422.44,79.4) .. controls (422.44,82.98) and (418.62,85.87) .. (413.91,85.87) .. controls (409.21,85.87) and (405.39,82.98) .. (405.39,79.4) -- cycle ;
              
            \draw  [fill={rgb, 255:red, 0; green, 0; blue, 0 }  ,fill opacity=1 ] (406.71,180.59) .. controls (406.71,177.02) and (410.53,174.12) .. (415.23,174.12) .. controls (419.94,174.12) and (423.76,177.02) .. (423.76,180.59) .. controls (423.76,184.16) and (419.94,187.06) .. (415.23,187.06) .. controls (410.53,187.06) and (406.71,184.16) .. (406.71,180.59) -- cycle ;
              
            \draw  [fill={rgb, 255:red, 0; green, 0; blue, 0 }  ,fill opacity=1 ] (406.71,130.5) .. controls (406.71,126.92) and (410.53,124.03) .. (415.23,124.03) .. controls (419.94,124.03) and (423.76,126.92) .. (423.76,130.5) .. controls (423.76,134.07) and (419.94,136.97) .. (415.23,136.97) .. controls (410.53,136.97) and (406.71,134.07) .. (406.71,130.5) -- cycle ;
            \draw  [fill={rgb, 255:red, 0; green, 0; blue, 0 }  ,fill opacity=1 ] (407.79,282.2) .. controls (407.79,278.63) and (411.61,275.73) .. (416.32,275.73) .. controls (421.03,275.73) and (424.84,278.63) .. (424.84,282.2) .. controls (424.84,285.77) and (421.03,288.67) .. (416.32,288.67) .. controls (411.61,288.67) and (407.79,285.77) .. (407.79,282.2) -- cycle ;
            \draw  [fill={rgb, 255:red, 0; green, 0; blue, 0 }  ,fill opacity=1 ] (407.79,232.11) .. controls (407.79,228.54) and (411.61,225.64) .. (416.32,225.64) .. controls (421.03,225.64) and (424.84,228.54) .. (424.84,232.11) .. controls (424.84,235.68) and (421.03,238.58) .. (416.32,238.58) .. controls (411.61,238.58) and (407.79,235.68) .. (407.79,232.11) -- cycle ;
            \draw    (66.99,79.4) .. controls (139.38,79.62) and (137.21,131.2) .. (68.31,130.5) ;
            \draw    (68.31,130.5) .. controls (140.7,130.72) and (137.21,181.29) .. (68.31,180.59) ;
            \draw    (68.31,180.59) .. controls (140.7,180.81) and (138.29,232.81) .. (69.39,232.11) ;
            \draw    (69.39,232.11) .. controls (141.78,232.33) and (138.29,282.9) .. (69.39,282.2) ;
            \draw    (66.99,79.4) .. controls (231.53,78.08) and (246.71,283.62) .. (69.39,282.2) ;
            \draw    (415.23,180.59) .. controls (487.62,180.81) and (486.05,283.16) .. (417.15,282.46) ;
            \draw    (413.91,79.4) .. controls (486.3,79.62) and (484.73,181.97) .. (415.83,181.27) ;
            \draw    (413.91,79.4) .. controls (569.78,80.39) and (591.46,232.04) .. (416.32,232.11) ;
            \draw    (415.23,130.5) .. controls (536.41,131.49) and (526.42,232.04) .. (416.32,232.11) ;
            \draw    (415.23,130.5) .. controls (578.45,131.2) and (583.87,284.39) .. (417.64,283.2) ;
            \draw  [fill={rgb, 255:red, 0; green, 0; blue, 0 }  ,fill opacity=1 ] (216.01,441.62) .. controls (216.01,438.05) and (219.83,435.15) .. (224.54,435.15) .. controls (229.25,435.15) and (233.06,438.05) .. (233.06,441.62) .. controls (233.06,445.2) and (229.25,448.09) .. (224.54,448.09) .. controls (219.83,448.09) and (216.01,445.2) .. (216.01,441.62) -- cycle ;
            \draw    (224.24,441.62) .. controls (361.12,376.76) and (358.95,504.55) .. (224.54,441.62) ;
            \draw    (106.34,317.87) -- (213.13,393.69) ;
            \draw [shift={(214.76,394.85)}, rotate = 215.38] [color={rgb, 255:red, 0; green, 0; blue, 0 }  ][line width=0.75]    (10.93,-3.29) .. controls (6.95,-1.4) and (3.31,-0.3) .. (0,0) .. controls (3.31,0.3) and (6.95,1.4) .. (10.93,3.29)   ;
            \draw    (442.43,317.87) -- (335.64,393.69) ;
            \draw [shift={(334.01,394.85)}, rotate = 324.62] [color={rgb, 255:red, 0; green, 0; blue, 0 }  ][line width=0.75]    (10.93,-3.29) .. controls (6.95,-1.4) and (3.31,-0.3) .. (0,0) .. controls (3.31,0.3) and (6.95,1.4) .. (10.93,3.29)   ;
            
            \draw (27.3,165.32) node [anchor=north west][inner sep=0.75pt] [xscale=1.7, yscale=1.7]   {$v_{2}$};
            \draw (27.3,119.13) node [anchor=north west][inner sep=0.75pt] [xscale=1.7, yscale=1.7]   {$v_{3}$};
            \draw (27.3,72.95) node [anchor=north west][inner sep=0.75pt] [xscale=1.7, yscale=1.7]   {$v_{4}$};
            \draw (27.3,273.09) node [anchor=north west][inner sep=0.75pt] [xscale=1.7, yscale=1.7]   {$v_{0}$};
            \draw (27.3,219.21) node [anchor=north west][inner sep=0.75pt]  [xscale=1.7, yscale=1.7]  {$v_{1}$};
            \draw (374.22,165.32) node [anchor=north west][inner sep=0.75pt] [xscale=1.7, yscale=1.7]   {$v'_{2}$};
            \draw (374.22,117.59) node [anchor=north west][inner sep=0.75pt] [xscale=1.7, yscale=1.7]   {$v'_{3}$};
            \draw (374.22,71.41) node [anchor=north west][inner sep=0.75pt] [xscale=1.7, yscale=1.7]   {$v'_{4}$};
            \draw (374.22,273.09) node [anchor=north west][inner sep=0.75pt] [xscale=1.7, yscale=1.7]   {$v'_{0}$};
            \draw (374.22,219.21) node [anchor=north west][inner sep=0.75pt] [xscale=1.7, yscale=1.7]   {$v'_{1}$};
            \draw (220,456.53) node [anchor=north west][inner sep=0.75pt] [xscale=1.7, yscale=1.7]   {$1$};
            \draw (195.66,433.44) node [anchor=north west][inner sep=0.75pt] [xscale=1.7, yscale=1.7]   {$v$};
            \draw (73.07,289.54) node [anchor=north west][inner sep=0.75pt] [xscale=1.7, yscale=1.7]   {$1$};
            \draw (73.07,239.45) node [anchor=north west][inner sep=0.75pt] [xscale=1.7, yscale=1.7]   {$1$};
            \draw (71.98,187.93) node [anchor=north west][inner sep=0.75pt] [xscale=1.7, yscale=1.7]   {$1$};
            \draw (71.98,137.84) node [anchor=north west][inner sep=0.75pt] [xscale=1.7, yscale=1.7]   {$1$};
            \draw (70.66,86.74) node [anchor=north west][inner sep=0.75pt] [xscale=1.7, yscale=1.7]   {$1$};
            \draw (419.99,289.54) node [anchor=north west][inner sep=0.75pt] [xscale=1.7, yscale=1.7]   {$1$};
            \draw (419.99,239.45) node [anchor=north west][inner sep=0.75pt] [xscale=1.7, yscale=1.7]   {$1$};
            \draw (418.91,187.93) node [anchor=north west][inner sep=0.75pt] [xscale=1.7, yscale=1.7]   {$1$};
            \draw (418.91,137.84) node [anchor=north west][inner sep=0.75pt] [xscale=1.7, yscale=1.7]   {$1$};
            \draw (417.59,86.74) node [anchor=north west][inner sep=0.75pt] [xscale=1.7, yscale=1.7]    {$1$};
            \draw (116.9,348.53) node [anchor=north west][inner sep=0.75pt] [xscale=1.7, yscale=1.7]    {$\pi _{1}$};
            \draw (402.03,354.69) node [anchor=north west][inner sep=0.75pt] [xscale=1.7, yscale=1.7]    {$\pi _{2}$};

            \draw (-45,164.4) node [anchor=north west][inner sep=0.75pt][xscale=1.7, yscale=1.7]    {$\tilde\Gamma_1$};
            \draw (590,164.4) node [anchor=north west][inner sep=0.75pt][xscale=1.7, yscale=1.7]    {$\tilde\Gamma_2$};
            
            \end{tikzpicture}
            }
        \caption{Two non-isomorphic $\Z/5$-covers.}
        \label{fig:exampleOfGCovering}
    \end{figure}
\end{example}

\subsection{The moduli space \texorpdfstring{$\Delta_{g,G}$}{} as a symmetric \texorpdfstring{$\Delta$}{}-complex}\label{subsec:SymmDeltaComplexes}

In \cite{ACP15} the authors introduced \emph{generalized cone complexes} as a combinatorial shadow of Deligne-Mumford stacks.
The link around the cone point captures important topological properties of the cone complex.
When the cones are simplicial, the link is a \emph{symmetric $\Delta$-complex}, which is a generalization of a $\Delta$-complex. Here, we briefly recall the definition of a symmetric $\Delta$-complexes and give a description of the moduli space of tropical $G$-covers $\Delta_{g,G}$ as a symmetric $\Delta$-complex. For a more detailed account we refer the reader to \cite[Section 3]{CGP1} and \cite{AllcockCoreyPayne}.

Let $I$ be the category whose objects are the finite sets
\[ 
    [n] \coloneqq \{0, \ldots, n\} 
\]
for $n \in \Z_{\geq -1}$ (with the convention that $[-1] \coloneqq \emptyset$) and whose morphisms are injections
\[
    [n] \longhookrightarrow [m], \quad \text{for } m \geq n.
\]
A \emph{symmetric $\Delta$-complex} is a functor $X \colon I^\op \to \mathrm{\bf Sets}$. For $n \geq 0$, we refer to $X_n \coloneqq X([n])$ as \emph{the set of $n$-simplices}.
Note that the functor $X$ sends any injection $\theta \colon [n] \hookrightarrow [m]$ to 
\[
    \theta^{\ast} \coloneqq X(\theta) \colon X_m \longrightarrow X_n,
\]
and in particular it induces an action of the symmetric group $S_{n+1}$ on $X_n$ by sending a permutation $\theta \colon [n] \to [n]$ to $\theta^\ast$. For any element $\sigma \in X_n$, denote by $\Aut(\sigma)$ the subgroup of $S_{n+1}$ stabilizing $\sigma$ and by $\Orb(\sigma) \coloneqq \{\theta^*(\sigma) \colon \theta \in S_{n+1}\}$ the symmetric orbit of $\sigma$. We say that $\tau \in X_n$ is a \emph{face} of $\sigma \in X_m$ and write $\tau \preccurlyeq \sigma$ if there is an injection $\theta \colon [n] \hookrightarrow [m]$ such that $\tau = \theta^\ast(\sigma)$. 

The \emph{geometric realization} of the symmetric $\Delta$-complex $X$ as a topological space is defined as
\begin{equation} \label{eq:geometricRealization}
    |X| \coloneqq \Big( \coprod_{n = 0}^\infty X_n \times \Delta^n \Big) / \sim \ , 
\end{equation}
where $(\theta^* \sigma,a) \sim (\sigma,\theta_* a)$ for any injection $\theta \colon[n] \hookrightarrow [m]$, simplex $\sigma \in X_m$, and any point $a$ in the standard $n$-simplex $\Delta^n$. Here,
\[
    \theta_\ast \colon \Delta^{n} \longrightarrow \Delta^{m},
\] 
is the map extending an element $a = (a_0, \dots, a_n) \in \Delta^n$ to $(b_0, \dots, b_m)$ with $b_{\theta(i)} = a_i$ and the remaining $b_j$'s for $j \in [m] \setminus \theta([n])$ are set to $0$. Alternatively to \eqref{eq:geometricRealization}, we can also describe $|X|$ as a union of cells, 
\[
    |X| = \coprod_n \coprod_{\sigma \in \underline{X}_n}(\Delta^n)^\circ/\Aut(\sigma) 
\] 
where $\underline{X}_n$ is a set of representatives of symmetric orbits. Note that in this description, each cell $(\Delta^n)^\circ/\Aut(\sigma)$ is homeomorphic to a quotient of the open simplex by the natural action of $\Aut(\sigma)$.

Now let $G$ be a finite abelian group. The moduli space $\Delta_{g, G}$ of tropical $G$-covers of genus $g$ tropical curves of volume $1$ can naturally be viewed as a symmetric $\Delta$-complex as follows. For $n\geq 0$, we define $\Delta_{g, G}([n])$ to be the set of equivalence classes of pairs $(\pi, \omega)$, where 
\begin{enumerate}
    \item $\pi \colon \tilde \Gamma \to \Gamma$ is a $G$-cover of graphs with $\big| E(\Gamma) \big| = n +1$ and genus of $\Gamma$ equal to $g$, and
    \item $\omega \colon E(\Gamma)\to [n]$ is an \emph{edge-labeling}, i.e. a bijection.
\end{enumerate}
Two pairs $(\pi_1:\tilde \Gamma_1 \to \Gamma_1, \omega_1)$ and $(\pi_2:\tilde \Gamma_2 \to \Gamma_2,\omega_2)$ are equivalent if there exists an isomorphism $(\tilde \varphi, \varphi)$ of $G$-covers such that
\begin{equation*}
\begin{tikzcd}
\tilde \Gamma_1 \arrow[rr,"\tilde \varphi"] \arrow[d,"\pi_1"']& & \tilde \Gamma_2 \arrow[d,"\pi_2"] \\
\Gamma_1 \arrow[rr, "\varphi"] \arrow[dr, "\omega_1"'] & & \Gamma_2 \arrow[dl, "\omega_2"] \\
& {[n]} &
\end{tikzcd}
\end{equation*} 
commutes. On the level of morphisms, an injection $\theta\colon [n] \to [m]$ sends an $m$-simplex $(\pi, \omega)$ to the pair $(\pi', \omega')$, where $\pi' \colon\tilde \Gamma'\to \Gamma'$ is obtained from $\pi \colon\tilde \Gamma \to \Gamma$ by contracting all edges $e \in E(\Gamma)$ with labels not in $\operatorname{Im}(\theta)$ and the corresponding edges in $\pi^{-1}(e)$. The function $\omega'$ is the induced edge-labeling, i.e. $\omega' \colon E(\Gamma') \to [n]$ is the unique bijection that preserves the total ordering on $E(\Gamma')$ given by $\omega$. Finally, we set $\Delta_{g,G}([-1]) = \emptyset$. 

\begin{figure}[ht]
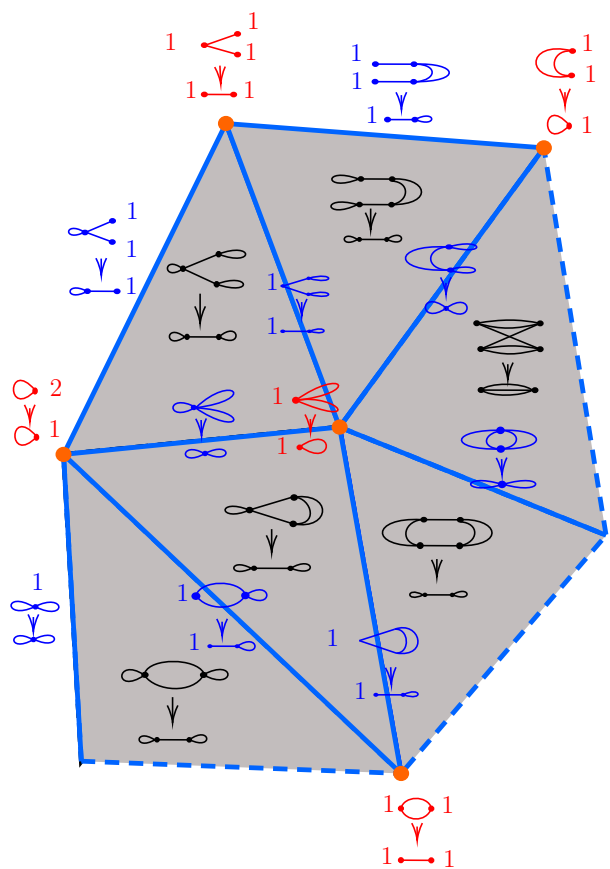

    \centering
    \scalebox{0.8}{    
        \tikzset{every picture/.style={line width=0.75pt}} 
        
        
    
    }
    \caption{The moduli space $\Delta_{2,2}$ as a $2$-dimensional symmetric $\Delta$-complex. The dashed edges correspond to a folding of a triangle due to automorphisms, hence there are no $2$-covers on dashed edges. In the picture we omitted any dilation flows, since in each $2$-cover there is a unique one.}
    \label{fig:Delta_2_2}
\end{figure}

For the rest of this article, we shall restrict our discussion to the case $G = \Z/p$ in which case we will simply refer to $\Z/p$-covers as $p$-covers. Note that, in a $p$-cover $\pi \colon \tilde \Gamma \to \Gamma$ every point in $X(\Gamma)$ has either 1 or $p$ preimages (with the dilation factor being $p$ or 1, respectively). Also, the condition that the stabilizer subgroup in $\Z/p$ of any point in $X(\tilde\Gamma)$ must be cyclic is void in this case (see \cref{def:Gcover}). Finally, the dilation flow is simply a nowhere-zero balanced flow. For ease of notation we denote the tropical moduli space $\Delta_{g, \Z/p}$ simply by $\MS$. The resulting space in $g = 2$ and $G = \Z/2$ is depicted in \cref{fig:Delta_2_2}.

In \cite[Section~4]{CGP2} the authors develop a contractibility criterion for subcomplexes of a symmetric $\Delta$-complex. This is the main technical ingredient we use in our proof of Theorem~\ref{thm:contractible_loci}. We summarize the essential results.

A \emph{property} $P$ is a subset $P \subseteq X_0$ of the set of vertices (i.e. $0$-simplices) of $X$. Given $P$, any vertex $v \in P$ is called a \emph{$P$-vertex}. The subcomplex $X_P \leq X$ is generated by the set of simplices with at least one $P$-vertex, i.e. $X_P$ is the closure of the star of $P$. The set of simplices of $X_P$ is denoted $P^*(X)$. We recall the following definitions from \cite[Section 4]{CGP2}:
\begin{enumerate}
	\item We say that $\tau$ is a \emph{co-$P$ face} of $\sigma$ and write $\tau \preccurlyeq_P \sigma$ if $\tau \preccurlyeq \sigma$ and the vertices of $\sigma$ which do not belong to $\tau$ are all in $P$.
	
	\item We say that $\tau$ has a \emph{co-$P$ maximal face} if $\{ \sigma \colon \tau \preccurlyeq_P \sigma\}$ has a maximal element with respect to face inclusion, $\sigma_0$ in some $X_m$, and which is also unique up to $S_{m+1}$-action. 
	
	\item A face inclusion $\tau \preccurlyeq \sigma$ is \emph{canonical} if for any two inclusions $\theta_1, \theta_2 \colon [n] \to [m]$ by which $\tau$ is a face of $\sigma$, there is a $\psi \in \Aut(\sigma)$ such that $\theta_1 = \psi \theta_2$.
	
	\item A collection $Y$ of simplices is \emph{co-$P$ saturated} if for all $\tau \in Y$ and all $\sigma$ with $\tau \preccurlyeq_P \sigma$ we have $\sigma \in Y$.
\end{enumerate}

\noindent
Now let $P$ and $Q$ be properties and consider the following conditions:

\begin{enumerate}[wide=12pt]
	\item[(CGP1)] The set of simplices $P^*(X)$ is co-$Q$-saturated.
	\item[(CGP2)] The simplices in $X \setminus X_P$ admit canonical co-$Q$ maximal faces.
	\item[(CGP3)] Every simplex whose vertices are all in $Q$, is in $X_P$.
\end{enumerate}

\noindent 
Under these conditions, by virtue of \cite[Proposition~4.11]{CGP2}, there is a strong deformation retraction $(X_P \cup X_Q) \searrow X_P$. Moreover, if we consider a sequence of properties $P_1, \ldots, P_N$ such that the conditions~(CGP1) to (CGP3) hold for $P = P_1 \cup \cdots \cup P_{i-1}$ and $Q = P_i$ for every $i = 2, \ldots, N$, then we can iterate the procedure and obtain $X_{P_1 \cup \cdots \cup P_N} \searrow X_{P_1}$. Finally, in the case where $P = \emptyset$ and $Q$ is a property, and the pair $(P,Q)$ satisfy conditions~(CGP1) and (CGP2), then $X_Q \searrow X_{Q,0}$ where $X_{Q,0} \leq X_Q$ is the subcomplex generated by simplices in $X$ that have no non-$Q$ vertices.

In the special case of $\Delta = \Delta_{g,G}$, we now reformulate the conditions  (CGP1) to (CGP3) in terms of $G$-covers. In this setting, simplices correspond to $G$-covers, face relations correspond to edge-contractions, and the \enquote{vertices} of a $G$-cover $\pi$ are the $G$-covers arising from $\pi$ by contracting all but one edge. Hence a property $P$ is a set of $G$-covers with a single edge. Given a $G$-cover $\pi$, an edge $e$ in $\pi$ is called a $P$-edge if contracting all the edges of $\pi$ except $e$ produces one of the covers in $P$. Using this dictionary, the conditions~(CGP1) to~(CGP3) translate as follows (we give the contrapositive of (CGP1) as this is what we use in our proofs below):

\begin{enumerate}[wide=12pt, leftmargin = 58pt]
    \item[{\crtcrossreflabel{(CGP1)}[condition1]}] If a $G$-cover neither contains a $P$-edge nor does it admit an uncontraction of a $P$-edge, then after contracting a $Q$-edge, the resulting $G$-cover continues to hold those properties.
    
    \item[{\crtcrossreflabel{(CGP2)}[condition2]}] Consider a $G$-cover $\pi$ which neither has a $P$-edge, nor does it allow for an uncontraction of a $P$-edge. Then there exists a unique (up to isomorphism) $\pi'$ such that $\pi = \pi' / \{e_1, \ldots, e_r\}$ for $Q$-edges $e_i$ in $\pi'$. Moreover, $r$ is maximal and every automorphism of $\pi$ lifts to an automorphism of~$\pi'$.
    
    \item[{\crtcrossreflabel{(CGP3)}[condition3]}] Consider a $G$-cover $\pi$ all of whose edges are $Q$-edges. Then $\pi$ has a $P$-edge or an uncontraction of a $P$-edge.
\end{enumerate}

\section{Contractible loci} \label{sec:contractible_loci}

In this section, we apply the tools from \cite{CGP2} to prove \cref{thm:contractible_loci}. To this end we work on the symmetric $\Delta$-complex $\MS$ that was defined in \cref{subsec:SymmDeltaComplexes}. Recall that in this space simplices correspond to $p$-covers of graphs $\pi \colon \tilde\Gamma \to \Gamma$ together with a labeling of the edges of $\Gamma$ and symmetric orbits of simplices correspond to isomorphism classes of $p$-covers without the edge-labeling. In this setting, two $p$-covers $\pi$ and $\rho$ are related by face inclusion $\pi \preccurlyeq \rho$ if and only if up to isomorphism $\pi$ is a contraction of $\rho$. This face inclusion is canonical if and only if for every automorphism $\theta \in \Aut(\pi)$ there exists an automorphism $\psi \in \Aut(\rho)$ such that $\psi$ induces $\theta$ under this contraction.
For the reader’s convenience, we have listed all vertices of $\MS$ in \cref{tab:vertices}. We will construct properties as defined in \cref{subsec:SymmDeltaComplexes} by giving collections of these vertices.

\begin{table}[ht]
    \centering
    \scalebox{0.9}{
    \begin{tabular}{|c|c|c|}
     \hline
      \makecell{Ring $R_i$ \\ for $1 \leq i \leq p-1$ }& Butterfly     &  \makecell{Spirals $S_a$ \\for $1 \leq a \leq \lceil \frac{p-1}{2} \rceil$, \\ see \cref{def:spirals}. }  \\
     \hline
        \begin{tikzpicture}            
            \draw (0.7,1.5) node [anchor=east] [xscale=0.8, yscale=0.8]  {$p(g-1)$};
            \fatvertex{1}{1.5};
            \selfloop[0.75]{1}{1.5};

            \draw[->] (1.2,1) -- (1.2,0.4);

            \draw (1.9,0.2) node [anchor=north west] [xscale=0.8, yscale=0.8]  {$i$};
            
            \draw (0.7,0) node [anchor=east] [xscale=0.8, yscale=0.8]  {$g-1$};
            \fatvertex{1}{0};
            \selfloop[0.75]{1}{0};
            \fatarrowup{(1.75,-0.1)}

        \end{tikzpicture} & 
        \begin{tikzpicture}
            \draw (0.7,1.5) node [anchor=east] [xscale=0.8, yscale=0.8]  {$p(g-2)+1$};
            \fatvertex{1}{1.5};
            \draw (1, 1.5) .. controls (2, 2.7) and (2, 1.7) .. (1, 1.5);
            
            \draw (1.5,2.01) node [anchor=north west] [xscale=1, yscale=1]  {$\vdots$};
            
            \draw (1, 1.5) .. controls (2, 0.3) and (2, 1.3) .. (1, 1.5);
                
            \draw[->] (1.2,1) -- (1.2,0.4);
            
            \draw (0.7,0) node [anchor=east] [xscale=0.8, yscale=0.8]  {$g-1$};
            \fatvertex{1}{0};
            \draw (1, 0) .. controls (2, 1) and (2, -1) .. (1, 0);
        \end{tikzpicture}
        &  
        \begin{tikzpicture}
            \draw (0.7, 1.5) node [anchor = east] [xscale=0.8, yscale=0.8]  {$g-1$};
            \fill (1, 1.5) circle (2pt);
            \draw (1, 2.0) node [anchor=center] [xscale=1, yscale=1]  {$\vdots$};
            \draw (0.7, 2.3) node [anchor = east] [xscale=0.8, yscale=0.8]  {$g-1$};
            \fill (1, 2.3) circle (2pt);
            \draw (1, 1.5) .. controls (2, 1.0) and (2, 2.8) .. (1, 2.3);
            \draw (1, 2.8) node [anchor=center] [xscale=1, yscale=1]  {$\vdots$};
            \draw (0.7, 3.1) node [anchor = east] [xscale=0.8, yscale=0.8]  {$g-1$};
            \fill (1, 3.1) circle (2pt);
            \draw (1, 2.3) .. controls (2, 1.8) and (2, 3.6) .. (1, 3.1);
            \begin{scope}
                \clip (1, 1.3) rectangle (2, 3.3);
                \draw (1, 3.1) .. controls (2, 2.6) and (2, 4.4) .. (1, 3.9);
                \draw (1, 1.5) .. controls (2, 2.0) and (2, 0.2) .. (1, 0.7);
            \end{scope}

            \draw [very thick,decorate,decoration = {calligraphic brace, amplitude=5pt}] (-0.4, 1.5) --  (-0.4, 2.3)
                node[pos=0.5,left=10pt]{$a$};

            \draw[->] (1.2,1) -- (1.2,0.4);
            
            \draw (0.7,0) node [anchor=east] [xscale=0.8, yscale=0.8]  {$g-1$};
            \fill (1, 0) circle (2pt);
            \draw (1, 0) .. controls (2, 1) and (2, -1) .. (1, 0);
         \end{tikzpicture}
        \\
        \hline
        \hline
        \makecell{Equivariant bridges $B_{h}$ \\ for $ 1\leq h \leq g-1$, \\ see \cref{def:coverB_h}.} & \multicolumn{2}{c|}{\makecell{Parallel bridges $P_h$ \\ for $1 \leq h \leq \lceil \frac{p-1}{2} \rceil $ }} 
        \\
        \hline
         \begin{tikzpicture}
            \draw (0.7,1.5) node [anchor=east] [xscale=0.8, yscale=0.8]  {$p(g-h-1) + 1$};
            \fatvertex{1}{1.5};
            \draw (1, 1.5) -- (2, 2);
            \fill (2, 2) circle (2pt);
            \draw (1, 1.5) -- (2, 1);
            \fill (2, 1) circle (2pt);

            \draw (2.3,2) node [anchor= west] [xscale=0.8, yscale=0.8]  {$h$};

            \draw (2.3,1) node [anchor= west] [xscale=0.8, yscale=0.8]  {$h$};
            \draw (1.8,2) node [anchor=north west] [xscale=1, yscale=1]  {$\vdots$};
            
            \draw[->] (1.2,1) -- (1.2,0.4);

            \draw (2.3,0) node [anchor= west] [xscale=0.8, yscale=0.8]  {$h$};
            \fill (2, 0) circle (2pt);
            \draw (0.7,0) node [anchor=east] [xscale=0.8, yscale=0.8]  {$g-h$};
            \fatvertex{1}{0};
            \draw (1, 0) --  (2, 0);
            
        \end{tikzpicture}

        &
        \multicolumn{2}{c|}{
        \begin{tikzpicture}
            \draw (0.7,1.5) node [anchor=east] [xscale=0.8, yscale=0.8]  {$p(g-h-1)+1$};
            \fatvertex{1}{1.5};
            \draw (1, 1.5) .. controls (1.2, 2) and (1.8, 2) .. (2, 1.5);
            \draw (1.5,1.6) node [anchor=center] [xscale=1, yscale=1]  {$\vdots$};
            \draw (1, 1.5) .. controls (1.2, 1) and (1.8, 1) .. (2, 1.5);
            \fatvertex{2}{1.5};
            \draw (2.3,1.5) node [anchor=west] [xscale=0.8, yscale=0.8]  {$p(h-1)+1$};
            
            \draw[->] (1.5,1) -- (1.5,0.4);

            \draw (2.3,0) node [anchor= west] [xscale=0.8, yscale=0.8]  {$h$};
            \fatvertex{2}{0};
            \draw (0.7,0) node [anchor=east] [xscale=0.8, yscale=0.8]  {$g-h$};
            \fatvertex{1}{0};
            \draw (1, 0) --  (2, 0);
        \end{tikzpicture}}
       \\
        \hline
    \end{tabular}
    }
    \caption{The $p$-covers corresponding to vertices in $\MS$. Indicated are vertex genera.}
    \label{tab:vertices}
\end{table}

\subsection{The equivariant (loop or) weight locus}

This subsection and the next follow the proof of \cite[Theorem~1.1]{CGP2} closely. Let us first recall the definition of the equivariant weight locus $(\MS)^w$:
\[ 
    \Big\{ \pi \colon \tilde \Gamma \to \Gamma \mathrel{\Big |} \text{there exists } v \in V(\Gamma) \text{ such that } g(v) \geq 0\text{ and }\sum_{\tilde v \in \pi^{-1}(v)} g(\tilde v) \geq p \Big\}, 
\]
and the equivariant loop or weight locus $\MS^{lw}$ defined as the closure of the locus of $p$-covers $\pi \colon \tilde \Gamma \to \Gamma$ such that $\pi \in \MS^{w}$ or there is a loop-edge $e$ in $\Gamma$ such that $\pi^{-1}(e)$ consists of $p$ loop edges.

\medskip

Consider the vertex $B_1 \in \MS$ corresponding to the $p$-cover $\tilde\Gamma \to \Gamma$ with
\begin{itemize}
	\item $\Gamma$ consisting of a single edge $e = (v_1, v_2)$ with $v_1 \neq v_2$ and $g(v_1) = g-1$ as well as $g(v_2) = 1$, and
	\item $\tilde \Gamma$ having a single vertex over $v_1$ and $p$ vertices over $v_2$.
\end{itemize}
See \cref{tab:vertices} for a sketch of $B_1$. 
Moreover, for any given $p$-cover of graphs $\tilde \Gamma \to \Gamma$ we will call an edge $e \in \Gamma$ an \emph{equivariant 1-bridge} if contracting all edges in $E(\Gamma) \setminus \{e\}$ takes $\pi$ to $B_1$.
In the following we will abuse notation and denote the property $\{B_1\}$ simply by $B_1$ as well.

Note that $B_1$ is a vertex of $\MS$, but it is also a vertex of the equivariant weight locus $\MS^w$ considered as a symmetric $\Delta$-complex in its own right. The closed stars of $B_1$ in these two complexes are
\begin{equation*}\label{eq:WeightLoopLocus}
	(\MS)_{B_1} = \MS^{lw} \quad \text{ and } \quad (\MS^w)_{B_1} = \MS^w
\end{equation*}
and therefore, Lemma~\ref{lem:loop_weight_locus} below, combined with \cite[Proposition~4.11]{CGP2} proves the parts of Theorem~\ref{thm:contractible_loci} about the contractibility of $\MS^w$ and $\MS^{lw}$, see \cref{cor:loop_weight}.

\begin{lemma} \label{lem:loop_weight_locus}
	The properties $P = \emptyset$ and $Q = B_1$, considered as properties in $\MS$ (resp. $\MS^w$) satisfy conditions~\ref{condition1} and \ref{condition2}.
\end{lemma}

The key to proving Lemma~\ref{lem:loop_weight_locus} is a good understanding of uncontractions of equivariant 1-bridges. This is achieved in the next lemma. The upshot is that the number of possible uncontractions of 1-bridges at any given vertex is an intrinsic property of that vertex, i.e. it is preserved under isomorphism.

\begin{lemma} \label{lem:uncontractions_br1}
	Let $\pi \colon \tilde \Gamma \to \Gamma$ be a $p$-cover and let $v \in \Gamma$ be a vertex. Denote by $l(v)$ the number of loop-edges rooted at $v$ which have $p$ loop-edges as preimages under $\pi$. Denote by $g^{\uparrow}(v) = \sum_{\tilde v \in \pi^{-1}(v)} g(\tilde v)$ the total vertex genus in the fiber over $v$. 
	Then the maximal number of possible uncontractions of equivariant 1-bridges at $v$, which still preserve stability of the cover is given by:
	\begin{equation}\label{eq:max1Bridges}
		\begin{cases}
			0 & \text{if } \ 2g(v) - 2 + \val(v) = 1 \\
			\min\left(g(v), \left\lfloor \frac{g^{\uparrow}(v)}{p} \right\rfloor \right) + l(v)  & \text{otherwise.}
		\end{cases}
	\end{equation}
\end{lemma}

\begin{proof}
    It is not so hard to see that if $2g(v) -2 + \val(v) = 1$, then any uncontraction (of any type) breaks the stability of the target graph at the vertex $v$. So without loss of generality we assume that $2g(v) - 2 + \val(v) > 1$. A priori we claim that for each loop with $p$ preimage loops and for each unit of $p$ in the total vertex genus in the fiber over $v$, an equivariant 1-bridge can be uncontracted. There are three things that can go wrong with this, namely the suggested uncontractions can:
 
    \begin{enumerate}
        \item  break stability (hence the case distinction),
        \item  produce negative genus on the target (hence the minimum), 
        \item  break any of the conditions of $p$-covers in \cref{def:Gcover}, i.e. the $\Z/p$-action, harmonicity, local Riemann-Hurwitz \eqref{eq:local-RH}, or the dilation flow condition.
    \end{enumerate}
    One can check that the third issue does not arise. In fact, the only part of this issue which requires some care, is that the uncontractions due to vertex genus are possible without violating the Riemann-Hurwitz condition. Indeed, if $v$ is free, i.e. it has $p$ preimages, then $\pi$ is a local isomorphism and the local Riemann-Hurwitz condition is trivial. If $v$ is dilated with $\pi^{-1}(v) = \{\tilde v\}$, then the proposed uncontraction lowers the genus of $\tilde v$ by $p$ and that of $v$ by 1. This subtracts $2p$ on both sides of Equation~\eqref{eq:local-RH}, so the local Riemann-Hurwitz condition still holds. The second issue is taken care of by the minimum in the formula.
    
    Concerning the first issue, note that after performing all the uncontractions suggested above, the remaining genus of $v$ is
	\[ 
    g(v) - \min \left( g(v), \left\lfloor \frac{g^{\uparrow}(v)}{p} \right\rfloor \right)
    \] 
	and the new valence of $v$ is
	\[ 
        \val(v) + \min \left( g(v), \left\lfloor \frac{g^{\uparrow}(v)}{p} \right\rfloor \right) - l(v)  .
    \]
	Plugging these new values into the definition of stability in Equation~\eqref{eq:stability} shows that the suggested uncontractions render the cover unstable if and only if:
 
	\[ 
     2 \underbrace{\left( g(v) - \min \left( g(v), \left\lfloor \frac{g^{\uparrow}(v)}{p} \right\rfloor \right) \right)}_{\geq 0}
	+ \underbrace{\val(v) - l(v)}_{\geq l(v) \geq 0}
	+ \underbrace{\min \left( g(v), \left\lfloor \frac{g^{\uparrow}(v)}{p} \right\rfloor \right)}_{\geq 0} \leq 2  . 
   \]
	Consider the term in the middle: the assumption $\val(v) - l(v) = 0$ leads to $\val(v) = 0$ which is not possible. So the middle term is positive and hence we also obtain $g = \min \left( g(v), \left\lfloor \frac{g^{\uparrow}(v)}{p} \right\rfloor \right)$ and are left with 
	\[ \underbrace{\val(v) - l(v)}_{\geq 1}
	+ \underbrace{g(v)}_{\geq 0} \leq 2  . \]
	Now distinguish cases for $g(v)$. 

 \begin{enumerate}[leftmargin = 40pt, itemindent = 65pt] 
     \item[\textbf{Case} \boldmath{$g(v) = 1$}\textbf{:}] In this case we have $\val(v) - l(v) = 1$ and the only possibilities are $\val(v) = 1$ and $l(v) = 0$ or $\val(v) = 2$ and $l(v) = 1$. If $\val(v) = 1$ and $l(v) = 0 $ then we were in the first case of \eqref{eq:max1Bridges} to begin with, which is a contradiction. So $\val(v) = 2$ and $l(v) = 1$ and in this case the graph $\Gamma$ has to be Butterfly cover in \cref{tab:vertices} and one can easily check that the formula in \cref{eq:max1Bridges} holds in this case. 

     \item[\textbf{Case} \boldmath{$g(v) = 0$}\textbf{:}] In this case we either have $\val(v) - l(v) = 1$ in which case the graph $\Gamma$ is not stable at $v$ which is a contradiction. So we deduce that $\val(v) - l(v) = 2$ and hence we have either $\val(v) = 3$ and $l(v) = 1$, in which case we are back in the first case of \cref{eq:max1Bridges} again. The other option is $\val(v) = 4$ and $l(v) = 2$, but this situation does not arise because of the local Riemann-Hurwitz condition \eqref{eq:local-RH}.
 \end{enumerate}
    To sum up, \cref{eq:max1Bridges} is correct in every case, thus finishing the proof.
\end{proof}

\begin{proof}[Proof of Lemma~\ref{lem:loop_weight_locus}]

    Notice that taking $P = \emptyset$ makes Condition~\ref{condition1} void, so we focus on proving Condition~\ref{condition2}.
	The argument is practically the same for $\MS$ and $\MS^{w}$, so let us work with the former first. Let $\pi \colon \tilde \Gamma \to \Gamma$ be an arbitrary $p$-cover, corresponding to some simplex $\sigma$ in $\MS$. We have to show that $\sigma$ has a canonical maximal co-$B_1$ face, or equivalently, that there is a unique (up to isomorphism) maximal stable uncontraction $\rho$ of $\pi$ by equivariant 1-bridges and that this is canonical.
	
	To this end, let $v \in \Gamma$ be a vertex. The exact statement about the possible uncontractions of said type at $v$ is given in Lemma~\ref{lem:uncontractions_br1}. In particular, we see that the number of possible uncontractions at $v$ is a property which is preserved under isomorphism. Note further, that all new vertices that are created in such an uncontraction are either 1-valent of genus 1 or 3-valent of genus 0, so that the new vertices cannot be uncontracted any further. As a consequence, we see that if we perform all the uncontractions suggested by Lemma~\ref{lem:uncontractions_br1} at each vertex individually, then we arrive at a unique maximal $\rho$. 
	
	Moreover, $\rho$ is canonical because any automorphism of $\pi$ can only permute vertices of the same genus and valence. Hence, again by Lemma~\ref{lem:uncontractions_br1} every automorphism of $\pi$ permutes only vertices with the same uncontraction behavior and thus we obtain a lift to an automorphism of $\rho$. This proves the claim when $P$ and $Q$ are considered properties of $\MS$.
	
	To see that $P$ and $Q$ still satisfy Condition~\ref{condition2} when considered as properties of $\MS^w$, we simply note that if the cover $\pi$ was in $\MS^{w}$, then after the uncontractions suggested by \cref{lem:uncontractions_br1} there is still a fiber of total vertex genus at least $p$. So the maximal uncontraction lives in $\MS^w$ as well and we have proved the entire claim.
\end{proof}

\begin{corollary} \label{cor:loop_weight}
    The loci $\MS^{w}$ and $\MS^{wl}$ are contractible in $\MS$.
\end{corollary}

\begin{proof}
    \cref{lem:loop_weight_locus} shows via \cite[Proposition~4.11]{CGP2} that $\MS^{lw}$ deformation retracts to the locus $(\MS)_{B_1, 0}$, which is the locus of $p$-covers all of whose edges are equivariant 1-bridges. There is a unique maximal cell in $(\MS)_{B_1, 0}$ and its automorphism group acts on the simplex as $S_{g-1}$ by permuting the vertices. In particular, \cite[Proposition~5.1]{AllcockCoreyPayne} shows that $(\MS)_{B_1, 0}$ is contractible. This shows the claim for $\MS^{lw}$. To see that $\MS^{w}$ is contractible, we simply note that $(\MS)_{B_1, 0} = (\MS^{w})_{B_1, 0}$.
\end{proof}

\subsection{The equivariant bridge locus}

We move on to showing contractibility of the equivariant bridge locus $\MS^{br}$ in $\MS$, still following the strategy used to prove \cite[Theorem~1.1]{CGP2}. 
Recall that $\MS^{br}$ is defined as the closure of the union of $\MS^{lw}$ and
	\begin{equation*}
		\left\{ 
		\pi \colon \tilde \Gamma \to \Gamma \mathrel{\bigg |} 
		\begin{minipage}{0.65\textwidth}
			there is a bridge-edge $e \in \Gamma$ such that $|\pi^{-1}(e)| = p$ and every $\tilde e \in \pi^{-1}(e)$ is a bridge as well
		\end{minipage}		
		\right\}.
	\end{equation*}

\medskip

We start by defining a generalization of the $p$-cover~$B_1$.

\begin{definition}\label{def:coverB_h}
	Let $1 \leq h \leq g-1$. Define the $p$-cover $B_h \colon \tilde \Gamma \to \Gamma$ (see \cref{tab:vertices}) with:
	\begin{itemize}
		\item $\Gamma$ consisting of a single edge $e = (v_1, v_2)$ with $v_1 \neq v_2$ and $g(v_1) = g-h$, $g(v_2) = h$, and
		\item $\tilde \Gamma$ having a single vertex covering $v_1$ and $p$ vertices covering $v_2$.
	\end{itemize}
	We will abuse notation and denote the vertex of $\MS$ corresponding to $B_h$, as well as the property consisting only of this vertex by $B_h$ again. Moreover, an edge $e \in \Gamma$ in a $p$-cover $\pi \colon \tilde \Gamma \to \Gamma$ is called \emph{equivariant $h$-bridge} if contracting all edges in $\Gamma \setminus \{e\}$ takes $\pi$ to $B_h$.
\end{definition}

The following \cref{lem:properties_are_good} is the technical core of this subsection. Combining this with \cite[Corollary~4.18]{CGP2} gives a strong deformation retraction from the equivariant bridge locus $\MS^{br} = (\MS)_{B_1 \cup \cdots \cup B_{g-1}}$ to the contractible subcomplex $\MS^{lw} = (\MS)_{B_1}$. 

\begin{lemma}
    \label{lem:properties_are_good}
    Consider the sequence
    \[ 
     B_1, B_2, \ldots, B_{g-1} 
    \]
    of properties in $\MS$. Let $1 < h < g$ be an integer and set $P = B_1 \cup \cdots \cup B_{h-1}$ and $Q = \{B_h\}$. Then $P$ and $Q$ satisfy conditions~\ref{condition1} to~\ref{condition3}. Consequently, there is a strong deformation retraction
    \[
        \MS^{br} = (\MS)_{B_1 \cup \cdots \cup B_{g-1}} \searrow (\MS)_{B_1} = \MS^{lw}. 
    \]
\end{lemma}

In order to verify Condition~\ref{condition2} in the proof of \cref{lem:properties_are_good}, we introduce the concepts of \enquote{cut components} and \enquote{bridge articulation points}. 

\begin{definition} \label{def:removing_vertex}
    Let $\pi \colon \tilde \Gamma \to \Gamma$ be a $p$-cover and let $v$ be a vertex in $\Gamma$. Then there exist connected subgraphs  $\Gamma_1, \dots, \Gamma_r$ of $\Gamma$, each containing $v$, such that 
    \[
        \Gamma = \Gamma_1 \vee_v \dots \vee_v \Gamma_r,
    \]
    where $\vee_v$ means gluing these subgraphs at $v$. Moreover, there exists a unique sequence of such subgraphs (up to relabeling) such that $r$ is maximal. These are what we call the \emph{cut components} of $\Gamma$ at $v$. In this situation we say that $\Gamma_i$ has \emph{trivial preimage away from $v$} if one the the following hold:
    \begin{enumerate}        
        \item the preimage $\pi^{-1}(\Gamma_i)$ consists of $p$ disjoint copies of $\Gamma_i$,

        \item the preimage $\pi^{-1}(\Gamma_i)$ is the gluing of $p$-many copies of $\Gamma_i$ at the vertex $v$.
    \end{enumerate}
     \cref{tab:cutComponentsExamples} depicts a few examples.
\end{definition}

\begin{table}[ht]
    \centering
    \begin{tabular}{|c|c|c|c|}
          \hline
          \small{Trivial preimage} & \small{Trivial preimage} & \small{Non-trivial preimage} &\small{Non-trivial preimage} \\ 
          \hline
          \scalebox{0.7}{         
                \tikzset{every picture/.style={line width=0.75pt}} 
                
                \begin{tikzpicture}[x=0.75pt,y=0.75pt,yscale=-1,xscale=1]
                
                \draw    (264.84,111.1) -- (313.14,111.1) ;
                \draw [color={rgb, 255:red, 0; green, 0; blue, 255 }  ,draw opacity=1 ][line width=1.5]    (313.14,111.1) .. controls (344.3,87.18) and (346.73,135.32) .. (313.32,111.1) ;
                \draw    (314.85,67.34) .. controls (346.02,43.43) and (348.45,91.57) .. (315.04,67.34) ;
                \draw    (315.5,20.53) .. controls (347.76,-3.24) and (348.24,44.9) .. (315.68,20.53) ;
                \draw    (292,73.67) -- (292,95.18) ;
                \draw [shift={(292,97.18)}, rotate = 270] [color={rgb, 255:red, 0; green, 0; blue, 0 }  ][line width=0.75]    (10.93,-3.29) .. controls (6.95,-1.4) and (3.31,-0.3) .. (0,0) .. controls (3.31,0.3) and (6.95,1.4) .. (10.93,3.29)   ;
                \draw    (315.5,43.53) .. controls (346.67,19.62) and (349.09,67.76) .. (315.68,43.53) ;
                \draw    (266.19,20.45) -- (315.5,20.53) ;
                \draw    (266.19,43.35) -- (315.5,43.53) ;
                \draw    (266.37,67.34) -- (314.85,67.34) ;
                \draw  [fill={rgb, 255:red, 0; green, 0; blue, 0 }  ,fill opacity=1 ] (313.43,20.53) .. controls (313.43,18.98) and (314.44,17.71) .. (315.68,17.71) .. controls (316.93,17.71) and (317.93,18.98) .. (317.93,20.53) .. controls (317.93,22.09) and (316.93,23.35) .. (315.68,23.35) .. controls (314.44,23.35) and (313.43,22.09) .. (313.43,20.53) -- cycle ;
                \draw  [fill={rgb, 255:red, 0; green, 0; blue, 0 }  ,fill opacity=1 ] (313.43,43.53) .. controls (313.43,41.98) and (314.44,40.71) .. (315.68,40.71) .. controls (316.93,40.71) and (317.93,41.98) .. (317.93,43.53) .. controls (317.93,45.09) and (316.93,46.35) .. (315.68,46.35) .. controls (314.44,46.35) and (313.43,45.09) .. (313.43,43.53) -- cycle ;
                \draw  [fill={rgb, 255:red, 0; green, 0; blue, 0 }  ,fill opacity=1 ] (314.85,67.34) .. controls (314.85,65.78) and (315.86,64.52) .. (317.1,64.52) .. controls (318.34,64.52) and (319.35,65.78) .. (319.35,67.34) .. controls (319.35,68.9) and (318.34,70.16) .. (317.1,70.16) .. controls (315.86,70.16) and (314.85,68.9) .. (314.85,67.34) -- cycle ;
                \draw  [color={rgb, 255:red, 0; green, 0; blue, 255 }  ,draw opacity=1 ][fill={rgb, 255:red, 0; green, 0; blue, 255 }  ,fill opacity=1 ] (313.14,111.1) .. controls (313.14,109.54) and (314.14,108.28) .. (315.38,108.28) .. controls (316.63,108.28) and (317.63,109.54) .. (317.63,111.1) .. controls (317.63,112.65) and (316.63,113.92) .. (315.38,113.92) .. controls (314.14,113.92) and (313.14,112.65) .. (313.14,111.1) -- cycle ;
                
                \draw (292,118) node [anchor=north west][inner sep=0.75pt]    {$v$};
                \draw (238,48.07) node [anchor=north west][inner sep=0.75pt]    {$\dotsc $};

                \end{tikzpicture}
                
        }
        &
        \scalebox{0.7}{

                \tikzset{every picture/.style={line width=0.75pt}} 
                
                \begin{tikzpicture}[x=0.75pt,y=0.75pt,yscale=-1,xscale=1]
                
                \draw    (306.54,43.91) .. controls (375.24,69.99) and (338.43,101.52) .. (306.73,43.91) ;
                \draw    (306.73,43.91) .. controls (368.84,21.27) and (367.24,65.7) .. (306.91,43.91) ;
                \draw    (298.19,73.47) -- (298.19,94.99) ;
                \draw [shift={(298.19,96.99)}, rotate = 270] [color={rgb, 255:red, 0; green, 0; blue, 0 }  ][line width=0.75]    (10.93,-3.29) .. controls (6.95,-1.4) and (3.31,-0.3) .. (0,0) .. controls (3.31,0.3) and (6.95,1.4) .. (10.93,3.29)   ;
                \draw    (306.91,43.91) .. controls (335.23,-15.28) and (372.04,16.25) .. (307.09,43.91) ;
                \draw    (256,43.91) -- (304.48,43.91) ;
                \draw    (258.84,113.1) -- (307.32,113.1) ;
                \draw  [fill={rgb, 255:red, 0; green, 0; blue, 0 }  ,fill opacity=1 ] (304.48,43.91) .. controls (304.48,42.36) and (305.48,41.09) .. (306.73,41.09) .. controls (307.97,41.09) and (308.97,42.36) .. (308.97,43.91) .. controls (308.97,45.47) and (307.97,46.73) .. (306.73,46.73) .. controls (305.48,46.73) and (304.48,45.47) .. (304.48,43.91) -- cycle ;
                \draw  [color={rgb, 255:red, 0; green, 0; blue, 255 }  ,draw opacity=1 ][fill={rgb, 255:red, 0; green, 0; blue, 255 }  ,fill opacity=1 ] (304.89,113.1) .. controls (304.89,111.54) and (305.89,110.28) .. (307.14,110.28) .. controls (308.38,110.28) and (309.38,111.54) .. (309.38,113.1) .. controls (309.38,114.65) and (308.38,115.92) .. (307.14,115.92) .. controls (305.89,115.92) and (304.89,114.65) .. (304.89,113.1) -- cycle ;
                \draw [color={rgb, 255:red, 0; green, 0; blue, 255 }  ,draw opacity=1 ][line width=1.5]    (307.95,113.1) .. controls (339.12,89.18) and (341.55,137.32) .. (308.14,113.1) ;
                
                \draw (286,118) node [anchor=north west][inner sep=0.75pt]    {$v$};
                \draw (226,50.07) node [anchor=north west][inner sep=0.75pt]    {$\dotsc $};

                \end{tikzpicture}
                
            }
            &
            \scalebox{0.7}{            
                
                \tikzset{every picture/.style={line width=0.75pt}} 
                
                \begin{tikzpicture}[x=0.75pt,y=0.75pt,yscale=-1,xscale=1]
                
                \draw    (318.74,65.34) .. controls (360.5,66.93) and (353.5,16.93) .. (318.55,21.35) ;
                \draw    (318.55,21.35) .. controls (341.5,24.93) and (333.5,40.93) .. (318.55,41.35) ;
                \draw    (302.83,76.91) -- (302.83,98.42) ;
                \draw [shift={(302.83,100.42)}, rotate = 270] [color={rgb, 255:red, 0; green, 0; blue, 0 }  ][line width=0.75]    (10.93,-3.29) .. controls (6.95,-1.4) and (3.31,-0.3) .. (0,0) .. controls (3.31,0.3) and (6.95,1.4) .. (10.93,3.29)   ;
                \draw    (318.55,41.35) .. controls (336.5,45.93) and (335.5,61.93) .. (318.55,65.34) ;
                \draw    (270.08,41.35) -- (318.55,41.35) ;
                \draw    (270.26,65.34) -- (318.74,65.34) ;
                \draw  [fill={rgb, 255:red, 0; green, 0; blue, 0 }  ,fill opacity=1 ] (316.31,41.35) .. controls (316.31,39.79) and (317.31,38.53) .. (318.55,38.53) .. controls (319.8,38.53) and (320.8,39.79) .. (320.8,41.35) .. controls (320.8,42.9) and (319.8,44.17) .. (318.55,44.17) .. controls (317.31,44.17) and (316.31,42.9) .. (316.31,41.35) -- cycle ;
                \draw  [fill={rgb, 255:red, 0; green, 0; blue, 0 }  ,fill opacity=1 ] (316.49,65.34) .. controls (316.49,63.78) and (317.5,62.52) .. (318.74,62.52) .. controls (319.98,62.52) and (320.99,63.78) .. (320.99,65.34) .. controls (320.99,66.9) and (319.98,68.16) .. (318.74,68.16) .. controls (317.5,68.16) and (316.49,66.9) .. (316.49,65.34) -- cycle ;
                \draw    (270.08,21.35) -- (318.55,21.35) ;
                \draw  [fill={rgb, 255:red, 0; green, 0; blue, 0 }  ,fill opacity=1 ] (316.31,21.35) .. controls (316.31,19.79) and (317.31,18.53) .. (318.55,18.53) .. controls (319.8,18.53) and (320.8,19.79) .. (320.8,21.35) .. controls (320.8,22.9) and (319.8,24.17) .. (318.55,24.17) .. controls (317.31,24.17) and (316.31,22.9) .. (316.31,21.35) -- cycle ;
                \draw    (273.39,109.1) -- (321.86,109.1) ;
                \draw  [color={rgb, 255:red, 0; green, 0; blue, 255 }  ,draw opacity=1 ][fill={rgb, 255:red, 0; green, 0; blue, 255 }  ,fill opacity=1 ] (319.43,109.1) .. controls (319.43,107.54) and (320.44,106.28) .. (321.68,106.28) .. controls (322.92,106.28) and (323.93,107.54) .. (323.93,109.1) .. controls (323.93,110.65) and (322.92,111.92) .. (321.68,111.92) .. controls (320.44,111.92) and (319.43,110.65) .. (319.43,109.1) -- cycle ;
                \draw [color={rgb, 255:red, 0; green, 0; blue, 255 }  ,draw opacity=1 ][line width=1.5]    (320.5,109.1) .. controls (351.67,85.18) and (354.09,133.32) .. (320.68,109.1) ;
                
                \draw (305,118) node [anchor=north west][inner sep=0.75pt]    {$v$};
                \draw (233,37.07) node [anchor=north west][inner sep=0.75pt]    {$\dotsc $};

                \end{tikzpicture}
                
            }
            &
            \scalebox{0.7}{
                \tikzset{every picture/.style={line width=0.75pt}} 
                
                \begin{tikzpicture}[x=0.75pt,y=0.75pt,yscale=-1,xscale=1]
                
                \draw    (231.84,115.1) -- (280.14,115.1) ;
                \draw    (259,77.67) -- (259,99.18) ;
                \draw [shift={(259,101.18)}, rotate = 270] [color={rgb, 255:red, 0; green, 0; blue, 0 }  ][line width=0.75]    (10.93,-3.29) .. controls (6.95,-1.4) and (3.31,-0.3) .. (0,0) .. controls (3.31,0.3) and (6.95,1.4) .. (10.93,3.29)   ;
                \draw    (233.19,24.45) -- (282.5,24.53) ;
                \draw    (233.19,47.35) -- (282.5,47.53) ;
                \draw    (233.37,71.34) -- (281.85,71.34) ;
                \draw  [fill={rgb, 255:red, 0; green, 0; blue, 0 }  ,fill opacity=1 ] (280.43,24.53) .. controls (280.43,22.98) and (281.44,21.71) .. (282.68,21.71) .. controls (283.93,21.71) and (284.93,22.98) .. (284.93,24.53) .. controls (284.93,26.09) and (283.93,27.35) .. (282.68,27.35) .. controls (281.44,27.35) and (280.43,26.09) .. (280.43,24.53) -- cycle ;
                \draw  [fill={rgb, 255:red, 0; green, 0; blue, 0 }  ,fill opacity=1 ] (280.43,47.53) .. controls (280.43,45.98) and (281.44,44.71) .. (282.68,44.71) .. controls (283.93,44.71) and (284.93,45.98) .. (284.93,47.53) .. controls (284.93,49.09) and (283.93,50.35) .. (282.68,50.35) .. controls (281.44,50.35) and (280.43,49.09) .. (280.43,47.53) -- cycle ;
                \draw  [fill={rgb, 255:red, 0; green, 0; blue, 0 }  ,fill opacity=1 ] (281.85,71.34) .. controls (281.85,69.78) and (282.86,68.52) .. (284.1,68.52) .. controls (285.34,68.52) and (286.35,69.78) .. (286.35,71.34) .. controls (286.35,72.9) and (285.34,74.16) .. (284.1,74.16) .. controls (282.86,74.16) and (281.85,72.9) .. (281.85,71.34) -- cycle ;
                \draw  [color={rgb, 255:red, 0; green, 0; blue, 255 }  ,draw opacity=1 ][fill={rgb, 255:red, 0; green, 0; blue, 255 }  ,fill opacity=1 ] (280.14,115.1) .. controls (280.14,113.54) and (281.14,112.28) .. (282.38,112.28) .. controls (283.63,112.28) and (284.63,113.54) .. (284.63,115.1) .. controls (284.63,116.65) and (283.63,117.92) .. (282.38,117.92) .. controls (281.14,117.92) and (280.14,116.65) .. (280.14,115.1) -- cycle ;
                \draw [color={rgb, 255:red, 0; green, 0; blue, 255 }  ,draw opacity=1 ][line width=1.5]    (282.38,115.1) .. controls (292.57,100.97) and (304.57,101.93) .. (315.2,115.1) ;
                \draw  [color={rgb, 255:red, 0; green, 0; blue, 255 }  ,draw opacity=1 ][fill={rgb, 255:red, 0; green, 0; blue, 255 }  ,fill opacity=1 ] (312.95,115.1) .. controls (312.95,113.54) and (313.96,112.28) .. (315.2,112.28) .. controls (316.44,112.28) and (317.45,113.54) .. (317.45,115.1) .. controls (317.45,116.65) and (316.44,117.92) .. (315.2,117.92) .. controls (313.96,117.92) and (312.95,116.65) .. (312.95,115.1) -- cycle ;
                \draw [color={rgb, 255:red, 0; green, 0; blue, 255 }  ,draw opacity=1 ][line width=1.5]    (282.45,116.1) .. controls (290.57,126.93) and (304.57,129.97) .. (315.2,115.1) ;
                \draw    (284.93,24.53) .. controls (302.88,23.12) and (311.86,32.93) .. (318.92,46.35) ;
                \draw    (282.5,24.53) .. controls (294.45,33.12) and (304.86,40.93) .. (318.92,46.35) ;
                \draw    (286.35,71.34) .. controls (301.11,70.93) and (316.86,61.93) .. (318.92,46.35) ;
                \draw    (286.35,71.34) .. controls (296.11,60.93) and (304.11,53.93) .. (318.92,46.35) ;
                \draw    (282.5,47.53) .. controls (286.45,40.12) and (302.86,40.93) .. (318.92,46.35) ;
                \draw    (282.68,47.53) .. controls (287.63,53.12) and (301.86,53.93) .. (318.92,46.35) ;
                \draw  [fill={rgb, 255:red, 0; green, 0; blue, 0 }  ,fill opacity=1 ] (316.67,46.35) .. controls (316.67,44.79) and (317.68,43.53) .. (318.92,43.53) .. controls (320.16,43.53) and (321.17,44.79) .. (321.17,46.35) .. controls (321.17,47.9) and (320.16,49.17) .. (318.92,49.17) .. controls (317.68,49.17) and (316.67,47.9) .. (316.67,46.35) -- cycle ;
                
                \draw (265,123) node [anchor=north west][inner sep=0.75pt]    {$v$};
                \draw (205,52.07) node [anchor=north west][inner sep=0.75pt]    {$\dotsc $};

                \draw (233,125) node [anchor=north west][inner sep=0.75pt]    {\color{white} $A$};
                
                \end{tikzpicture}
                
            }\\
          \hline
    \end{tabular}
     \caption{In the first and second $3$-cover the cut component $\Gamma_1$ at $v$ (colored in blue) has trivial preimage away from $v$. This fails to be the case in the remaining two examples.  }
    \label{tab:cutComponentsExamples}
\end{table}

\begin{definition} \label{def:bridge_articulation}
	Let $\pi \colon \tilde \Gamma \to \Gamma$ be a $p$-cover and let $v$ be a vertex in $\Gamma$.
	We call $v$ a \emph{bridge articulation point of type $h$} in $\pi$ if and only if there is a cut component $\Gamma_0$ of $\Gamma$ at $v$ which has trivial preimage away from $v$ (see  Definition~\ref{def:removing_vertex}) and $g(\Gamma_0) = h + g(v)$.
\end{definition}

The idea of \cref{def:bridge_articulation} is as follows. Let $\pi \colon \tilde\Gamma \to \Gamma$ be a $p$-cover subject to the following:
\begin{assumption} \label{assumption}
    There are no equivariant $h'$-bridges and no vertices at which an uncontraction of an equivariant $h'$-bridge is possible, for any $1 \leq h' < h$.
\end{assumption}
Under this assumption, a vertex $v$ of $\Gamma$ is a bridge articulation point of type $h$ if and only if it is incident to an equivariant $h$-bridge or an uncontraction of an equivariant $h$-bridge is possible at $v$. It is important to note that the vertex $v$ cannot supplement any genus to the uncontraction of an equivariant bridge: if this were possible, then an equivariant $1$-bridge could be uncontracted at $v$, which contradicts \cref{assumption}.

The upshot of the following lemma is similar to Lemma~\ref{lem:uncontractions_br1}; that is, the number of possible uncontractions of equivariant bridges of type $h$ at any given vertex $v$ is a local and intrinsic property of $v$.

\begin{lemma} \label{lem:bridge_articulation_points}
	Let $\pi \colon \tilde \Gamma \to \Gamma$ be a $p$-cover and let $2 \leq h \leq g-1$. Suppose that $\pi$ satisfies \cref{assumption}.
	Let $v$ be a vertex in $\Gamma$ and denote by $\Gamma_1, \ldots, \Gamma_r, \Gamma_{r+1} \dots, \Gamma_n$ the cut components of $\Gamma$ at $v$. Without loss of generality we order them such that $\Gamma_1, \ldots, \Gamma_r$ are all of the cut components with trivial preimage away from $v$, with genus $g(\Gamma_i) = h + g(v)$, and $\val_{\Gamma_i}(v) \geq 2$. Finally, we denote by $\Gamma'$ the remaining graph
    $\Gamma_{r+1} \vee_v \dots \vee_v \Gamma_{n}$. Then the maximal number of equivariant $h$-bridges that can be uncontracted at $v$ is precisely:
	\begin{equation} \label{eq:hBridges}
		\begin{cases}
			r - 1  & \text{if } \ 2g(v) - 2  +  r + \val_{\Gamma'}(v) \leq 0 \\
			  r      & \text{otherwise.}
		\end{cases}
	\end{equation}
	Moreover, we have $\pi \in B_h^\ast(\MS)$ if and only if there is a bridge articulation point of type $h$ in $\pi$.
\end{lemma}

\begin{proof}
	Analogous to Lemma~\ref{lem:uncontractions_br1}: a priori every cut component of genus $h$ with trivial preimage away from $v$ is a candidate for uncontracting an equivariant $h$-bridge. These uncontractions do not break the local Riemann-Hurwitz condition nor do they impact the dilation flow conditions in \cref{def:Gcover}. The only thing that might prevent some of these uncontractions is the stability condition. For example, uncontracting an equivariant bridge from a cut component $\Gamma_i$ with $\val_{\Gamma_i}(v) = 1$ would produce a 2-valent vertex of genus 0, since by \cref{assumption} the vertex $v$ cannot supplement any genus to the uncontraction of an equivariant bridge. This explains the definition of $r$.
 
    Besides this scenario, stability can only fail if we reduce the valence of $v$ to $\leq 2$ while having $g(v) = 0$ at the same time. This is equivalent to saying that after the suggested uncontractions, $v$ is still supposed to be stable, which by definition translates to
    \[2g(v) - 2  +  r + \val_{\Gamma'}(v) > 0.\]
    This explains the case distinction in \cref{eq:hBridges}. 

    Now we turn our attention to the last statement. The \enquote{only if} part of the statement is clear. For the other direction, suppose that $\pi$ has a bridge articulation point $v$ of type $h$. This means that there is a cut component $\Gamma_0$ at $v$ with trivial preimage away from $v$ and genus $h + g(v)$. If \cref{eq:hBridges} gives a positive number of uncontractions, then we are done. So let us assume that either $r = 0$ or $r = 1$ and $2g(v) - 1 + \val_{\Gamma'}(v) \leq 0$. In the first case, $\val_{\Gamma_0}(v) = 1$ and therefore $\Gamma_0$ is connected to $v$ via an equivariant $h$-bridge and we are done. Otherwise, $r = 1$ and $g(v) = 0$ as well as $\val_{\Gamma'}(v) \leq 1$. Note that every dilated edge incident to $v$ necessarily lies in $\Gamma'$ and moreover that the dilation flow conditions of \cref{def:Gcover} force $\val_{\Gamma^{\dil}}(v) \neq 1$. Hence, the one edge $e$ in $\Gamma'$ incident to $v$ has to be free. But then the local Riemann-Hurwitz condition forces $v$ to be free as well. Combining these observations we conclude that $e$ is already an equivariant $h$ bridge. 
\end{proof}

\begin{proof}[Proof of \cref{lem:properties_are_good}]
    Let $h \in \{ 2, \ldots, g-1 \}$ and set $P = B_1 \cup \cdots \cup B_{h-1}$ and $Q = B_h$.
    We check that $P, Q \subseteq \MS$ satisfy conditions~\ref{condition1} to~\ref{condition3}.
    \medskip

    \textbf{Condition~\ref{condition1}.} Let $\pi \colon \tilde \Gamma \to \Gamma$ be a $p$-cover. 
    We prove that, if $\pi$ is not in $P^\ast(\MS)$ then contracting an equivariant $h$-bridge will still not take us to $P^\ast(\MS)$.
    Let $e = (v_1, v_2)$ be an equivariant $h$-bridge of $\pi$ and denote the two connected components of $\Gamma \setminus \{e\}$ by $\Gamma_1$ and $\Gamma_2$. Assume for a contradiction that $\pi / e \in P^\ast(\MS)$. 
    This means that there is a $1 \leq h' < h$ such that:
    \begin{enumerate}
        \item there is an equivariant $h'$-bridge in $\pi/e$, or
        \item there is a vertex $v_0$ in $\Gamma/e$ where such a bridge can be uncontracted. 
    \end{enumerate}
    In the first case, this means that $\pi$ already contained that equivariant $h'$-bridge, contradicting $\pi \not\in P^\ast(\MS)$. If we are in the second case, then we show that the uncontraction is possible in $\pi$ as well. This is clear, if the vertex $v_0$ is a vertex of $\pi$ as well. It is therefore without loss of generality that we assume that $v_0$ is the image of $e$ under the edge contraction $\pi \to \pi/e$. 

    Suppose the $h'$-bridge which we can uncontract at $v_0$ is a 1-bridge. At least one of the endpoints of $e$ -- without loss of generality assume it is $v_1$ -- is a free vertex in $\pi$ (this is because $e$ is an $h$-bridge). If $\pi$ admits an uncontraction of a 1-bridge at $v_1$, then we are done. So without loss of generality, it does not. But this means that $g(v_1) = 0$ and similar for every preimage of $v_1$. This implies that the genus profile in the fiber over $v_0$ is the same as the genus profile in the fiber of $v_2$ and similarly the $\val_{\Gamma^\dil}(v_2) = \val_{\Gamma^\dil}(v_0)$. Therefore all parameters that determine whether an 1-bridge can be uncontracted are the same for $v_0$ and $v_2$, leading to the desired contradiction.

    It remains to argue that uncontractibility of an $h'$-bridge with $1 < h'$ at $v_0$ implies the same for either $v_1$ or $v_2$.     
    Since $v_0$ is an articulation point of $\pi/e$ of type $h'$, then at least one of the cut components $\Gamma_1, \dots, \Gamma_n$ of $\pi/e$ at $v_0$ has trivial preimage and genus $g(v_0) + h'$. Now, uncontracting the bridge $e$, we deduce one of $v_1$ or $v_2$ is an articulation point of $\pi$ of type $h'$. Then from \cref{lem:bridge_articulation_points} we conclude that $\pi \in B_{h'}^\ast(\MS) \subseteq P^\ast(\MS)$ which is the desired contradiction. 
    
    \medskip

    \textbf{Condition~\ref{condition2}.} 
    Let $\pi \colon \tilde \Gamma \to \Gamma$ be a $p$-cover which does not lie in $P^\ast(\MS)$. We have to show that $\pi$ admits a unique canonical maximal uncontraction $\rho$ by equivariant $h$-bridges. We are in the setting of Lemma~\ref{lem:bridge_articulation_points} which shows that at every vertex $v \in \Gamma$ there is an intrinsically defined maximal number of uncontractions by equivariant $h$-bridges. These uncontractions at various vertices do not interfere, in particular the order in which we perform these steps does not matter. Combining these uncontractions for all vertices yields a $p$-cover $\rho$ corresponding to the unique maximal co-$Q$ face of $\pi$. Furthermore, the face inclusion is canonical because any automorphism of $\pi$ only permutes $h$-bridge articulation points with isomorphic cut components and therefore it lifts to an automorphism of the uncontraction $\rho$ as well. 

    \medskip

    \textbf{Condition~\ref{condition3}.} Let $\pi \colon \tilde \Gamma \to \Gamma$ be a $p$-cover such that every edge is an equivariant $h$-bridge. In particular, $\Gamma$ is a tree and there is a leaf vertex $v$ which is free of genus $h$. But then, each preimage of $v$ has genus $h$ as well and therefore $\pi \in B_1^\ast(\MS) \subseteq P^\ast(\MS)$.
    \medskip
    
    This finishes the proof of the main part of the statement. The second assertion is obtained by virtue of \cite[Proposition~4.11]{CGP2}.
\end{proof}

\subsection{The sparsely connected locus} \label{subsec:scon}

In this subsection we assume that $p$ is an odd prime number. We show in this case that the sparsely connected (\emph{scon} for short) locus $\MS^{scon}$ deformation retracts onto the equivariant bridge locus. Combined with the results of the previous subsection, this shows that $\MS^{scon}$ is contractible. The proof we give does not hold for $p = 2$, see \cref{rem:scon_in_p_2}. As of now we do not know if the sparsely connected locus is contractible for $p = 2$ or not. Recall that $\MS^{scon}$ is defined as the closure of the union of $\MS^{br}$ and
	\begin{equation*}
		\left\{
		\pi \colon \tilde \Gamma \to \Gamma \mathrel{\bigg |}
		\begin{minipage}{0.65\textwidth}
			there is an edge $e \in \Gamma$ and a connected component $\Gamma_0 \subseteq \Gamma \setminus \{e\}$ such that $\pi^{-1}(\Gamma_0)$ is disconnected
		\end{minipage}
		\right\}.
	\end{equation*}

\medskip
 
The key idea in this section is to study \enquote{spiral properties}. To that end we define:

\begin{definition}\label{def:spirals}
    Let $a \in (\Z/p)^{\times}$ and define the \emph{spiral of type} $a$ as the free $p$-cover $S_a \colon \tilde \Gamma_a \to \Gamma_a$ such that:
    \begin{itemize}
        \item The target graph $\Gamma_a$ consists of a single loop edge $e$ rooted at a vertex $v$ of genus $g - 1$.
        \item Choose an orientation for $e$ and a preimage $\tilde v \in S_a^{-1}(v)$. In particular, this gives a $\Z/p$-equivariant identification of $S_a^{-1}(v)$ and $\Z/p$. The source graph $\tilde \Gamma_a$ is now defined by connecting every $\tilde w \in S_a^{-1}(v) = \Z/p$ with $\tilde w + a$ along the orientation of $e$.
    \end{itemize}
    We refer the reader to \cref{tab:vertices}.
    
    \medskip
    
    We say that $a$ is the \emph{ascent} of the spiral $S_a$. By abuse of notation, we denote the vertex in $\MS$ corresponding to $S_a$ by $S_a$ again. Furthermore, we set $S = \{S_1, \ldots, S_{(p-1)/2}\}$. We call an edge $e$ in a $p$-cover $\pi\colon \tilde \Gamma \to \Gamma$ a \emph{spiral edge of type $a$} if contracting all edges in $E(\Gamma) \setminus \{e\}$ yields the $p$-cover $S_a$. See \cref{fig:exampleOfGCovering} for an example.
\end{definition}

\begin{remark} \label{rem:OnSpirals}
    We collect some basic properties of spirals.
    
    \begin{enumerate}

        \item We stress that an edge $e$ for which $\pi|_{\pi^{-1}(e)} = S_a$ is \textbf{NOT} necessarily a spiral edge! For example, none of the edges in \cref{fig:G-contraction} is a spiral edge, not even those in the middle column.
        
        \item The idea behind \enquote{spiral ascents} is the following observation. For a finite abelian group $G$, and a $G$-torsor $T$ there is no canonical identification of $T$ and $G$. However, there is a canonical $G$-equivariant identification of $\Aut(T)$ and $G$. Now every simple closed oriented cycle $\gamma$ rooted at a vertex $v$ in the target graph of a free $p$-cover defines a $\Z/p$-equivariant automorphism of the fiber over $v$. The ascent of $\gamma$ is the element of $\Z/p$ to which this automorphism corresponds. If $\gamma$ is taken to be a loop-edge, we recover \cref{def:spirals}.
    
        \item We emphasize that spirals of different types are all isomorphic as Hurwitz covers, but not as $p$-covers. To see this, suppose there is a $\Z/p$-equivariant isomorphism $\theta \colon S_a \to S_b$. Choosing a vertex $\tilde v \in \tilde \Gamma_a$ to be identified with $0 \in \Z/p$ automatically gives an identification $V(\tilde \Gamma_b) \cong \Z/p$ based on setting $\theta(\tilde v) = 0$. Now $\theta$ being an isomorphism $\tilde \Gamma_a \to \tilde \Gamma_b$ means that the endpoint of a lift of $e$ starting at $\tilde v$ is mapped to the endpoint of a lift of $\theta(e)$ starting at $\theta(\tilde v)$, i.e. $\theta(a) = \pm b$, depending on the choice of lift of $\theta(e)$.         
        But $\theta$ being a $\Z/p$-equivariant isomorphism means that under the identifications $\Z/p \cong V(\tilde \Gamma_a)$ and $\Z/p \cong V(\tilde \Gamma_b)$ we have $\theta = \pm \id_{\Z/p}$. This shows that $a = \pm b$.

        \item From the previous observation, there are precisely $\frac{p-1}{2}$ different spiral types up to isomorphism: $a = 0$ is not a valid ascent as this would correspond to a disconnected cover. Moreover, ascents $a$ and $-a$ define isomorphic covers corresponding to the automorphism of the target graph given by flipping the loop (i.e. reversing the choice of orientation in the definition). From now on we will always consider ascents of spiral edges in $\{1,2, \dots, \frac{p-1}{2}\}$.
    \end{enumerate}
\end{remark}

The core technical statement of this subsection is the following lemma.

\begin{lemma}\label{lem:properties_are_good2}
	For all $g \geq 2$ and $p > 2$, the properties $P = \{B_1, \ldots , B_{g-1}\}$ and $Q = S = \{S_1, \ldots, S_{(p-1)/2}\}$ satisfy the conditions~\ref{condition1} to~\ref{condition3} in $\MS$. Consequently, there exists a strong deformation retraction 
    \[
        \MS^{scon} = (\MS)_{P \cup Q} \searrow (\MS)_P = \MS^{br}.
    \]
\end{lemma}

In the proof of \cref{lem:properties_are_good2}, given a $p$-cover $\pi\colon \tilde \Gamma \to \Gamma$, we will describe maximal sequences of uncontractions by spiral edges in $\pi$. Similar to the treatment of equivariant bridges we work with a suitable notion of \enquote{articulation point}, i.e. vertices $v\in V(\Gamma)$ at which an uncontraction of a spiral edge can be carried out. However, similar to \cref{assumption}, in the proof of \cref{lem:properties_are_good2} we will have the liberty to assume that $\pi$ is not in the equivariant bridge locus. The following definition already makes use of this, i.e. we claim that up to breaking stability, at every spiral articulation point a spiral edge can be uncontracted (Lemma~\ref{lem:uncontracting_spirals}), but we do \emph{not} claim that uncontractions of spirals are possible only at spiral articulation points (Lemma~\ref{lem:contracting_spiral_edges}).

\begin{figure}[ht]
		\centering
        \scalebox{0.95}{
            \tikzset{every picture/.style={line width=0.75pt}}      
            
            \begin{tikzpicture}[x=0.75pt,y=0.75pt,yscale=-1,xscale=1]
            
            \draw  [pattern= north east lines] (34.93,49.45) -- (34.93,130.14) -- (133.17,130.14) -- (133.17,49.45) -- (148.29,49.45) -- (148.29,144.81) -- (19.81,144.81) -- (19.81,49.45) -- cycle ;
            
            \draw  [pattern=north east lines] (223.78,49.61) -- (223.78,130.3) -- (349.04,130.3) -- (349.04,49.61) -- (365.55,49.61) -- (365.55,144.98) -- (208.67,144.98) -- (208.67,49.61) -- cycle ;
            
            \draw [color={rgb, 255:red, 0; green, 0; blue, 255 }  ,draw opacity=1 ]   (35.09,59.08) -- (84.57,82.06) ;
            
            \draw [color={rgb, 255:red, 0; green, 0; blue, 255 }  ,draw opacity=1 ]   (36.33,86.13) -- (84.57,82.06) ;
            
            \draw [color={rgb, 255:red, 0; green, 0; blue, 255 }  ,draw opacity=1 ]   (35.09,111.17) -- (84.57,82.06) ;
            
            \draw    (86.64,82.06) -- (133.35,58.35) ;
            
            \draw    (84.57,82.06) -- (133.35,80.36) ;
            
            \draw    (84.57,82.06) -- (133.35,106.77) ;
            
            \draw  [fill={rgb, 255:red, 0; green, 0; blue, 0 }  ,fill opacity=1 ] (82.51,82.06) .. controls (82.51,80.86) and (83.43,79.88) .. (84.57,79.88) .. controls (85.72,79.88) and (86.64,80.86) .. (86.64,82.06) .. controls (86.64,83.26) and (85.72,84.23) .. (84.57,84.23) .. controls (83.43,84.23) and (82.51,83.26) .. (82.51,82.06) -- cycle ;
            
            \draw    (59.83,70.57) -- (54.91,73.02) ;
            
            \draw    (57.38,65.69) -- (59.83,70.57) ;
            
            \draw    (59.83,96.61) -- (57.21,100.52) ;
            
            \draw    (55.63,95.15) -- (59.83,96.61) ;
            
            \draw    (56.25,81.16) -- (60.45,84.09) ;
            
            \draw    (60.45,84.09) -- (56.78,87.02) ;
            
            \draw    (223.94,63.16) -- (273.94,86.38) ;
            
            \draw    (223.94,85.9) -- (273.94,86.38) ;
            
            \draw    (223.94,109.37) -- (273.94,86.38) ;
            
            \draw    (313.32,86.38) -- (348.62,63.46) ;
            
            \draw    (313.32,86.38) -- (349.04,89.96) ;
            
            \draw    (313.32,86.38) -- (348.62,111.15) ;
            
            \draw    (273.94,86.38) -- (313.32,86.38) ;
            
            \draw  [fill={rgb, 255:red, 0; green, 0; blue, 0 }  ,fill opacity=1 ] (271.87,86.38) .. controls (271.87,85.18) and (272.8,84.2) .. (273.94,84.2) .. controls (275.08,84.2) and (276.01,85.18) .. (276.01,86.38) .. controls (276.01,87.57) and (275.08,88.55) .. (273.94,88.55) .. controls (272.8,88.55) and (271.87,87.57) .. (271.87,86.38) -- cycle ;
            
            \draw  [fill={rgb, 255:red, 0; green, 0; blue, 0 }  ,fill opacity=1 ] (311.25,86.38) .. controls (311.25,85.18) and (312.17,84.2) .. (313.32,84.2) .. controls (314.46,84.2) and (315.38,85.18) .. (315.38,86.38) .. controls (315.38,87.57) and (314.46,88.55) .. (313.32,88.55) .. controls (312.17,88.55) and (311.25,87.57) .. (311.25,86.38) -- cycle ;
            
            \draw    (288.64,82.38) -- (293.63,86.38) ;
            
            \draw    (293.63,86.38) -- (288.64,89.71) ;
            
            \draw    (249.56,75.33) -- (245.36,78.14) ;
            
            \draw    (246.15,71.83) -- (249.56,75.33) ;
            
            \draw    (249.56,97.34) -- (246.93,101.25) ;
            
            \draw    (245.36,95.87) -- (249.56,97.34) ;
            
            \draw    (245.36,83.4) -- (249.56,86.33) ;
            
            \draw    (249.56,86.33) -- (245.88,89.26) ;
            
            \draw  [pattern=north east lines] (423.73,48.7) -- (423.73,129.39) -- (548.99,129.39) -- (548.99,48.7) -- (565.5,48.7) -- (565.5,144.06) -- (408.61,144.06) -- (408.61,48.7) -- cycle ;
            
            \draw    (423.89,62.5) -- (451.5,85.82) ;
            
            \draw    (423.89,85.24) -- (451.5,85.82) ;
            
            \draw    (423.89,107.98) -- (451.5,85.82) ;
            
            \draw    (513.26,85.46) -- (548.57,62.55) ;
            
            \draw    (513.26,85.46) -- (548.99,89.04) ;
            
            \draw    (513.26,85.46) -- (548.57,110.23) ;
            
            \draw    (453.57,85.82) -- (513.26,85.46) ;
            
            \draw  [fill={rgb, 255:red, 0; green, 0; blue, 0 }  ,fill opacity=1 ] (449.43,85.82) .. controls (449.43,84.62) and (450.36,83.64) .. (451.5,83.64) .. controls (452.64,83.64) and (453.57,84.62) .. (453.57,85.82) .. controls (453.57,87.02) and (452.64,87.99) .. (451.5,87.99) .. controls (450.36,87.99) and (449.43,87.02) .. (449.43,85.82) -- cycle ;
            \draw  [fill={rgb, 255:red, 0; green, 0; blue, 0 }  ,fill opacity=1 ] (511.2,85.46) .. controls (511.2,84.26) and (512.12,83.29) .. (513.26,83.29) .. controls (514.41,83.29) and (515.33,84.26) .. (515.33,85.46) .. controls (515.33,86.66) and (514.41,87.63) .. (513.26,87.63) .. controls (512.12,87.63) and (511.2,86.66) .. (511.2,85.46) -- cycle ;
            \draw    (478.43,81.65) -- (483.42,85.64) ;
            \draw    (483.42,85.64) -- (478.43,88.98) ;
            
            \path[draw, <-, decorate, decoration ={snake, amplitude = 1.5}] (160,98) -- (194,98);
            
            \draw (136.78,24.51) node [anchor=north west][inner sep=0.75pt]    {$\Gamma $};
            
            \draw (42,55) node [anchor=north west][inner sep=0.75pt]    {$a$};
            \draw (42,75) node [anchor=north west][inner sep=0.75pt]    {$a$};
            \draw (42,93) node [anchor=north west][inner sep=0.75pt]    {$a$};

            \draw (81.23,58.49) node [anchor=north west][inner sep=0.75pt]    {$v$};
            
            \draw (40,30) node [anchor=north west][inner sep=0.75pt] {$\textcolor[rgb]{0,0,1}{H'}$};
            \draw (263.53,65.88) node [anchor=north west][inner sep=0.75pt]    {$w_{1}$};
            \draw (305.16,64.9) node [anchor=north west][inner sep=0.75pt]    {$w_{2}$};
            \draw (289.01,63.46) node [anchor=north west][inner sep=0.75pt]    {$e$};
            \draw (290.24,90.18) node [anchor=north west][inner sep=0.75pt]    {$0$};
            \draw (227,89.3) node [anchor=north west][inner sep=0.75pt]    {$a$};
            \draw (227,51.22) node [anchor=north west][inner sep=0.75pt]    {$a$};
            \draw (227,70.09) node [anchor=north west][inner sep=0.75pt]    {$a$};
            
            \draw (446.48,63.97) node [anchor=north west][inner sep=0.75pt]    {$w_{1}$};
            \draw (500.11,62.99) node [anchor=north west][inner sep=0.75pt]    {$w_{2}$};
            \draw (478.61,62.35) node [anchor=north west][inner sep=0.75pt]    {$e$};
            \draw (474.19,90.27) node [anchor=north west][inner sep=0.75pt]    {$a$};
            \draw (351.11,27.41) node [anchor=north west][inner sep=0.75pt]    {$\Gamma '$};
            \draw (554.4,25.83) node [anchor=north west][inner sep=0.75pt]    {$\Gamma '$};
            \draw (376,91.4) node [anchor=north west][inner sep=0.75pt]    {$=$};
            
            \end{tikzpicture}
            
        }

		\caption{Uncontracting a spiral edge of type $a$ at a spiral articulation point of ascent $a$. Here, $v$ has only one (spiral) cut component. Depicted is a $p$-cover in gain graph presentation, see \cref{def:gain}. The shaded area denotes an arbitrary connected graph and all edges without a specified gain have gain equal to $0$, including those in the shaded area.}
		\label{fig:spiral_uncontraction}
	\end{figure}
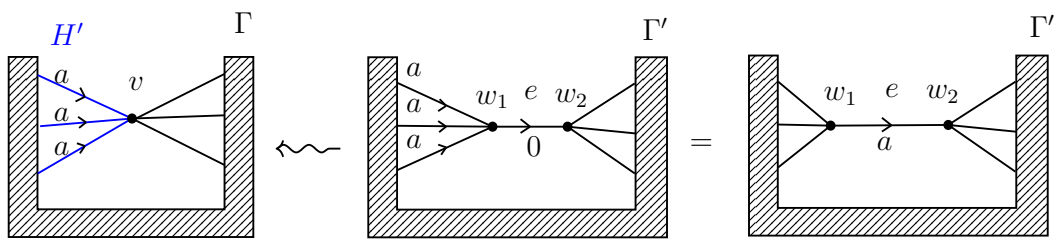

\begin{definition} \label{def:spiral_articulation}
	Let $\pi \colon \tilde \Gamma \to \Gamma$ be a free $p$-cover. 
    Let $v \in V(\Gamma)$ be a vertex and $\Gamma_0$ a cut-component at $v$. We say that $\Gamma_0$ is a \emph{spiral cut component of ascent $a \in \Z/p$} if there exists a non-empty set of half-edges $H' \subseteq T_v(\Gamma_0)$ with associated set of edges $E' = \big\{  \{h, \iota(h)\} \colon h \in H' \big\}$ such that the following hold:
	\begin{itemize}

        \item $\Gamma_0 \setminus E'$ is still connected and has trivial preimage, and
 
		\item for each $h \in H'$ with associated edge $e = \{h, \iota(h)\}$ and a (equivalently any) simple closed cycle $\gamma$ in $(\Gamma_0 \setminus E') \cup \{e\}$ passing through $e$ has ascent $a$ (with $\gamma$ oriented such that it enters $v$ through $h$).
	\end{itemize}
    We call $v$ a \emph{spiral articulation point} if there exists $a \in (\Z/p)^\times$ such that every cut component $\Gamma_i$ of $v$ is a spiral cut component of ascent $a$, see \cref{fig:spiral_uncontraction}.
\end{definition}

\begin{remark} \label{rem:complementyield-a}
    In \cref{def:spiral_articulation}, taking the complement of $H'$ in $T_{v}(\Gamma_0)$ instead of $H'$ makes $\Gamma_0$ a spiral cut component of ascent $-a$. Hence the ascent of a spiral cut component, and that of a spiral articulation point, is well defined up to sign.
\end{remark}

\begin{lemma} \label{lem:contracting_spiral_edges}
	Let $\pi' \colon \tilde \Gamma' \to \Gamma'$ be a $p$-cover and $e \in E(\Gamma')$ a spiral edge of type $a$. Then $\pi'$ is free and one of the following holds:
 
    \begin{enumerate}
        \item There is exactly one dilated vertex $v$ in $\pi' / e$ (which is the image of $e$) and every cut component of $\pi'/e$ at $v$ has trivial preimage.

        \item Or $\pi' / e $ is still free and every cut component of $\Gamma'/e$ at the image vertex of $e$ is either a spiral cut component of ascent $a$, or has trivial preimage.
 
    \end{enumerate}
\end{lemma}

\begin{proof}
	Performing an edge contraction on a dilated cover produces again a dilated cover. So, since $S_a$ is free, the same has to hold for $\pi'$. If $e$ is a loop, then contracting $e$ produces a dilated vertex. Since $e$ was a spiral edge, $\pi' \setminus \{e\}$ has trivial preimage and we are done. 
 
    Thus let us assume that $e = \{h_1, h_2\}$ is not a loop. Since $e$ is a spiral edge, it is clear that $\Gamma' \setminus \{e\}$ is connected and that $(\pi')^{-1}(\Gamma' \setminus \{e\})$ consists of $p$ disjoint copies of $\Gamma' \setminus \{e\}$. Let $v_i$ be the root vertex of $h_i$ for $i=1,2$ and define $H' \coloneqq T_{v_1}(\Gamma') \setminus \{h_1\}$. Then it is clear that any lift of a simple closed cycle $\gamma'$ in $\Gamma'$ that uses $e$ has to pass through exactly one half-edge in $H'$ and furthermore, the ascent of $\gamma'$ is necessarily equal to $a$. But now keeping the set $H'$, any cycle $\gamma$ in $\Gamma' / e$ as in Definition~\ref{def:spiral_articulation} is necessarily the image of some cycle $\gamma'$ in $\Gamma'$ and hence the ascent of $\gamma$ is $a$ as well. This shows that every cut component of $\Gamma'/e$ at $v$, which is the image of a subgraph of $\Gamma'$ containing both end vertices of $e$ is a spiral cut component of ascent $a$ in $\pi'/e$. Any cut components of $v_1$ and $v_2$ with trivial preimage in $\pi'$ are still cut components at $v$ with trivial preimage in $\pi'/e$. This proves the claim.
\end{proof}

Before we continue our study of spiral edges, we would like to introduce some helpful language, see \cite[Chapter~2]{GT01_TopGraphTheory}.

\begin{definition}\label{def:gain}
    A \emph{gain} or \emph{voltage graph} is a graph $\Gamma$ together with an orientation of its edges and an element $x_e \in \Z/p$ for each oriented edge $e$, which is called the \emph{gain}\footnote{We view the gain of an edge as the topological counterpart of dilation flow.} of $e$. The opposite edge is thought to have gain $-x_e$.   
\end{definition}

This is related to $p$-covers as follows. If $\pi : \tilde \Gamma \to \Gamma$ is a free $p$-cover, then we can represent it as a gain graph by identifying each vertex-fiber $\Z/p$-equivariantly with $\{0, 1, \ldots, p-1\}$ and choose an orientation for each edge in $\Gamma$. The gain of an edge $e$ is then the unique $x_e \in \Z/p$ such that the preimages of $e$ are precisely the edges $(i, i + x_e)$ for $i$ ranging over $\Z/p$. It is easy to see that $\pi$ can be reconstructed from the gain graph presentation\footnote{In the topological graph theory literature this is called the \emph{derived cover} of the gain graph, see \cite{GT01_TopGraphTheory}.}. We remark that the gain graph presentation of a free $p$-cover is not unique: indeed, changing the identification of a vertex fiber with $\Z/p$ will add a constant to the gain of every incident half-edge. This procedure is called \emph{switching}. To contract an edge in the gain graph presentation we switch to the presentation where that edge has gain $0$, then we contract the edge in question while keeping the gains on the other edges.

\begin{lemma} \label{lem:uncontracting_spirals}
	Let $\pi \colon \tilde \Gamma \to \Gamma$ be free and $v$ a spiral articulation point of ascent $a$. Then a spiral edge $e$ of type $a$ can be uncontracted at $v$ unless this uncontraction results in an unstable graph. Moreover, if $p$ is odd, this uncontraction is unique.
\end{lemma}

\begin{proof}
    A spiral edge can be uncontracted at $v$ by the following construction. Let $\Gamma_i$ be the spiral cut components of $\Gamma$ at $v$ and choose sets $H_i'$ as in \cref{def:spiral_articulation} such that all $\Gamma_i$ have the same ascent $a$. This might involve replacing some of the $H_i'$ with $T_{v}(\Gamma_i) \setminus H_i'$. Let $F \coloneqq \bigsqcup_i H_i'$ and define the uncontracted $p$-cover $\rho \colon \tilde \Gamma' \to \Gamma'$ where:
    \begin{enumerate}
        \item $\Gamma'$ is obtained from $\Gamma$ by replacing the vertex $v$ with a new edge $e=(w_1,w_2)$ and defining $w_1$ as the root vertex of all half-edges in $F$ and $w_2$ as the root vertex of the half-edges in $T_v(\Gamma) \setminus F$.

        \item We describe the cover $\rho$ based on a gain graph presentation of $\pi$: the gains of the old edges in $\Gamma$ are kept the same and the gain of the new edge $e$ is set to $0$. 
    \end{enumerate}
    The resulting $p$-cover $\rho$ is depicted in the middle of \cref{fig:spiral_uncontraction}. To show that the new edge $e$ is a spiral edge of type $a$, notice that by switching at the vertex $w_1$ we can rewrite $\rho$ as in the right side of \cref{fig:spiral_uncontraction}, i.e. with every edge having gain 0 except for $e$. In this picture it is now clear that contracting every edge except for $e$ produces the spiral $S_a$. 

    Now we assume $p$ is odd and we turn to proving uniqueness. We start by fixing a cut component $\Gamma_0$ of $v$. We show that the set $H'$ as in \cref{def:spiral_articulation} is unique up to complement in $T_{v}(\Gamma_0)$. This is clear if $\Gamma_0$ consists of a single loop edge, so assume this is not the case. Choose such a set $H'$ and let $E' = \big\{ \{h,\iota(h)\} \colon h \in H' \big\}$. By contracting all the edges not in $\Star_{\Gamma_0}(v)$ and choosing a gain graph presentation, we may reduce to the case where the cut component $\Gamma_0$ is the graph with two vertices $v,w$, the edges in $E'$ with gain $0$, and the edges in $\Star_{\Gamma_0}(v) \setminus E'$ which have gain $a$. Now let $H''$ be any other set satisfying \cref{def:spiral_articulation} and $E'' = \big\{ \{h,\iota(h)\} \colon h \in H'' \big\}$. Without loss of generality we may assume that there is a half-edge $h_0 \in H' \cap H''$ wih associated edge $e_0 = \{h_0, \iota(h_0)\}$. Then for any $h \neq h_0$ we have $h \in H'$ if and only if the cycle consisting of $\{h, \iota(h)\}$ and $e_0$ has trivial preimage. But the same holds for $H''$. Hence, up to taking the complement of $H''$ we have $H' = H''$.

    To finish the proof, we note that $\pi$ is free and hence by \cref{lem:contracting_spiral_edges} the only possible uncontraction of a spiral edge is given by the construction carried out at the begining of the proof with respect to some set $F \subseteq T_v(\Gamma)$. We have just shown that the sets $H_i' = F \cap T_v(\Gamma_i)$ are unique for any cut component $\Gamma_i$. From this it follows that the set $F$ is unique for the ascent $a$ (see also \cref{rem:complementyield-a} and note that for odd $p$ we have $a \neq -a$).
\end{proof}

\begin{figure}[ht]
    \centering
    \begin{tikzpicture}[scale = 1.3, decoration={
    markings,
    mark=
      at position 0.5
      with
      {
        \draw (-2pt,-2pt) -- (2pt,2pt);
        \draw (2pt,-2pt) -- (-2pt,2pt);
      }
    }]
        \vertex{0.5}{1}
        \vertex{0}{0}
        \vertex{1}{0}
        \vertex{0.5}{-1}
        \path[draw] (0.5, 1) to[bend left = 20] (1, 0) to[bend left = 20] (0.5, 1) -- (0, 0);
        \path[draw, postaction = decorate, color = blue] (0, 0) -- (1, 0);
        \path[draw] (1, 0) -- (0.5, -1) to[bend left = 20] (0,0) to[bend left = 20] (0.5, -1);
        
        \path[draw, ->, decorate, decoration ={snake, amplitude = 1.5}] (2,0) -- (3.5,0);

        \vertex{4}{1}
        \vertex{4}{0}
        \vertex{4}{-1}
        \path[draw] (4, 1) -- (4,0) to[bend right = 60] (4, 1);
        \path[draw] (4, 0) to[bend left = 60] (4, -1);
        \path[draw, postaction = decorate] (4, 1) to[bend right = 60] (4, 0); 
        \path[draw, postaction = decorate] (4, 0) -- (4, -1);
        \path[draw, postaction = decorate] (4, -1) to[bend left = 60] (4, 0);
        
        \node at (5, 0) {$=$};

        \vertex{6}{1}
        \vertex{6}{0}
        \vertex{6}{-1}
        \path[draw] (6, 1) -- (6,0) to[bend right = 60] (6, 1);
        \path[draw, postaction = decorate] (6, 0) to[bend left = 60] (6, -1);
        \path[draw, postaction = decorate] (6, 1) to[bend right = 60] (6, 0); 
        \path[draw] (6, 0) -- (6, -1) to[bend left = 60] (6, 0);

        \path[draw, ->, decorate, decoration ={snake, amplitude = 1.5}] (8.5,0) -- (7,0);

        \vertex{9.5}{1}
        \vertex{9}{0}
        \vertex{10}{0}
        \vertex{9.5}{-1}
        \path[draw] (9.5, 1) -- (9, 0) -- (9.5, -1);
        \path[draw, postaction = decorate, color = blue] (9, 0) -- (10, 0);
        \path[draw] (10, 0) to[bend left = 20] (9.5, 1) to[bend left = 20] (10, 0) to[bend left = 20] (9.5, -1) to[bend left = 20] (10,0);
    \end{tikzpicture}
    
    \caption{Two non-isomorphic uncontractions of a spiral edge (blue) in $p = 2$. Depicted are gain graph presentations, where edges with a cross are understood to have gain 1 (orientation does not matter in this case) and the remaining edges have gain 0.}
    \label{fig:scon_in_p_2}
\end{figure}
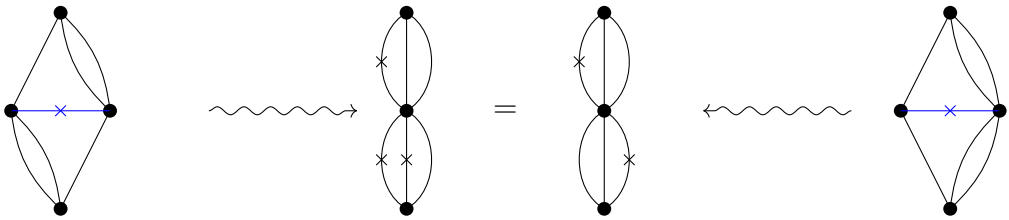

\begin{remark} \label{rem:scon_in_p_2}
The uniqueness statement in \cref{lem:uncontracting_spirals} does not hold for $p = 2$, see \cref{fig:scon_in_p_2} for a counter example.
\end{remark}

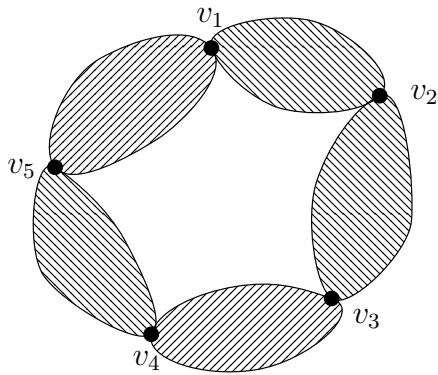
\begin{figure}[ht]
    \centering
    \scalebox{0.6}{         
            \tikzset{every picture/.style={line width=0.75pt}}     
            
            \begin{tikzpicture}[x=0.75pt,y=0.75pt,yscale=-1,xscale=1]

            \draw  [pattern=north east lines] (199.23,50.9) .. controls (229.47,33.8) and (275.23,20.9) .. (292,42) .. controls (308.77,63.1) and (268.47,104.8) .. (246.23,118.9) .. controls (224,133) and (176.77,160.1) .. (162,142) .. controls (147.23,123.9) and (169,68) .. (199.23,50.9) -- cycle ;
            
            \draw  [pattern=north west lines] (374.23,18.9) .. controls (405.47,25.8) and (449.77,65.1) .. (432,82) .. controls (414.23,98.9) and (368.47,98.8) .. (347.23,91.9) .. controls (326,85) and (298.77,61.1) .. (292,42) .. controls (285.23,22.9) and (343,12) .. (374.23,18.9) -- cycle ;
            
            \draw  [pattern=north west lines] (152,232) .. controls (137.23,208.9) and (143.77,130.1) .. (162,142) .. controls (180.23,153.9) and (206.23,177.9) .. (219.23,201.9) .. controls (232.23,225.9) and (255.23,271.9) .. (242,282) .. controls (228.77,292.1) and (166.77,255.1) .. (152,232) -- cycle ;

            \draw  [pattern=north west lines] (379.23,154.9) .. controls (388.23,125.9) and (410.23,91.9) .. (432,82) .. controls (453.77,72.1) and (461.23,160.9) .. (458.23,187.9) .. controls (455.23,214.9) and (406.23,262.9) .. (392,252) .. controls (377.77,241.1) and (370.23,183.9) .. (379.23,154.9) -- cycle ;

            \draw  [pattern=north east lines] (312,242) .. controls (343.23,236.9) and (363.23,241.9) .. (392,252) .. controls (420.77,262.1) and (373.23,304.9) .. (328.23,311.9) .. controls (283.23,318.9) and (236.77,302.1) .. (242,282) .. controls (247.23,261.9) and (280.77,247.1) .. (312,242) -- cycle ;

            \draw  [fill={rgb, 255:red, 0; green, 0; blue, 0 }  ,fill opacity=1 ] (285.95,42) .. controls (285.95,38.66) and (288.66,35.95) .. (292,35.95) .. controls (295.34,35.95) and (298.05,38.66) .. (298.05,42) .. controls (298.05,45.34) and (295.34,48.05) .. (292,48.05) .. controls (288.66,48.05) and (285.95,45.34) .. (285.95,42) -- cycle ;

            \draw  [fill={rgb, 255:red, 0; green, 0; blue, 0 }  ,fill opacity=1 ] (155.95,142) .. controls (155.95,138.66) and (158.66,135.95) .. (162,135.95) .. controls (165.34,135.95) and (168.05,138.66) .. (168.05,142) .. controls (168.05,145.34) and (165.34,148.05) .. (162,148.05) .. controls (158.66,148.05) and (155.95,145.34) .. (155.95,142) -- cycle ;

            \draw  [fill={rgb, 255:red, 0; green, 0; blue, 0 }  ,fill opacity=1 ] (235.95,282) .. controls (235.95,278.66) and (238.66,275.95) .. (242,275.95) .. controls (245.34,275.95) and (248.05,278.66) .. (248.05,282) .. controls (248.05,285.34) and (245.34,288.05) .. (242,288.05) .. controls (238.66,288.05) and (235.95,285.34) .. (235.95,282) -- cycle ;

            \draw  [fill={rgb, 255:red, 0; green, 0; blue, 0 }  ,fill opacity=1 ] (385.95,252) .. controls (385.95,248.66) and (388.66,245.95) .. (392,245.95) .. controls (395.34,245.95) and (398.05,248.66) .. (398.05,252) .. controls (398.05,255.34) and (395.34,258.05) .. (392,258.05) .. controls (388.66,258.05) and (385.95,255.34) .. (385.95,252) -- cycle ;

            \draw  [fill={rgb, 255:red, 0; green, 0; blue, 0 }  ,fill opacity=1 ] (425.95,82) .. controls (425.95,78.66) and (428.66,75.95) .. (432,75.95) .. controls (435.34,75.95) and (438.05,78.66) .. (438.05,82) .. controls (438.05,85.34) and (435.34,88.05) .. (432,88.05) .. controls (428.66,88.05) and (425.95,85.34) .. (425.95,82) -- cycle ;
            
            \draw (277,6) node [anchor=north west][inner sep=0.75pt] [xscale=1.7, yscale=1.7]    {$v_{1}$};
            \draw (455,70.4) node [anchor=north west][inner sep=0.75pt]    [xscale=1.7, yscale=1.7]    {$v_{2}$};
            \draw (407,258.4) node [anchor=north west][inner sep=0.75pt]    [xscale=1.7, yscale=1.7]    {$v_{3}$};
            \draw (224,293.4) node [anchor=north west][inner sep=0.75pt]    [xscale=1.7, yscale=1.7]    {$v_{4}$};
            \draw (120,132.4) node [anchor=north west][inner sep=0.75pt]    [xscale=1.7, yscale=1.7]    {$v_{5}$};
            
            \end{tikzpicture}
            
    }
    \caption{Cartoon of the shape of a $p$-cover with $5$ spiral articulation points which does not lie in the equivariant bridge locus. The shaded areas denote arbitrary connected graphs with trivial preimage.}
    \label{fig:multiple_spiral_art}
\end{figure}

If a $p$-cover has multiple spiral articulation points, then this imposes severe restrictions on the geometry of the cover, see \cref{fig:multiple_spiral_art}. As a consequence, we find the following lemma.

\begin{lemma} \label{lem:ascent_well_def}
    Let $\pi \colon \tilde \Gamma \to \Gamma$ be a $p$-cover. Then all spiral articulation points in $\pi$ have the same ascent.
\end{lemma}
\begin{proof}
    Let $v_1,v_2$ be two spiral articulation points in $\Gamma$. Denote by $\rho \colon \tilde \Gamma' \to \Gamma'$ the $p$-cover resulting from uncontracting spiral edges $e_1$ and $e_2$ at $v_1$ and $v_2$ and then contracting every edge except for $e_1, e_2$. Assume for a contradiction that one of the $e_i$ (and hence both of them) are loops. The $e_i$ are spiral edges, so both have the property that contracting the other edge a cover of the form $S_a$. But that means that the other edge has to have ascent zero which is a contradiction. We conclude that $\Gamma'$ in fact consists of two vertices $w_1,w_2$ with the two spiral edges $e_1, e_2$ joining them. The $p$-cover $\rho$ has an automorphism exchanging the edges $e_1, e_2$ and the vertices $w_1, w_2$. This implies that whatever the ascent of $e_1$ after contracting $e_2$ may be, it is the same as the ascent for $e_2$ after contracting $e_1$. 
\end{proof}

\begin{proof}[Proof of \cref{lem:properties_are_good2}]
	Recall that $(\MS)_P$ is the equivariant bridge locus $\MS^{br}$ and that by Lemma~\ref{lem:bridge_articulation_points} a $p$-cover $\pi$ is in $P^\ast(\MS)$ if and only if $\pi$ possesses a bridge articulation point.
	\medskip

    	\textbf{Check Condition~\ref{condition1}.}
	Let $\pi \colon \tilde \Gamma \to \Gamma$ be a $p$-cover and $e \in \Gamma$ a spiral edge. In particular, $\pi$ is a free cover. Assume that $\pi/e \in P^\ast(\MS)$. We want to show that this implies $\pi \in P^\ast(\MS)$. 
    Consider the possible shapes of $\pi / e$ as worked out in \cref{lem:contracting_spiral_edges}.
    In the first case we see that $\pi / e$ has a $(g-1)$-bridge articulation point (the image of $e$). But in this case, $e$ was a loop to begin with and therefore, the root vertex of $e$ already was a $(g-1)$-bridge articulation point in $\pi$. So $\pi \in P^\ast(\MS)$.

    In the second case of \cref{lem:contracting_spiral_edges}, the contracted cover $\pi / e$ is in $P^\ast(\MS)$ if and only if there is at least one cut component $\Gamma_0$ of some vertex $v$ with trivial preimage. Without loss of generality, $v$ is the image of $e$ under the edge contraction, because otherwise we are already done. Let $\Gamma_1 \subseteq \Gamma$ be the preimage of $\Gamma_0$ under the contraction of $e$ (so in particular, $\Gamma_1$ contains $e$). Then $\pi^{-1}(\Gamma_1)$ is disconnected if and only if $e$ is a bridge in $\Gamma_1$. On the other hand, $(\pi/e)^{-1}(\Gamma_0) = \pi^{-1}(\Gamma_1) / \pi^{-1}(e)$ is disconnected. This means that in fact $e$ must be a bridge in $\Gamma_1$, and hence $\Gamma_1 \setminus \{e\}$ is simply a copy of $\Gamma_0$ rooted at one of the ends of $e$. But this means that $\Gamma_0$ is a cut component with trivial preimage in $\pi$ as well. This implies that $\pi$ has a bridge articulation point (see \cref{def:bridge_articulation}), hence $\pi \in P^\ast(\MS)$.

    \medskip
 
	\textbf{Check Condition~\ref{condition2}.}
	Let $\pi \colon \tilde \Gamma \to \Gamma$ be a $p$-cover which is not in $P^\ast(\MS)$. We have to show that $\pi$ admits a unique maximal uncontraction $\rho$ by spiral edges (of any type) and that the corresponding face inclusion $\pi \preccurlyeq \rho$ is canonical.

    Combining Lemmas~\ref{lem:contracting_spiral_edges} and \ref{lem:uncontracting_spirals} and the assumption that $\pi$ is not in the bridge locus, we see that spiral edges can only be uncontracted at spiral articulation points. By Lemma~\ref{lem:uncontracting_spirals} this is in indeed possible, unless it breaks stability. By Lemma~\ref{lem:ascent_well_def} all of these possible uncontractions produce spiral edges of the same type $a$. Note that uncontracting a spiral edge at any vertex does not affect any other spiral articulation points in the $p$-cover. Summing up: by going over all the spiral articulation points and uncontracting wherever possible, we obtain a unique maximal uncontraction $\rho$ or $\pi$ by spiral edges.
	
	The corresponding face inclusion $\pi \preccurlyeq \rho$ is canonical by the same argument as before: any automorphism of $\pi$ would permute the spiral articulation points and therefore it can be lifted to $\rho$.

    \medskip

	\textbf{Check Condition~\ref{condition3}.}
    A simplex in $\MS$ with all vertices in $Q$ corresponds to a stable $p$-cover all of whose edges are spiral edges (of whatever type). 
    From \cref{lem:ascent_well_def} we can see that the only such covers are in fact free $p$-covers of a cycle graph with positive genus at every vertex. All of these are in the equivariant bridge locus $(\MS)_P$ since the fiber over any vertex has enough genus to allow for the uncontraction of an equivariant bridge.

    \medskip
    
    This finishes the proof of the main part of the statement. The second assertion is obtained by virtue of \cite[Proposition~4.11]{CGP2}.
\end{proof}

\subsection{Proof of \cref{thm:contractible_loci}}
\label{subsec:proofOfThm}

    By \cref{cor:loop_weight} the loci $\MS^{w}$ and $\MS^{lw}$ are contractible. Combining this with \cref{lem:properties_are_good} we deduce that $\MS^{br}$ is contractible as well. Note that so far no assumption on the prime $p$ has been made. 
    For $p > 2$, \cref{lem:properties_are_good2} shows that the sparsely connected locus $\MS^{scon}$ deformation retracts onto $\MS^{br}$ and is therefore contractible as well. 

    Finally, recall that the equivariant parallel edge locus $\MS^{par}$ is defined as the closure of the union of $\MS^{scon}$ and
	\begin{equation*}
		\left\{
		\pi \colon \tilde \Gamma \to \Gamma \mathrel{\Bigg |} \ 
		\begin{minipage}{0.75\textwidth}
			there exists a pair of parallel edges $e_1, e_2 \in \Gamma$ such that
			$|\pi^{-1}(e_1)| = |\pi^{-1}(e_2)| = p$ and every lift of the simple closed cycle defined by $e_1$ and $e_2$ to $\tilde \Gamma$ is again a simple closed cycle
		\end{minipage}
		\ \right\}.
	\end{equation*}    
    For this locus, we use the following observation from \cite[Sections~5 and~6]{AllcockCoreyPayne}. Suppose we have a topological space $X$ defined by gluing a quotient of a simplex $\Delta^n / G$ for a subgroup $G \leq S_{n+1}$ along its surface $\partial(\Delta^n / G) = \partial(\Delta^n) / G$ to some topological space $\tilde{X}$. If $G$ contains a reflection, then $\partial (\Delta^n/G)$ is contractible and the glued space $X$ is the mapping cone of $\partial \Delta^n / G \hookrightarrow \tilde{X}$. So if $\tilde{X}$ was contractible, then the same is true for $X$. 
    
    Now let $\pi \colon \tilde \Gamma \to \Gamma$ be a $p$-cover in $\MS^{par}$. If $\pi$ has equivariant parallel edges, i.e. a pair of parallel edges $e_1$, $e_2$ such that both $e_i$ have $p$ preimages and every lift of the simple closed cycle defined by $e_1$ and $e_2$ is a simple closed cycle, then we notice two things. First, $\pi$ has an automorphism (given by exchanging $e_1$ and $e_2$) which acts on the corresponding simplex of the moduli space as a reflection. Second, if we contract any arbitrary edge in $\pi$, then the result will have equivariant parallel edges or equivariant loops. More generally, this shows that for any $\pi$ in $\MS^{par}$, any edge contraction will yield a $p$-cover in $\MS^{scon}$. Combined, these observations show that $\MS^{par}$ is an iterated mapping cone over $\MS^{scon}$. For $p > 2$ the latter is contractible, and hence so is $\MS^{par}$.                    $\hfill \square$

\begin{remark} \label{rem:moreLoci}
    We currently do not believe that our findings in \cref{thm:contractible_loci} can be extended any further for large $p$. Indeed, looking at the list of vertices in $\MS$ given in \cref{tab:vertices}, the stars of the ring vertices, the butterfly, and the parallel bridges all contribute to the non-trivial topology of $\Delta_{2, 5}$, see \cref{sec:genus_2} below. However, we leave the following open question:

    \begin{question}
        For $g \geq 2$ and $p = 2$, is the locus generated by the parallel bridges (\cref{tab:vertices}) contractible? And if so, can this be combined with the contractibility of $\Delta_{g, 2}^{br}$?
    \end{question}
\end{remark}

\section{The moduli space \texorpdfstring{$\MS$}{} is simply connected}
\label{sec:simply_connected}

In this section we prove that the geometric realization of $\MS$ is simply connected for all $g \geq 3$ and all primes $p$. This constitutes most of the proof of \cref{thm:simply_connected}. The remaining case $g = 2$  will follow from the explicit computations in the \cref{sec:genus_2}. The key for our discussion here is \cite[Theorem~3.1]{AllcockCoreyPayne} which states that the fundamental group of a symmetric $\Delta$-complex $X$ is generated by the fundamental group of its 1-skeleton $X^{(1)}$. More precisely, let $x \in |X^{(1)}|$ be a base point. Then there is a natural surjection 
\[ \pi_1\big( |X^{(1)}|, x \big) \longrightarrow \pi_1\big( |X|, x \big). \]
If $\gamma \subseteq \big| \MS^{(1)} \big|$ is a simple closed cycle which is fully contained in the equivariant bridge locus, then by Theorem~\ref{thm:contractible_loci} it is clear that $\gamma$ is contractible in $\big| \MS \big|$. Therefore, a cycle $\gamma$ can only be non-trivial in $\pi_1\big( |\MS|, x \big)$ if is is partly supported on the geometric realization of 1-simplices which are not contained in the equivariant bridge locus. In order to give a complete list of the combinatorial types of $p$-covers corresponding to these 1-simplices, we first state a handy criterion to decide if a given cover is in the equivariant weight locus.

\begin{lemma} \label{lem:weight_locus_criterion}
	A $p$-cover $\tilde \Gamma \to \Gamma$ is in the equivariant weight locus if and only if one of the following holds:
	\begin{enumerate}
		\item There is a free vertex $v \in \Gamma$ with $g(v) \geq 1$.
		\item There is a dilated vertex $v \in \Gamma$ with $d$ incident dilated half-edges such that:
            \setlist[enumerate]{wide=40pt}
	   		\begin{enumerate}
				\item $d \geq 2\Big(\frac{p}{p+1} + 1\Big)$, or
				\item $d \geq 2$ and $g(v) \geq 1$, or
				\item $g(v) \geq 2$.
			\end{enumerate}
	\end{enumerate}
\end{lemma}

\begin{proof}
	If $v \in \Gamma$ is free then each of the $p$ preimages of $v$ has genus equal to $g(v)$. If $v$ is dilated then by the local Riemann-Hurwitz condition~\eqref{eq:local-RH} the genus of the preimage $\tilde v$ of $v$ is
	\begin{equation} \label{eq:genus_cover_vertex}
		g(\tilde v) =  p\big( g(v) - 1 \big) + 1 + \frac{1}{2}d(p-1).
	\end{equation}
	The claim follows by a case-by-case analysis which is left to the reader.
\end{proof}

\begin{proposition}
	Let $p$ be a prime number and suppose that $g\geq{4}$. Then $\big| \MS \big|$ is simply connected.
\end{proposition}

\begin{proof}
	 Suppose that $g\geq{4}$. The first Betti number of the target graph $\Gamma$ is $\beta = |E(\Gamma)| - |V(\Gamma)| + 1$. If the cover corresponds to a $1$-simplex in $\MS$, then $|E(\Gamma)|=2$, so that $\beta = 3 - |V(\Gamma)|$. Going through all possible values for $|V(\Gamma)|$, we see that in every case there must be a vertex $v$ with $g(v) \geq 2$. Consequently, in every case we find the cover to lie in the equivariant weight locus by \cref{lem:weight_locus_criterion}, so we see that all $1$-simplices lie in the contractible part of $\MS$. Hence the moduli space is simply connected.
\end{proof}

\begin{proposition} 
Let $p$ be a prime number. Then $\big| \Delta_{3,p} \big|$ is simply connected. 
\end{proposition}

\begin{proof}
	Let $\pi \colon \tilde \Gamma \to \Gamma$ be a $p$-cover such that $g(\Gamma) = 3$  and $\Gamma$ has exactly 2 edges. We first show that the only such covers which are not contained in the equivariant bridge locus are
 
    \begin{equation}\label{eq:genus3GraphA}
        \vcenter{\hbox{
		\begin{tikzpicture}
			\draw (0.5, 1.5) ellipse (0.5 and 0.1);
			\draw (0.5, 1.5) ellipse (0.5 and 0.25);
			\draw (1.5, 1.5) ellipse (0.5 and 0.1);
			\draw (1.5, 1.5) ellipse (0.5 and 0.25);
			\vertex{0}{1.5}
			\vertex{1}{1.5}
			\vertex{2}{1.5}
			
			\path[draw, ->] (1, 1.1) -- (1, 0.6);
			
			\draw (0,0) -- (2,0);
			\vertex[1]{0}{0}
			\vertex[1]{1}{0}
			\vertex[1]{2}{0}
		\end{tikzpicture} 
        }}
        \quad \text{ and } \quad
        \vcenter{\hbox{
        \begin{tikzpicture}
            \draw (0.5, 2.5) ellipse (0.5 and 0.1);
            \draw (0.5, 2.5) ellipse (0.5 and 0.25);
            \draw (0.5, 2.5) ellipse (0.5 and 0.4);
            \vertex{0}{2.5}
            \vertex{1}{2.5}
            
            \path[draw, ->] (0.5, 2) -- (0.5, 1.5);
            
            \draw (0.5, 1) ellipse (0.5 and 0.25);
            \vertex[1]{0}{1}
            \vertex[1]{1}{1}
        \end{tikzpicture}
        }}
        \text{.}
        \vspace{2mm}
    \end{equation}
    
    \noindent
    To this end, note first that $\Gamma$ has to be one of the graphs depicted in \cref{fig:1-skeleton_targets}. For each of these we now discuss the possible covers which are not in the equivariant bridge locus.

    \begin{table}[ht]
        \centering
        \begin{tabular}{|c|c|c|c|c|}
            \hline
             $\Gamma_1$ & $\Gamma_2$ & $\Gamma_3$ & $\Gamma_4$ & $\Gamma_5$ \\
             \hline
    		\begin{tikzpicture}
    			\draw (0,0) -- (1,0);
    			\vertex[2]{0}{0}
    			\vertex[0]{1}{0}
    			\selfloop{1}{0}
    		\end{tikzpicture}&
    		\begin{tikzpicture}
    			\draw (0,0) -- (1,0);
    			\vertex[1]{0}{0}
    			\vertex[1]{1}{0}
    			\selfloop{1}{0}
    		\end{tikzpicture}&
    		\begin{tikzpicture}
    			\vertex[1]{0}{0}
    			\selfloop{0}{0}
    			\selfloopleft{0}{0}
    		\end{tikzpicture}&
    		\begin{tikzpicture}
    			\draw (0,1) -- (2,1);
    			\vertex[1]{0}{1}
    			\vertex[1]{1}{1}
                \vertex[1]{2}{1}
    			{\color{white} \vertex{2}{0}}
    		\end{tikzpicture}& 
    		\begin{tikzpicture}
    			\draw (0.5, 1) ellipse (0.5 and 0.25);
    			\vertex[1]{0}{1}
    			\vertex[1]{1}{1}
                {\color{white} \vertex{1}{0}}
    		\end{tikzpicture}\\
            \hline
        \end{tabular}
        \caption{All stable weighted graphs with two edges and genus 3.}
    		\label{fig:1-skeleton_targets}
    \end{table}

	\begin{itemize}
		\item {Target graph $\Gamma_1$:} The genus 2 vertex automatically places any cover of $\Gamma_1$ in the equivariant weight locus.
		
		\item {Target graph $\Gamma_2$:} By Lemma~\ref{lem:weight_locus_criterion} both vertices need to be dilated and the loop needs to be free. But then the cover has an equivariant loop and thus lives in the (contractible) equivariant weight or loop locus.
		
		\item {Target graph $\Gamma_3$:} Same argument: in order to avoid the equivariant weight locus, the vertex needs to be dilated and the loops need to be free. But then these are equivariant loops.

        \item {Target graph $\Gamma_4$:} The vertices all have to be dilated or else the cover would be in the equivariant weight locus by Lemma~\ref{lem:weight_locus_criterion}. The edges all have to be free because they are bridges. This leaves the first cover in \eqref{eq:genus3GraphA} as only option.
  
		\item {Target graph $\Gamma_5$:} Again by Lemma~\ref{lem:weight_locus_criterion} both vertices need to be dilated and at least one edge needs to be free. But then both edges need to be free or else there is no dilation flow. It remains the second cover in \eqref{eq:genus3GraphA}.
	\end{itemize}
	This shows the first claim. To show that the covers in \eqref{eq:genus3GraphA} do not support any non-trivial cycle we note the following general fact. If a combinatorial type $\pi \colon \tilde \Gamma \to \Gamma$ of $X^{(1)}$ has an automorphism which exchanges the two edges of $\Gamma$, then the 1-simplex $\sigma_\pi$ in $X^{(1)}$ corresponding to $\pi$ is folded in half and has only one zero-dimensional face. Therefore, $\sigma_\pi$ is not contained in the support of any simple closed cycle in $X^{(1)}$. This applies to both covers in $\eqref{eq:genus3GraphA}$ and we are done.
\end{proof}

\section{Homotopy type in genus 2}  \label{sec:genus_2}

\begin{figure}[ht]
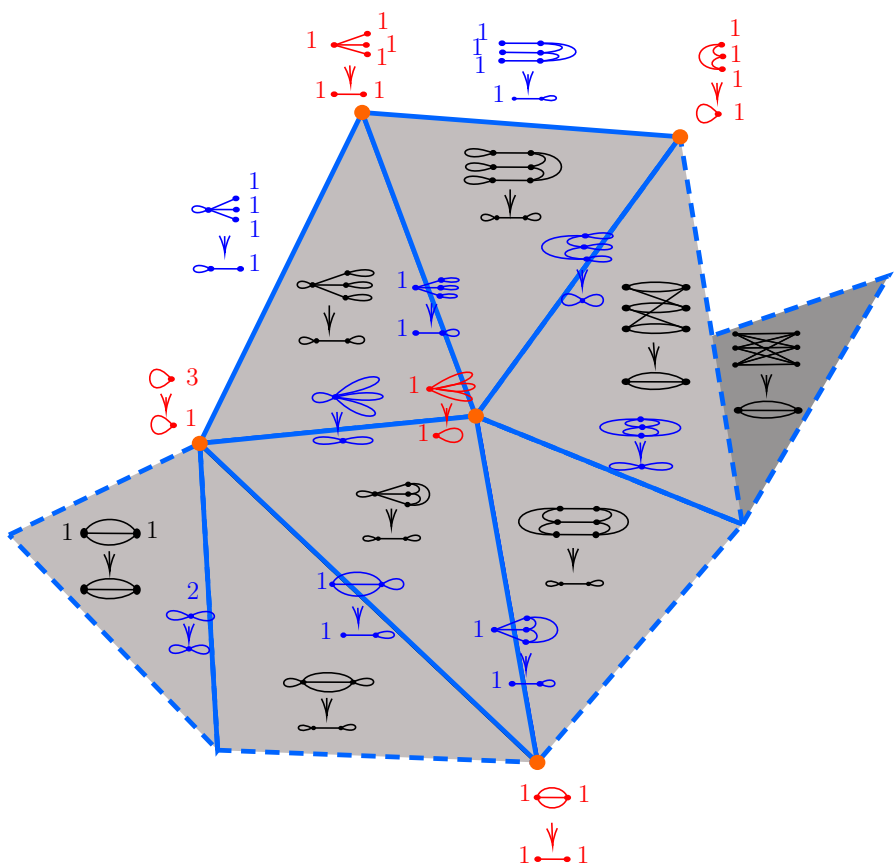

    \centering
    \scalebox{0.8}{  
    
            \tikzset{every picture/.style={line width=0.75pt}}       
            

            
    }
    \caption{The moduli space $\Delta_{2,3}$. The same explanations as in \cref{fig:Delta_2_2} apply. Again we do not draw dilation flows: it turns out that every $3$-cover in the picture has up to isomorphism a unique dilation flow.}
    \label{fig:Delta_2_3}
\end{figure}

\begin{figure}[ht]
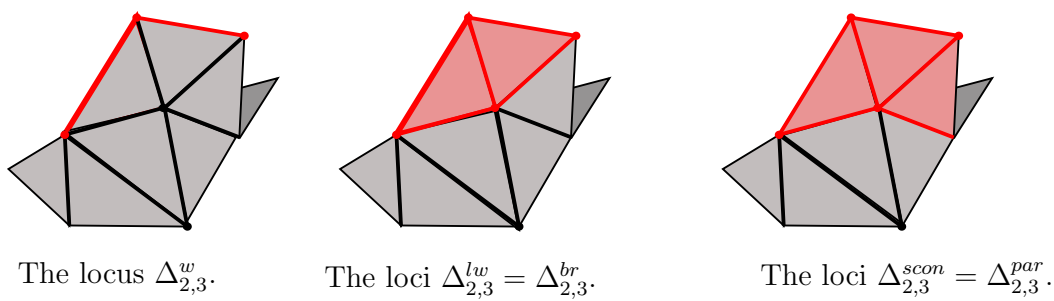

    \centering
    \scalebox{0.95}{
        \tikzset{every picture/.style={line width=0.75pt}} 
            
            %
}
    \caption{The different loci discussed in \cref{sec:contractible_loci} highlighted in the case $g=2$ and $p =3$.}
    \label{fig:loci}
\end{figure}

In this section, we prove \cref{thm:MainThmGenus2} and thereby conclude the proof of \cref{thm:simply_connected}. We start by verifying the claim for $p = 2,3$. This can be easily done by constructing the moduli spaces explicitly by gluing the cells they consist of, see Figures~\ref{fig:Delta_2_2} and \ref{fig:Delta_2_3}. We use this opportunity to highlight the various contractible loci discussed in \cref{sec:contractible_loci}, see \cref{fig:loci}. However, from the pictures it is clear that $\Delta_{2,2}$ and $\Delta_{2, 3}$ are contractible. From now on we work in $p \geq 5$.

In order to determine the homotopy type of $\Delta_{2, p}$, we are going to list all $p$-covers corresponding to the maximal cells of $\Delta_{2,p}$. For this and the remainder of the section we use the following convention for drawing these. This is a refinement of the gain graph presentations that we used in \cref{subsec:scon}. 

\begin{convention}\label{convention}
    We depict a $p$-cover $\pi \colon \tilde \Gamma \to \Gamma$ by simply drawing the target graph $\Gamma$ with the following decorations:
    \begin{enumerate}
        \item the dilation subgraph $\Gamma_\dil$ together with the dilation flow with respect to some orientation in bold, 
        \item the free part of $\pi$ in gain graph presentation, i.e. with gains in $\Z/p$ with respect to some orientation indicated in thin print. In case one or both endpoints of a free edge are dilated, there is no gain to be specified. 
    \end{enumerate}
    \cref{tab:p-decoratedGraph} gives a few examples of decorated graphs and the $p$-covers they represent.
\end{convention}

\begin{table}[ht]
    \centering
    \begin{tabular}{|c|c|c|}     
    \hline
    \scalebox{0.6}{
        \tikzset{every picture/.style={line width=0.75pt}}       
        
        \begin{tikzpicture}[x=0.75pt,y=0.75pt,yscale=-1,xscale=1]
        \draw    (280.93,121.96) -- (337.72,121.96) ;
        \draw    (337.51,121.96) .. controls (374.03,93.02) and (376.87,151.27) .. (337.72,121.96) ;
        \draw    (280.93,121.96) .. controls (242.83,95.29) and (244.53,150.13) .. (281.14,121.96) ;
        \draw  [fill={rgb, 255:red, 0; green, 0; blue, 0 }  ,fill opacity=1 ] (278.02,121.96) .. controls (278.02,120.81) and (279.42,119.88) .. (281.14,119.88) .. controls (282.87,119.88) and (284.27,120.81) .. (284.27,121.96) .. controls (284.27,123.1) and (282.87,124.04) .. (281.14,124.04) .. controls (279.42,124.04) and (278.02,123.1) .. (278.02,121.96) -- cycle ;
        \draw  [fill={rgb, 255:red, 0; green, 0; blue, 0 }  ,fill opacity=1 ] (334.6,121.96) .. controls (334.6,120.81) and (336,119.88) .. (337.72,119.88) .. controls (339.45,119.88) and (340.85,120.81) .. (340.85,121.96) .. controls (340.85,123.1) and (339.45,124.04) .. (337.72,124.04) .. controls (336,124.04) and (334.6,123.1) .. (334.6,121.96) -- cycle ;
        \draw    (276.91,47.84) -- (334.19,76.96) ;
        \draw    (334.19,76.96) .. controls (370.71,48.02) and (373.55,106.27) .. (334.41,76.96) ;
        \draw  [fill={rgb, 255:red, 0; green, 0; blue, 0 }  ,fill opacity=1 ] (274.01,47.84) .. controls (274.01,46.69) and (275.4,45.76) .. (277.13,45.76) .. controls (278.85,45.76) and (280.25,46.69) .. (280.25,47.84) .. controls (280.25,48.99) and (278.85,49.92) .. (277.13,49.92) .. controls (275.4,49.92) and (274.01,48.99) .. (274.01,47.84) -- cycle ;
        \draw  [fill={rgb, 255:red, 0; green, 0; blue, 0 }  ,fill opacity=1 ] (331.07,47.92) .. controls (331.07,46.77) and (332.47,45.84) .. (334.19,45.84) .. controls (335.92,45.84) and (337.32,46.77) .. (337.32,47.92) .. controls (337.32,49.07) and (335.92,50) .. (334.19,50) .. controls (332.47,50) and (331.07,49.07) .. (331.07,47.92) -- cycle ;
        \draw    (277.13,47.84) -- (334.19,20.23) ;
        \draw  [fill={rgb, 255:red, 0; green, 0; blue, 0 }  ,fill opacity=1 ] (331.07,20.23) .. controls (331.07,19.08) and (332.47,18.15) .. (334.19,18.15) .. controls (335.92,18.15) and (337.32,19.08) .. (337.32,20.23) .. controls (337.32,21.38) and (335.92,22.31) .. (334.19,22.31) .. controls (332.47,22.31) and (331.07,21.38) .. (331.07,20.23) -- cycle ;
        \draw    (333.98,20.23) .. controls (371.77,-8.53) and (372.34,49.71) .. (334.19,20.23) ;
        \draw    (310.34,76.96) -- (310.34,103.41) ;
        \draw [shift={(310.34,105.41)}, rotate = 270] [color={rgb, 255:red, 0; green, 0; blue, 0 }  ][line width=0.75]    (10.93,-3.29) .. controls (6.95,-1.4) and (3.31,-0.3) .. (0,0) .. controls (3.31,0.3) and (6.95,1.4) .. (10.93,3.29)   ;
        \draw    (277.13,47.84) -- (334.19,47.92) ;
        \draw  [fill={rgb, 255:red, 0; green, 0; blue, 0 }  ,fill opacity=1 ] (331.07,76.96) .. controls (331.07,75.81) and (332.47,74.88) .. (334.19,74.88) .. controls (335.92,74.88) and (337.32,75.81) .. (337.32,76.96) .. controls (337.32,78.11) and (335.92,79.04) .. (334.19,79.04) .. controls (332.47,79.04) and (331.07,78.11) .. (331.07,76.96) -- cycle ;
        \draw    (334.19,47.92) .. controls (370.71,18.99) and (373.55,77.23) .. (334.41,47.92) ;
        \draw    (276.93,48.17) .. controls (238.83,21.5) and (240.53,76.34) .. (277.14,48.17) ;
        \draw  [color={rgb, 255:red, 0; green, 0; blue, 0 }  ,draw opacity=1 ][fill={rgb, 255:red, 0; green, 0; blue, 0 }  ,fill opacity=1 ] (254.28,118.59) -- (258.12,126.22) -- (250.68,126.34) -- cycle ;
        \draw (236.5,114.07) node [anchor=north west][inner sep=0.75pt]    {$1$};
        \end{tikzpicture}
        }
        &
        \scalebox{0.6}{
            \tikzset{every picture/.style={line width=0.75pt}} 
            \begin{tikzpicture}[x=0.75pt,y=0.75pt,yscale=-1,xscale=1]
            \draw    (271.93,127.59) -- (328.72,127.59) ;
            \draw    (328.51,127.59) .. controls (365.03,98.65) and (367.87,156.9) .. (328.72,127.59) ;
            \draw    (271.93,127.59) .. controls (233.83,100.92) and (235.53,155.76) .. (272.14,127.59) ;
            \draw  [fill={rgb, 255:red, 0; green, 0; blue, 0 }  ,fill opacity=1 ] (269.02,127.59) .. controls (269.02,126.44) and (270.42,125.51) .. (272.14,125.51) .. controls (273.87,125.51) and (275.27,126.44) .. (275.27,127.59) .. controls (275.27,128.74) and (273.87,129.67) .. (272.14,129.67) .. controls (270.42,129.67) and (269.02,128.74) .. (269.02,127.59) -- cycle ;
            \draw  [fill={rgb, 255:red, 0; green, 0; blue, 0 }  ,fill opacity=1 ] (325.6,127.59) .. controls (325.6,126.44) and (327,125.51) .. (328.72,125.51) .. controls (330.45,125.51) and (331.85,126.44) .. (331.85,127.59) .. controls (331.85,128.74) and (330.45,129.67) .. (328.72,129.67) .. controls (327,129.67) and (325.6,128.74) .. (325.6,127.59) -- cycle ;
            \draw    (328.26,84.67) .. controls (364.78,55.74) and (367.62,113.98) .. (328.47,84.67) ;
            \draw    (271.19,56.55) .. controls (246.5,61.38) and (247.5,81.38) .. (271.19,85.59) ;
            \draw  [fill={rgb, 255:red, 0; green, 0; blue, 0 }  ,fill opacity=1 ] (325.14,55.63) .. controls (325.14,54.49) and (326.53,53.55) .. (328.26,53.55) .. controls (329.98,53.55) and (331.38,54.49) .. (331.38,55.63) .. controls (331.38,56.78) and (329.98,57.71) .. (328.26,57.71) .. controls (326.53,57.71) and (325.14,56.78) .. (325.14,55.63) -- cycle ;
            \draw  [fill={rgb, 255:red, 0; green, 0; blue, 0 }  ,fill opacity=1 ] (325.14,27.94) .. controls (325.14,26.79) and (326.53,25.86) .. (328.26,25.86) .. controls (329.98,25.86) and (331.38,26.79) .. (331.38,27.94) .. controls (331.38,29.09) and (329.98,30.02) .. (328.26,30.02) .. controls (326.53,30.02) and (325.14,29.09) .. (325.14,27.94) -- cycle ;
            \draw    (328.26,27.94) .. controls (366.05,-0.82) and (366.62,57.42) .. (328.47,27.94) ;
            \draw    (301.34,91.59) -- (301.34,118.04) ;
            \draw [shift={(301.34,120.04)}, rotate = 270] [color={rgb, 255:red, 0; green, 0; blue, 0 }  ][line width=0.75]    (10.93,-3.29) .. controls (6.95,-1.4) and (3.31,-0.3) .. (0,0) .. controls (3.31,0.3) and (6.95,1.4) .. (10.93,3.29)   ;
            \draw    (271.19,55.55) -- (328.26,55.63) ;
            \draw  [fill={rgb, 255:red, 0; green, 0; blue, 0 }  ,fill opacity=1 ] (325.14,84.67) .. controls (325.14,83.52) and (326.53,82.59) .. (328.26,82.59) .. controls (329.98,82.59) and (331.38,83.52) .. (331.38,84.67) .. controls (331.38,85.82) and (329.98,86.75) .. (328.26,86.75) .. controls (326.53,86.75) and (325.14,85.82) .. (325.14,84.67) -- cycle ;
            \draw    (328.26,55.63) .. controls (364.78,26.7) and (367.62,84.95) .. (328.47,55.63) ;
            \draw  [fill={rgb, 255:red, 0; green, 0; blue, 0 }  ,fill opacity=1 ] (268.07,55.55) .. controls (268.07,54.41) and (269.47,53.47) .. (271.19,53.47) .. controls (272.92,53.47) and (274.32,54.41) .. (274.32,55.55) .. controls (274.32,56.7) and (272.92,57.63) .. (271.19,57.63) .. controls (269.47,57.63) and (268.07,56.7) .. (268.07,55.55) -- cycle ;
            \draw  [fill={rgb, 255:red, 0; green, 0; blue, 0 }  ,fill opacity=1 ] (268.07,27.86) .. controls (268.07,26.71) and (269.47,25.78) .. (271.19,25.78) .. controls (272.92,25.78) and (274.32,26.71) .. (274.32,27.86) .. controls (274.32,29.01) and (272.92,29.94) .. (271.19,29.94) .. controls (269.47,29.94) and (268.07,29.01) .. (268.07,27.86) -- cycle ;
            \draw  [fill={rgb, 255:red, 0; green, 0; blue, 0 }  ,fill opacity=1 ] (268.07,84.59) .. controls (268.07,83.44) and (269.47,82.51) .. (271.19,82.51) .. controls (272.92,82.51) and (274.32,83.44) .. (274.32,84.59) .. controls (274.32,85.74) and (272.92,86.67) .. (271.19,86.67) .. controls (269.47,86.67) and (268.07,85.74) .. (268.07,84.59) -- cycle ;
            \draw    (271.19,27.86) .. controls (246.5,32.68) and (247.5,51.35) .. (271.19,55.55) ;
            \draw    (271.19,27.86) .. controls (220.5,29.38) and (222.5,86.38) .. (271.19,84.59) ;
            \draw    (271.19,84.59) -- (328.26,84.67) ;
            \draw    (271.19,27.86) -- (328.26,27.94) ;

            \end{tikzpicture}
            }
            &
            \scalebox{0.6}{
                \tikzset{every picture/.style={line width=0.75pt}}
                \begin{tikzpicture}[x=0.75pt,y=0.75pt,yscale=-1,xscale=1]
                
                \draw    (264.93,126.17) -- (321.72,126.17) ;
                \draw    (321.51,126.17) .. controls (358.03,97.23) and (360.87,155.48) .. (321.72,126.17) ;
                \draw    (264.93,126.17) .. controls (226.83,99.5) and (228.53,154.34) .. (265.14,126.17) ;
                \draw  [fill={rgb, 255:red, 0; green, 0; blue, 0 }  ,fill opacity=1 ] (262.02,126.17) .. controls (262.02,125.02) and (263.42,124.09) .. (265.14,124.09) .. controls (266.87,124.09) and (268.27,125.02) .. (268.27,126.17) .. controls (268.27,127.31) and (266.87,128.25) .. (265.14,128.25) .. controls (263.42,128.25) and (262.02,127.31) .. (262.02,126.17) -- cycle ;
                \draw  [fill={rgb, 255:red, 0; green, 0; blue, 0 }  ,fill opacity=1 ] (318.6,126.17) .. controls (318.6,125.02) and (320,124.09) .. (321.72,124.09) .. controls (323.45,124.09) and (324.85,125.02) .. (324.85,126.17) .. controls (324.85,127.31) and (323.45,128.25) .. (321.72,128.25) .. controls (320,128.25) and (318.6,127.31) .. (318.6,126.17) -- cycle ;
                \draw    (294.34,90.17) -- (294.34,116.62) ;
                \draw [shift={(294.34,118.62)}, rotate = 270] [color={rgb, 255:red, 0; green, 0; blue, 0 }  ][line width=0.75]    (10.93,-3.29) .. controls (6.95,-1.4) and (3.31,-0.3) .. (0,0) .. controls (3.31,0.3) and (6.95,1.4) .. (10.93,3.29)   ;
                \draw    (320.51,55.17) .. controls (357.03,26.23) and (359.87,84.48) .. (320.72,55.17) ;
                \draw    (263.93,55.17) .. controls (225.83,28.5) and (227.53,83.34) .. (264.14,55.17) ;
                \draw  [fill={rgb, 255:red, 0; green, 0; blue, 0 }  ,fill opacity=1 ] (261.02,55.17) .. controls (261.02,54.02) and (262.42,53.09) .. (264.14,53.09) .. controls (265.87,53.09) and (267.27,54.02) .. (267.27,55.17) .. controls (267.27,56.31) and (265.87,57.25) .. (264.14,57.25) .. controls (262.42,57.25) and (261.02,56.31) .. (261.02,55.17) -- cycle ;
                \draw  [fill={rgb, 255:red, 0; green, 0; blue, 0 }  ,fill opacity=1 ] (317.6,55.17) .. controls (317.6,54.02) and (319,53.09) .. (320.72,53.09) .. controls (322.45,53.09) and (323.85,54.02) .. (323.85,55.17) .. controls (323.85,56.31) and (322.45,57.25) .. (320.72,57.25) .. controls (319,57.25) and (317.6,56.31) .. (317.6,55.17) -- cycle ;
                \draw    (320.51,55.17) .. controls (315.5,23.92) and (269.38,21.84) .. (264.14,53.09) ;
                \draw    (263.93,55.17) .. controls (261.5,82.92) and (319.5,87.92) .. (320.72,55.17) ;
                \draw    (264.14,55.17) -- (320.94,55.17) ;
                \draw  [color={rgb, 255:red, 0; green, 0; blue, 0 }  ,draw opacity=1 ][fill={rgb, 255:red, 0; green, 0; blue, 0 }  ,fill opacity=1 ] (238.13,122.81) -- (242.27,130.44) -- (234.82,130.56) -- cycle ;
                \draw  [color={rgb, 255:red, 0; green, 0; blue, 0 }  ,draw opacity=1 ][fill={rgb, 255:red, 0; green, 0; blue, 0 }  ,fill opacity=1 ] (349.28,122.81) -- (353.12,130.44) -- (345.68,130.56) -- cycle ;
                
                \draw (224.5,118.07) node [anchor=north west][inner sep=0.75pt]    {$1$};
                \draw (358.5,114.93) node [anchor=north west][inner sep=0.75pt]    {$1$};
                \end{tikzpicture}
            }
            
        \\
        \hline
        \scalebox{0.6}{   
                \tikzset{every picture/.style={line width=0.75pt}} 
                \begin{tikzpicture}[x=0.75pt,y=0.75pt,yscale=-1,xscale=1]
                \draw    (253.02,66.85) -- (324.62,66.85) ;
                \draw    (324.35,66.85) .. controls (370.39,22.12) and (373.97,112.16) .. (324.62,66.85) ;
                \draw [line width=3]    (252.75,66.85) .. controls (204.72,25.63) and (206.87,110.41) .. (253.02,66.85) ;
                \draw  [fill={rgb, 255:red, 0; green, 0; blue, 0 }  ,fill opacity=1 ] (247.04,66.85) .. controls (247.04,63.63) and (249.6,61.02) .. (252.75,61.02) .. controls (255.9,61.02) and (258.46,63.63) .. (258.46,66.85) .. controls (258.46,70.07) and (255.9,72.68) .. (252.75,72.68) .. controls (249.6,72.68) and (247.04,70.07) .. (247.04,66.85) -- cycle ;
                \draw  [fill={rgb, 255:red, 0; green, 0; blue, 0 }  ,fill opacity=1 ] (320.69,66.85) .. controls (320.69,65.07) and (322.45,63.63) .. (324.62,63.63) .. controls (326.8,63.63) and (328.56,65.07) .. (328.56,66.85) .. controls (328.56,68.63) and (326.8,70.07) .. (324.62,70.07) .. controls (322.45,70.07) and (320.69,68.63) .. (320.69,66.85) -- cycle ;
                \draw    (356.5,73.53) -- (360.5,66.53) ;
                \draw    (360.5,66.53) -- (364.5,72.53) ;
                \draw  [color={rgb, 255:red, 0; green, 0; blue, 0 }  ,draw opacity=1 ][fill={rgb, 255:red, 0; green, 0; blue, 0 }  ,fill opacity=1 ] (218.02,61.7) -- (223.67,71.17) -- (212.68,71.34) -- cycle ;

                \draw (368,62.4) node [anchor=north west][inner sep=0.75pt]    {$0$};
                \draw (213,79.4) node [anchor=north west][inner sep=0.75pt]    {$1$};
                \end{tikzpicture}
                }
            & 
            \scalebox{0.6}{
                \tikzset{every picture/.style={line width=0.75pt}} 
                
                \begin{tikzpicture}[x=0.75pt,y=0.75pt,yscale=-1,xscale=1]
                \draw    (246.52,76.52) -- (318.12,76.52) ;
                \draw    (317.85,76.52) .. controls (363.89,31.79) and (367.47,121.83) .. (318.12,76.52) ;
                \draw [line width=0.75]    (246.25,76.52) .. controls (198.22,35.3) and (200.37,120.07) .. (246.52,76.52) ;
                \draw  [fill={rgb, 255:red, 0; green, 0; blue, 0 }  ,fill opacity=1 ] (314.19,76.52) .. controls (314.19,74.74) and (315.95,73.3) .. (318.12,73.3) .. controls (320.3,73.3) and (322.06,74.74) .. (322.06,76.52) .. controls (322.06,78.29) and (320.3,79.73) .. (318.12,79.73) .. controls (315.95,79.73) and (314.19,78.29) .. (314.19,76.52) -- cycle ;
                \draw    (348.5,82.42) -- (354,76.2) ;
                \draw    (354,76.2) -- (358,82.2) ;
                \draw    (206.5,82.42) -- (211,76.2) ;
                \draw    (211,76.2) -- (216.5,81.42) ;
                \draw  [fill={rgb, 255:red, 0; green, 0; blue, 0 }  ,fill opacity=1 ] (242.58,76.52) .. controls (242.58,74.74) and (244.35,73.3) .. (246.52,73.3) .. controls (248.7,73.3) and (250.46,74.74) .. (250.46,76.52) .. controls (250.46,78.29) and (248.7,79.73) .. (246.52,79.73) .. controls (244.35,79.73) and (242.58,78.29) .. (242.58,76.52) -- cycle ;
                \draw    (275.5,72.65) -- (282.32,76.52) ;
                \draw    (282.32,76.52) -- (275.5,80.65) ;
                \draw (361,76.07) node [anchor=north west][inner sep=0.75pt]    {$0$};
                \draw (202.5,88.07) node [anchor=north west][inner sep=0.75pt]    {$1$};
                \draw (272,53.07) node [anchor=north west][inner sep=0.75pt]    {$0$};
                \end{tikzpicture}
                
            }
            &
            \scalebox{0.6}{  
                \tikzset{every picture/.style={line width=0.75pt}} 
                \begin{tikzpicture}[x=0.75pt,y=0.75pt,yscale=-1,xscale=1]
                \draw    (265.52,75.52) -- (337.12,75.52) ;
                \draw [line width=3]    (336.85,75.52) .. controls (382.89,30.79) and (386.47,120.83) .. (337.12,75.52) ;
                \draw [line width=3]    (265.25,75.52) .. controls (217.22,34.3) and (219.37,119.07) .. (265.52,75.52) ;
                \draw  [fill={rgb, 255:red, 0; green, 0; blue, 0 }  ,fill opacity=1 ] (259.81,75.52) .. controls (259.81,72.3) and (262.37,69.68) .. (265.52,69.68) .. controls (268.67,69.68) and (271.23,72.3) .. (271.23,75.52) .. controls (271.23,78.74) and (268.67,81.35) .. (265.52,81.35) .. controls (262.37,81.35) and (259.81,78.74) .. (259.81,75.52) -- cycle ;
                \draw  [fill={rgb, 255:red, 0; green, 0; blue, 0 }  ,fill opacity=1 ] (331.14,75.52) .. controls (331.14,72.3) and (333.7,69.68) .. (336.85,69.68) .. controls (340.01,69.68) and (342.56,72.3) .. (342.56,75.52) .. controls (342.56,78.74) and (340.01,81.35) .. (336.85,81.35) .. controls (333.7,81.35) and (331.14,78.74) .. (331.14,75.52) -- cycle ;
                \draw  [color={rgb, 255:red, 0; green, 0; blue, 0 }  ,draw opacity=1 ][fill={rgb, 255:red, 0; green, 0; blue, 0 }  ,fill opacity=1 ] (230.02,70.45) -- (235.67,79.91) -- (224.68,80.08) -- cycle ;
                \draw  [color={rgb, 255:red, 0; green, 0; blue, 0 }  ,draw opacity=1 ][fill={rgb, 255:red, 0; green, 0; blue, 0 }  ,fill opacity=1 ] (372.02,69.65) -- (377.67,79.11) -- (366.68,79.28) -- cycle ;
                \draw (382.5,78.93) node [anchor=north west][inner sep=0.75pt]    {$1$};
                \draw (220.5,84.07) node [anchor=north west][inner sep=0.75pt]    {$1$};
                \end{tikzpicture}
                
            }
            \\
            \hline
    \end{tabular}
    \caption{Examples of 3-covers and their respective depictions using \cref{convention}.}
    \label{tab:p-decoratedGraph}
\end{table}

In \cref{tab:pCoversGenusTwo}, we list all maximal cells in the moduli space $\Delta_{2,p}$ for $p \geq 5$. These are obtained by listing all possible $p$-covers of the dumbbell graph $D$ and theta graph $T$, i.e.
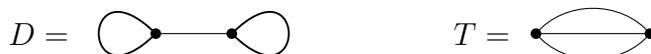
\vspace{-7mm}
\begin{center}
 \begin{tikzpicture}
          \vertex[]{0}{0};
	       \vertex[]{1}{0};

		   \draw (0,0) -- (1,0);
     
		  \selfloop{1}{0};
          \selfloopleft{0}{0};

          \draw (-1,0) node [anchor=east]  {$D = $};

          \vertex[]{5}{0};
	       \vertex[]{6.5}{0};
          \path[draw] (5, 0) to[bend left = 50](6.5,0);
          \draw (5, 0) -- (6.5,0);
          \path[draw] (5, 0) to[bend right = 50](6.5,0);
          
          \draw (4.8,0) node [anchor = east]  {$T = $};
      \end{tikzpicture}
\end{center}
\vspace{-4mm}
\noindent since these correspond to the two maximal cells of $M^\trop_2$.

\begin{table}[ht]
    \centering
    \scalebox{0.8}{
    \begin{tabular}{|c|c|c|}
     \hline
    \makecell{Covers in $\Delta_{2,p}^{par}$} &  \makecell{Covers with automorphism \\ acting as reflection on \\ the cell in $\Delta_{2, p}$}   &  \makecell{All other $p$-covers in $\Delta_{2,p}$}\\
     \hline
       \begin{tikzpicture}
          \fatvertex[]{0}{0};
	       \vertex[]{1}{0};
          \draw (0, 0) -- (1, 0);
    
		   \selfloop{1}{0};
          \thinarrowup{(1.75, 0.1)}
          \draw (2,0) node [xscale=1, yscale=1]  {$0$};

          \selfloopleft[0.75]{0}{0};
          \fatarrowup{(-0.75,-0.05)}
          \draw (-1,0) node [xscale=1, yscale=1]  {{\boldmath $i$}};

          \draw (0.4,-0.8) node [xscale=1, yscale=1]  {$1 \leq i \leq \frac{p-1}{2}$};
          
      \end{tikzpicture}  & 
     \begin{tikzpicture}
          \fatvertex[]{0}{0};
	       \fatvertex[]{1}{0};

		   \draw (0,0) -- (1,0);
     
		   \selfloop[0.75]{1}{0};
          \fatarrowup{(-0.75,-0.05)}
          \draw (2,0) node [xscale=1, yscale=1]  {{\boldmath $i$}};

          \selfloopleft[0.75]{0}{0};
          \fatarrowup{(1.75,-0.05)}
          
          \draw (-1,0) node [xscale=1, yscale=1]  {{\boldmath $i$}};

          \draw (0.4,-0.8) node [xscale=1, yscale=1]  {$1 \leq i \leq \frac{p-1}{2}$};
      \end{tikzpicture}
      & 
      \begin{tikzpicture}
          \draw (-2,1) node [xscale=1, yscale=1]  {$\boxed{D_{i}^{j}}$};
          
          \fatvertex[]{0}{0};
	       \vertex[]{1}{0};

		   \draw (0,0) -- (1,0);
     
		   \selfloop{1}{0};
          \fatarrowup{(-0.75,-0.05)}
          \draw (2,0) node [xscale=1, yscale=1]  {$j$};
          
          \selfloopleft[0.75]{0}{0};
          \thinarrowup{(1.75, 0.1)}
          \draw (-1,0) node [xscale=1, yscale=1]  {{\boldmath $i$}};

          \draw (0.4,-0.8) node [xscale=1, yscale=1]  {$1 \leq i, j \leq \frac{p-1}{2}$};

      \end{tikzpicture}
      
    \\

      
      \hline
      \begin{tikzpicture}
          \vertex[]{0}{0};
	       \vertex[]{1}{0};
     
		   \draw (0,0) -- (1,0);
     
		   \selfloop{1}{0};
          \thinarrowup{(-0.75, 0.1)}
          \draw (2,0) node [xscale=1, yscale=1]  {$0$};

          \selfloopleft{0}{0};
          \thinarrowup{(1.75, 0.1)}
          \draw (-1,0) node [xscale=1, yscale=1]  {$i$};

          \draw (0.4,-0.8) node [xscale=1, yscale=1]  {$1 \leq i \leq \frac{p-1}{2}$};
          
      \end{tikzpicture}
       & 
       \begin{tikzpicture}
       
          \vertex[]{0}{0};
	       \vertex[]{1}{0};

		   \draw (0,0) -- (1,0);
     
		   \selfloop{1}{0};
          \thinarrowup{(-0.75, 0.1)}
          \draw (2,0) node [xscale=1, yscale=1]  {$i$};

          \selfloopleft{0}{0};
          \thinarrowup{(1.75, 0.1)}
          \draw (-1,0) node [xscale=1, yscale=1]  {$i$};

          \draw (0.4,-0.8) node [xscale=1, yscale=1]  {$1 \leq i \leq \frac{p-1}{2}$};
      \end{tikzpicture} 
      & 
      \begin{tikzpicture}
          \draw (-2,1) node [xscale=1, yscale=1]  {$\boxed{D^{i,j}}$};
      
          \vertex[]{0}{0};
	       \vertex[]{1}{0};

		   \draw (0,0) -- (1,0);
     
		   \selfloop{1}{0};
          \thinarrowup{(-0.75, 0.1)}
          \draw (2,0) node [xscale=1, yscale=1]  {$j$};

          \selfloopleft{0}{0};
          \thinarrowup{(1.75, 0.1)}
          \draw (-1,0) node [xscale=1, yscale=1]  {$i$};

          \draw (0.4,-0.8) node [xscale=1, yscale=1]  {$1 \leq i < j \leq \frac{p-1}{2}$};
          
      \end{tikzpicture}  \\

      
      \hline
       \begin{tikzpicture}
          \vertex[]{0}{0};
	       \vertex[]{1.5}{0};

          \draw (0.8,0.8) node  [xscale=0.8, yscale=0.8] {$0$};
          \thinarrowright{(0.75, 0.6)} 
          \draw (0, 0) .. controls (0.5, 0.8) and (1, 0.8) .. (1.5,0);

          \draw (0.8,0.2) node  [xscale=0.8, yscale=0.8] {$0$};
          \thinarrowright{(0.75, 0)} 
          \draw (0, 0) -- (1.5,0);

          \draw (0.8,-0.4) node  [xscale=0.8, yscale=0.8] {$i$};
          \thinarrowright{(0.75, -0.6)} 
          \draw (0, 0) .. controls (0.5, -0.8) and (1, -0.8) .. (1.5,0);

          \draw (0.4,-1.2) node [xscale=1, yscale=1]  {$1 \leq i \leq \frac{p-1}{2}$};
      \end{tikzpicture} 
      & 
      \begin{tikzpicture}
          \vertex[]{0}{0};
	       \vertex[]{1.5}{0};

          \draw (0.8,0.8) node  [xscale=0.8, yscale=0.8] {$0$};
          \thinarrowright{(0.75, 0.6)} 
          \draw (0, 0) .. controls (0.5, 0.8) and (1, 0.8) .. (1.5,0);

          \draw (0.8,0.2) node  [xscale=0.8, yscale=0.8] {$i$};
          \thinarrowright{(0.75, 0)} 
          \draw (0, 0) -- (1.5,0);

          \draw (0.7,-0.4) node  [xscale=0.8, yscale=0.8] {$-i$};
          \thinarrowright{(0.75, -0.6)} 
          \draw (0, 0) .. controls (0.5, -0.8) and (1, -0.8) .. (1.5,0);

          \draw (0.4,-1.2) node [xscale=1, yscale=1]  {$1 \leq i \leq \frac{p-1}{2}$};
      \end{tikzpicture} 
      & 
      \begin{tikzpicture}
          \draw (-2,1) node [xscale=1, yscale=1]  {$\boxed{D_{i,j}}$};
          \fatvertex[]{0}{0};
	       \fatvertex[]{1}{0};

		   \draw (0,0) -- (1,0);
     
		   \selfloop[0.75]{1}{0};
          \fatarrowup{(-0.75, -0.05)}
          \draw (2,0) node [xscale=1, yscale=1]  {{\boldmath$j$}};

          \selfloopleft[0.75]{0}{0};
          \fatarrowup{(1.75, -0.05)}
          \draw (-1,0) node [xscale=1, yscale=1]  {{\boldmath$i$}};

          \draw (0.4,-0.8) node [xscale=1, yscale=1]  {$1 \leq i < j \leq \frac{p-1}{2}$};
      \end{tikzpicture} 
     \\

      
      \hline
       & 
       \begin{tikzpicture}
          \fatvertex[]{0}{0};
	       \fatvertex[]{1.5}{0};

          \draw (0.9,0.8) node  [xscale=0.8, yscale=0.8] {$0$};
          \thinarrowright{(0.75, 0.6)}
          \draw (0, 0) .. controls (0.5, 0.8) and (1, 0.8) .. (1.5,0);

          \draw (0.9,0.2) node  [xscale=0.8, yscale=0.8] {{\boldmath $i$}};
          \fatarrowright{(0.75, 0)}
          \draw [line width = 0.75 mm] (0, 0) -- (1.5,0);

          \draw (0.8,-0.35) node  [xscale=0.8, yscale=0.8] {{\boldmath $-i$}};
          \fatarrowright{(0.75, -0.6)}
          \draw [line width = 0.75 mm] (0, 0) .. controls (0.5, -0.8) and (1, -0.8) .. (1.5,0);

          \draw (0.4,-1.2) node [xscale=1, yscale=1]  {$1 \leq i \leq \frac{p-1}{2}$};
      \end{tikzpicture}
      & 
        \begin{tikzpicture}
          \draw (-2,1) node [xscale=1, yscale=1] {$\boxed{T^{0,i,j}}$};
          \vertex[]{-1}{0};
	       \vertex[]{0.5}{0};

          \draw (-0.2,0.8) node  [xscale=0.8, yscale=0.8] {$0$};
          \thinarrowright{(-0.25, 0.6)}
          \draw (-1, 0) .. controls (-0.5, 0.8) and (0, 0.8) .. (0.5,0);
          
          \draw (-0.2,0.2) node  [xscale=0.8, yscale=0.8] {$i$};
          \thinarrowright{(-0.25, 0)}
          \draw (-1, 0) -- (0.5,0);

          \draw (-0.2,-0.4) node  [xscale=0.8, yscale=0.8] {$j$};
          \thinarrowright{(-0.25, -0.6)}
          \draw (-1, 0) .. controls (-0.5, -0.8) and (0, -0.8) .. (0.5,0);


          \draw (-0.3,-1.2) node [xscale=1, yscale=1]  {see \cref{eq:SystemOfReps}};
          
      \end{tikzpicture}
      \\


      \hline
      
       &
       
        \begin{tikzpicture}
          \fatvertex[]{0}{0};
	       \fatvertex[]{1.5}{0};

          \draw (0.9,0.8) node  [xscale=0.8, yscale=0.8] {{\boldmath $i$}};
          \fatarrowright{(0.75, 0.6)}
          \draw [line width = 0.75 mm] (0, 0) .. controls (0.5, 0.8) and (1, 0.8) .. (1.5,0);
          
          \draw (0.9,0.2) node  [xscale=0.8, yscale=0.8] {{\boldmath $i$}};
          \fatarrowright{(0.75, 0)}
          \draw [line width = 0.75 mm] (0, 0) -- (1.5,0);

          \draw (0.8,-0.9) node  [xscale=0.8, yscale=0.8] {{\boldmath $-2i$}};
          \fatarrowright{(0.75, -0.6)}
          \draw [line width = 0.75 mm] (0, 0) .. controls (0.5, -0.8) and (1, -0.8) .. (1.5,0);
          
          \draw (0.4,-1.6) node [xscale=1, yscale=1]  {$1 \leq i \leq \frac{p-1}{2}$};
      \end{tikzpicture}
       
       &
       
      \begin{tikzpicture}
          \draw (-2,1) node [xscale=1, yscale=1] {$\boxed{T_{i,j,k}}$};
          \fatvertex[]{-1}{0};
	       \fatvertex[]{0.5}{0};

          \draw (-0.1,0.8) node  [xscale=0.8, yscale=0.8] {{\boldmath $i$}};
          \fatarrowright{(-0.25, 0.6)}
          \draw [line width = 0.75 mm] (-1, 0) .. controls (-0.5, 0.8) and (0, 0.8) .. (0.5,0);

          \draw (-0.1,0.25) node  [xscale=0.8, yscale=0.8] {{\boldmath $j$}};
          \fatarrowright{(-0.25, 0)}
          \draw [line width = 0.75 mm] (-1, 0) -- (0.5,0);

          \draw (-0.1,-0.30) node  [xscale=0.8, yscale=0.8] {{\boldmath $k$}};
          \fatarrowright{(-0.25, -0.6)}
          \draw [line width = 0.75 mm] (-1, 0) .. controls (-0.5, -0.8) and (0, -0.8) .. (0.5,0);

          \draw (-0.3,-1.2) node [xscale=1, yscale=1]  {see \cref{eq:sysRepDilated}};

      \end{tikzpicture}\\
      \hline
    \end{tabular}
    }
    \caption{The $p$-covers corresponding to maximal cells in $\Delta_{2, p}$. We use \cref{convention}.}
    \label{tab:pCoversGenusTwo}
    \end{table}

Let $X$ be the symmetric $\Delta$-complex generated by the cells corresponding to covers in the first and last column of \cref{tab:pCoversGenusTwo}. Then $\Delta_{2, p}$ is homotopy equivalent to $X$. Indeed, the covers in the second column of \cref{tab:pCoversGenusTwo} all have an automorphism which acts as reflection on the corresponding simplex in $\Delta_{2, p}$. This makes it possible to deformation retract their interior down to their $1$-skeleton in $\Delta_{2,p}$. After further contracting some edges, we arrive at the subcomplex $X$.

Let us now focus on a careful study of the cells in the third column of \cref{tab:pCoversGenusTwo}.

\begin{lemma}\label{lem:NumberOfDistinctThetaCovers}
    Up to isomorphism, there are $\frac{p^2-1}{12}$ free covers of $T$ such that every edge in the graph has a distinct gain. This counts the total number of covers in the entries (3,2) and (4,3) in \cref{tab:pCoversGenusTwo}.
\end{lemma}
\begin{proof}
    Every free cover of $T$ can be represented as some $T^{i,j,k}$ where $i,j,k \in \Z/p$ are pairwise distinct values for the gains of the three edges. This representation is not unique, since the same $p$-cover can also be written as:
    \begin{enumerate}
        \item $T^{\sigma(i),\sigma(j),\sigma(k)}$ where $\sigma$ is any permutation of $\{i,j,k\}$,
        \item $T^{i + a , j + a , k + a }$ for $a \in \Z/p$ (this is switching of the gain graph presentation),
        \item $T^{-i ,-j, -k }$,
    \end{enumerate}
    or a combination of these. Let $\mathcal{T} = \big\{ [T^{i,j,k}] \colon i,j,k \in \Z/p \text{ distinct} \big\}$ where $[ - ]$ denotes isomorphism class. Note that there is a natural bijection between $\mathcal{T}$ and the set $\mathcal{B}$ of binary $3$-bracelets of length $p$, i.e. configurations of subsets of size $3$ in a regular $p$-gon, up to the natural action of the dihedral group $D_{2p}$. Counting the number of elements in $\mathcal{B}$ can be carried out using the Redfield-P\'olya counting theorem, see the theorem on page 98 of \cite{PolyaRead87}. This theorem uses the \emph{cycle index polynomial} $\bm{Z}_{D_{2p}}(T_1, \dots, T_p)$ of the dihedral group $D_{2p}$ acting on the regular $p$-gon by rotations and reflections, which is given by:
    \begin{align*}
        \bm{Z}_{D_{2p}}(z_1, \dots, z_p)
        \coloneqq{}& \frac{1}{|D_{2p}|} \sum_{\sigma \in D_{2p}} \ z_1^{c_1(\sigma)} \cdots z_p^{c_p(\sigma)}\\
        ={}& \frac{1}{2p} \left( z_1^p + p \ z_1 z_2^{\frac{p-1}{2}} + (p-1)z_p \right),
    \end{align*}
    where $c_k(\sigma)$ is the number of $k$-cycles in the cycle decomposition of the permutation of the $p$-gon induced by $\sigma \in D_{2p}$. By P\'olya's theorem, the size of $\mathcal{B}$ for $p \geq 5$ is the $t^3$-coefficient in the following expression:
    \begin{equation*} 
         | \mathcal{B} | = |\mathcal{T}|  = [t^3]\Big( \bm{Z}_{D_p}(1+t, 1+t^2, \dots, 1+t^{p}) \Big) = \frac{p^2 - 1}{12}.
    \end{equation*}    
    Moreover, it is not so difficult to see that the following form a system of representatives in $\mathcal{T}$:
    \begin{equation} \label{eq:SystemOfReps}
    T^{0, i, j} \quad \text{for } 1 \leq i < \frac{p + 1}{3} \text{ and } 2i \leq j \leq \left\lfloor \frac{p+i}{2} \right\rfloor.
    \end{equation}
    This finishes the proof.
\end{proof}

Some of the free covers $T^{i,j,k}$ we just counted in \cref{lem:NumberOfDistinctThetaCovers} do actually have automorphisms acting as reflections. The following lemma shows which of these covers have such an automorphism, justifying the grouping in \cref{tab:pCoversGenusTwo}.

\begin{lemma}\label{lem:AutTheta}
    Suppose $p \geq 5$ and let $i,j,k$ be distinct elements in $\Z/p$ and $T^{i,j,k}$ the corresponding free cover. Let $\epsilon_{i,j} \in \Aut(T)$ denote the automorphism that swaps the two vertices and the two edges decorated by $i,j$. Then $\epsilon_{i,j}$ lifts to an automorphism $\tilde{\epsilon}_{i,j}$ of the $p$-cover $T^{i,j,k}$ if and only if $i+j = 2k$. Moreover, we have the following:
    \[
    \Aut(T^{i,j,k}) \cong \begin{cases} D_{2p} & \text{if } i + j + k \in \{3i, 3j, 3k\}, \\ \Z/p  & \text{otherwise.} \end{cases}
    \]
    In particular, the 2-simplex in $\Delta_{2,p}$ corresponding to the $p$-cover $T^{i,j,k}$ is folded in half if and only if $i+j+k \in \{3i,3j,3k\}$.
    
\end{lemma}
\begin{proof}
    We trivialize the cover $T^{i,j,k}$ along the edge with gain equal to $k$, hence without loss of generality we take $k = 0$. To simplify notation we denote the $p$-cover $T^{i,j,k}$ by $\pi \colon \tilde{T} \to T$. Choosing a vertex $v_0 \in \pi^{-1}(v)$, the $\Z/p$-action gives a canonical labeling of the preimages of $v$ by $v_0, \ldots, v_{p-1}$ and by compatibility with the trivialization $k = 0$ for the vertices above $v'$ by $v'_0, \ldots, v'_{p-1}$.

    Since $\pi$ is a $p$-cover, $\Z/p$ is a normal subgroup of $\Aut(\pi)$. We show that the quotient $\Aut(\pi) / (\Z/p)$ is either trivial or $\Z/2$. For this let $(\tilde f, f)$ be an automorphism of $\pi$ with $\tilde{f}(v_{0}) \in \{v_{0},v_{0}'\}$.

     Suppose first that $\tilde{f}(v_0) = v_0$. In this case we claim that $\tilde{f}$ is the identity. Indeed, since the gains of the three edges of $T$ are distinct, $\tilde f$ is fully determined as a graph automorphism of $\tilde{T}$ by its action on the fiber $\pi^{-1}(v')$. Moreover, $\tilde f|_{\pi^{-1}(v')}$ is an automorphism of the $\Z/p$-torsor $\pi^{-1}(v')$, and hence the action of $\tilde f$ on that fiber is fully determined by its action on any vertex over $v'$. In particular, this implies that the order of $\tilde f$ divides $p$. Furthermore, since $\tilde{f}$ maps edges to edges it permutes the vertices $\{v_{0}', v_{i}', v_{j}'\}$ so its order divides $6$. Combined, this implies that the order of $\tilde f$ is $1$ and hence $\tilde{f} = \id$ as claimed.
    
     Now suppose that $\tilde{f}(v_0) = v_{0}'$. Since $\tilde{f}$ must map edges to edges we then have $\tilde{f}(v_{0}') \in \{v_0, v_{-i}, v_{-j} \}$. We distinguish the following three cases:
     
    \begin{enumerate}[leftmargin = 40pt, itemindent = 30pt] 
    
        \item[\textbf{Case 1:}] Suppose that $i+j \not \in \{0,3i,3j\}$. We show that this assumption gives a contradiction. Indeed, if $\tilde{f}(v_0') = v_{0}$ then by $\Z/p$-equivariance we have $\tilde{f}(v_i') = v_i$. Since edges are mapped to edges the edge $(v_0, v_{i}')$ is mapped to the edge $(v_0', v_i)$ which implies that $0 \in \{i, 2i, i+j\}$ which is a contradiction. If $\tilde{f}(v_0') = v_{-i}$ then by $\Z/p$-equivariance we find $\tilde{f}(v_{j}') = v_{j-i}$ and the edge $(v_0, v'_{j})$ is then mapped to the edge $(v_{0}',v_{j-i})$ which implies that $0 \in \{j-i, j, 2j-i\}$ which is a contradiction. Similarly if $\tilde{f}(v_{0}') = v_{-j}$ we obtain a contradiction. So in this case we have $\Aut(\pi) = \Z/p$.

        \item[\textbf{Case 2:}] Suppose that $i+j = 0$. We claim that in this case $\tilde{f}(v'_0) = v_{0}$ and $\tilde f$ lifts $\epsilon_{i,j}$. To see why, suppose for instance that $\tilde{f}(v'_0) = v_{-i}$. Then the edge $(v_0, v_j')$ is mapped to  $(v_{j-i}, v_{0}')$ which implies that $0 \in \{j-i, 2j -i , j\}$. This implies that $2j -i = 0$. Since $i+j = 0$ we then deduce that $3j = 0$. Because $3$ is invertible in $\Z/p$ we then find $j = 0$ which is a contradiction. So $\tilde{f}(v'_{0}) \neq v_{-i}$ and similarly we can also obtain $\tilde{f}(v'_{0}) \neq v_{-j}$. Hence $\tilde{f}(v'_0) = v_0$. By $\Z/p$-equivariance we get:
            \[
                \tilde f(v_\ell) = v'_{\ell} \quad \text{ and } \quad \tilde f(v_\ell') = v_{\ell}
            \]
            and one can check that $(\tilde{f}, \epsilon_{i,j})$ is a valid $\Z/p$-equivariant automorphism of $T^{0,i,j}$ that maps edges $\Z/p$-equivariantly as
            \begin{equation*}
                (v_0, v_0') \longmapsto (v_0',v_0), \quad  (v_0, v_i') \longmapsto (v_0',v_i), \quad  (v_0, v_j') \longmapsto (v_0',v_j)
            \end{equation*}
            and we denote $\tilde{f}$ by $\tilde{\epsilon}_{i,j}$. So, in this case we deduce that $\Aut(\pi)$ is the dihedral group $D_{2p}$ generated by the reflection $\tilde \epsilon_{i,j}$ and the rotations in $\Z/p$.
                
        \item[\textbf{Case 3:}] Suppose that $i+j = 3i$. This case reduces to the previous case by replacing $\{i,j,0\}$ with $\{i', j', 0\}$ where $i' = -i$ and $j' = j-i$.
    \end{enumerate} 
    All the cases have been exhausted so this finishes the proof.
\end{proof}

\begin{lemma}\label{lem:AutomorphismsFreeTheta}
    There are $\frac{(p-1)(p-5)}{12}$ isomorphism classes of fully dilated covers of $T$ which have no nontrivial automorphisms (entry $(5, 3)$ in \cref{tab:pCoversGenusTwo}) and $\frac{p-1}{2}$ which have an automorphism that acts on $T$ by swapping two edges (entry $(5, 2)$ in \cref{tab:pCoversGenusTwo}). 
\end{lemma}

\begin{proof}
    After choosing an orientation for the edges of $T$, the fully dilated $p$-covers of $T$ are obtained by decorating the $3$ edges with elements $(i,j,k) \in ((\Z/p)^{\times})^3$ such that $i+j+k = 0$. It is clear that $p$-covers of the form $T_{i,i,-2i}$ have an automorphism that swaps the two edges of $T$ with the same dilation flow. There are $\frac{p-1}{2}$ isomorphism classes of covers of the form $T_{i, i, -2i}$, one for each $ 1 \leq i \leq \frac{p-1}{2}$.
    
    On the other hand, it is not so hard to see that if $i,j,k$ are three distinct elements in $(\Z/p)^\times$ such that $i+j+k=0$ then the cover $T_{i,j,k}$ has no non-trivial automorphisms. Now we shall count the isomorphism classes of such covers. The covers in the isomorphism class of $T_{i,j,k}$ are the $p$-covers represented by:
    
    \begin{enumerate}
        \item $T_{-i,-j,-k}$, and
        
        \item $T_{\sigma(i),\sigma(j),\sigma(k)}$ for $\sigma$ a permutation of $\{i,j,k\}$.
    \end{enumerate}
    So the number of these isomorphism classes is $| A | / (2 \cdot 6)$,
    where $A$ is the following set:
    \begin{equation} \label{eq:sysRepDilated}
        \begin{aligned}
            A ={}& \big\{ (i,j,k) \in ((\Z/p)^\times)^3 \colon i+j+k = 0 \text{ and } i,j,k \text{ are distinct}\big\} \\
            ={}& \big\{ (i,j,-i-j) \in ((\Z/p)^\times)^3 \colon j \not \in \{ i,-i, -2i, -i/2 \} \text{ are distinct} \big\}.
        \end{aligned}
    \end{equation}
    So we easily see that $|A| = (p-1)(p-5)$ and hence the count $\frac{(p-1)(p-5)}{12}$.
\end{proof}

This finishes the enumeration of top-dimensional cells in \cref{tab:pCoversGenusTwo}. We now work towards the proof of \cref{thm:MainThmGenus2}.

\begin{lemma} \label{lem:CompYContractible}
    We define a subcomplex $Y$ of $X$ as the complex generated by the following cells (see \cref{tab:pCoversGenusTwo}):
    \begin{enumerate}
        \item $D_{1}^{i}$ and $D_{i}^{1}$ for $1 \leq i \leq (p-1)/2$, and 
        \item $D_{i,j}$ and $D^{i,j}$ for $1 \leq i<j \leq (p-1)/2$.
    \end{enumerate}
    The complex $Y$ is contractible.
\end{lemma}

\begin{proof}
    We define the following nested subcomplexes in $Y$:
    \begin{enumerate}
        \item $Y_1$ is generated by the cells $D_{1}^{i}$ for $1 \leq i \leq (p-1)/2$, and
        \item $Y_2$ is generated by $Y_1$ and $D_{i}^{1}$ for $2 \leq i \leq (p-1)/2$.
    \end{enumerate}
    To see why $Y$ is contractible, we argue repeatedly as follows: if in a complex $W$, there is a 2-dimensional cell $\sigma$ with a facet $\tau$ which is not a facet of any other top dimensional cell in $W$, then $W$ deformation retracts to the subcomplex of $W$ which omits $\sigma$.
    
    Start with the complex $Y$. Each of the cells $D_{i,j}$ and $D^{i, j}$ have a unique such facet, which correspond to the covers (drawn in \cref{convention})

    \vspace{-6mm}
    
    \begin{center}
        \begin{tikzpicture}
          \fatvertex[]{0}{0};

		  \selfloop[0.75]{0}{0};
          \fatarrowup{(-0.75, -0.05)}
          \draw (1,0) node [xscale=1, yscale=1]  {{\boldmath$j$}};
          
          \selfloopleft[0.75]{0}{0};
          \fatarrowup{(0.75, -0.05)}
          \draw (-1,0) node [xscale=1, yscale=1]  {{\boldmath$i$}};

      \end{tikzpicture} 
        \qquad \raisebox{2em}{and} \qquad
      \begin{tikzpicture}
              \vertex[]{0}{0};
    		  \selfloop{0}{0};
              \thinarrowup{(-0.75, 0.1)}
              \draw (1,0) node [xscale=1, yscale=1]  {$j$};

              \selfloopleft{0}{0};
              \thinarrowup{(0.75, 0.1)}
              \draw (-1,0) node [xscale=1, yscale=1]  {$i$};
          \end{tikzpicture}
    \end{center}
    \vspace{-7mm}
    for $1 \leq i < j \leq (p-1)/2$, respectively. Therefore, we have a deformation retraction from $Y$ onto $Y_2$. 
        
    Now we repeat the same argument with the cells $D^{1}_{i}$ for $2 \leq i \leq \frac{p-1}{2}$, each of which has the facet corresponding to
    
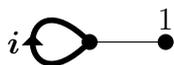
\vspace{-7mm}
    \begin{center}
        \begin{tikzpicture}  
               \fatvertex[]{0}{0};
    	       \fatvertex[1]{1}{0};
    
    		   \draw (0,0) -- (1,0);
               \fatarrowup{(-0.75,-0.05)}
               
               \selfloopleft[0.75]{0}{0};
               \draw (-1,0) node [xscale=1, yscale=1]  {{\boldmath $i$}};  
          \end{tikzpicture}
    \end{center}
    \vspace{-7mm}
    which does not belong to any other two dimensional cell of $Y_2$ (in the picture both vertices are dilated). Therefore, $Y_2$ deformation retracts to $Y_1$. 

    Finally, the complex $Y_1$ is a cone made out of $(p-1)/2$ triangles that are joined along two edges (its cone point corresponds to the cover $R_1$, see \cref{tab:vertices}). So $Y_1$ is contractible. In summary, $Y$ is contractible.
\end{proof}

\begin{proof}[Proof of \cref{thm:MainThmGenus2}]
The statements about $\Delta_{2,2}$ and $\Delta_{2,3}$ can be seen from Figures~\ref{fig:Delta_2_2} and~\ref{fig:Delta_2_3}, respectively. Now let $p \geq 5$. The total number of 2-dimensional cells in $\Delta_{2, p}$ can be counted from \cref{tab:pCoversGenusTwo} and Lemmas~\ref{lem:NumberOfDistinctThetaCovers} and~\ref{lem:AutomorphismsFreeTheta}. Concerning the homotopy type, recall that $X$ is the complex generated by all the maximal cells in the first and third columns of \cref{tab:pCoversGenusTwo} and that $\Delta_{2,p}$ is homotopy equivalent to $X$.

The complex $Y$ of \cref{lem:CompYContractible} is contractible and $\Delta_{2,p}^{par}$ is contractible by virtue of \cref{thm:contractible_loci}. Their intersection is the star-shaped 1-dimensional subcomplex depicted in \cref{fig:intersectionOfLoci} with the butterfly (see \cref{tab:vertices}) as its central vertex. Clearly the intersection is contractible and hence $Y \cup \Delta_{2,p}^{par}$ is contractible as well. 

\begin{figure}[H]
    \begin{center}
    \scalebox{0.6}{
            \tikzset{every picture/.style={line width=0.75pt}} 
            
            \begin{tikzpicture}[x=0.75pt,y=0.75pt,yscale=-1,xscale=1]
            
            \draw [color={rgb, 255:red, 0; green, 0; blue, 255 }  ,draw opacity=1 ][line width=3]    (318.73,89.77) -- (64.29,286.07) ;
            
            \draw [color={rgb, 255:red, 0; green, 0; blue, 255 }  ,draw opacity=1 ][line width=3]    (318.73,89.77) -- (192.47,285.13) ;
            
            \draw [color={rgb, 255:red, 0; green, 0; blue, 255 }  ,draw opacity=1 ][line width=3]    (318.73,89.77) -- (529.17,286.07) ;
            
            \draw  [color={rgb, 255:red, 255; green, 0; blue, 0 }  ,draw opacity=1 ][fill={rgb, 255:red, 255; green, 0; blue, 0 }  ,fill opacity=1 ] (312.04,89.77) .. controls (312.04,86.16) and (315.04,83.23) .. (318.73,83.23) .. controls (322.43,83.23) and (325.43,86.16) .. (325.43,89.77) .. controls (325.43,93.39) and (322.43,96.31) .. (318.73,96.31) .. controls (315.04,96.31) and (312.04,93.39) .. (312.04,89.77) -- cycle ;
            
            \draw  [color={rgb, 255:red, 255; green, 0; blue, 0 }  ,draw opacity=1 ][fill={rgb, 255:red, 255; green, 0; blue, 0 }  ,fill opacity=1 ] (522.48,286.07) .. controls (522.48,282.45) and (525.48,279.52) .. (529.17,279.52) .. controls (532.87,279.52) and (535.87,282.45) .. (535.87,286.07) .. controls (535.87,289.68) and (532.87,292.61) .. (529.17,292.61) .. controls (525.48,292.61) and (522.48,289.68) .. (522.48,286.07) -- cycle ;
            
            \draw  [color={rgb, 255:red, 255; green, 0; blue, 0 }  ,draw opacity=1 ][fill={rgb, 255:red, 255; green, 0; blue, 0 }  ,fill opacity=1 ] (185.78,285.13) .. controls (185.78,281.52) and (188.77,278.59) .. (192.47,278.59) .. controls (196.17,278.59) and (199.17,281.52) .. (199.17,285.13) .. controls (199.17,288.75) and (196.17,291.68) .. (192.47,291.68) .. controls (188.77,291.68) and (185.78,288.75) .. (185.78,285.13) -- cycle ;
            
            \draw  [color={rgb, 255:red, 255; green, 0; blue, 0 }  ,draw opacity=1 ][fill={rgb, 255:red, 255; green, 0; blue, 0 }  ,fill opacity=1 ] (57.6,286.07) .. controls (57.6,282.45) and (60.6,279.52) .. (64.29,279.52) .. controls (67.99,279.52) and (70.99,282.45) .. (70.99,286.07) .. controls (70.99,289.68) and (67.99,292.61) .. (64.29,292.61) .. controls (60.6,292.61) and (57.6,289.68) .. (57.6,286.07) -- cycle ;
            
            \draw [color={rgb, 255:red, 255; green, 0; blue, 0 }  ,draw opacity=1 ]   (299.15,64.29) .. controls (342.37,15.5) and (343.75,111.15) .. (299.35,64.29) ;
            
            \draw [color={rgb, 255:red, 255; green, 0; blue, 0 }  ,draw opacity=1 ]   (332.37,62.98) -- (325.43,67.34) ;
            
            \draw [color={rgb, 255:red, 255; green, 0; blue, 0 }  ,draw opacity=1 ]   (338.82,68.27) -- (332.37,62.98) ;
            
            \draw  [color={rgb, 255:red, 255; green, 0; blue, 0 }  ,draw opacity=1 ][fill={rgb, 255:red, 255; green, 0; blue, 0 }  ,fill opacity=1 ] (295.5,64.29) .. controls (295.5,62.08) and (297.23,60.29) .. (299.35,60.29) .. controls (301.48,60.29) and (303.2,62.08) .. (303.2,64.29) .. controls (303.2,66.51) and (301.48,68.3) .. (299.35,68.3) .. controls (297.23,68.3) and (295.5,66.51) .. (295.5,64.29) -- cycle ;
            
            \draw [color={rgb, 255:red, 0; green, 0; blue, 255 }  ,draw opacity=1 ]   (150.29,175.37) .. controls (170.26,98.21) and (231.04,185.57) .. (150.53,175.22) ;
            \draw [color={rgb, 255:red, 0; green, 0; blue, 255 }  ,draw opacity=1 ][line width=3]    (150.29,175.37) .. controls (125.67,251.84) and (71.53,162.46) .. (150.06,175.52) ;
            \draw  [color={rgb, 255:red, 0; green, 0; blue, 255 }  ,draw opacity=1 ][fill={rgb, 255:red, 0; green, 0; blue, 255 }  ,fill opacity=1 ] (146.25,177.96) .. controls (144.71,175.65) and (145.37,172.56) .. (147.74,171.05) .. controls (150.1,169.54) and (153.26,170.18) .. (154.81,172.48) .. controls (156.35,174.79) and (155.68,177.88) .. (153.32,179.39) .. controls (150.96,180.9) and (147.79,180.26) .. (146.25,177.96) -- cycle ;
            \draw [color={rgb, 255:red, 0; green, 0; blue, 255 }  ,draw opacity=1 ]   (185.93,147.71) -- (183.86,155.81) ;
            \draw [color={rgb, 255:red, 0; green, 0; blue, 255 }  ,draw opacity=1 ]   (194.38,149.27) -- (185.93,147.71) ;
            \draw  [color={rgb, 255:red, 0; green, 0; blue, 255 }  ,draw opacity=1 ][fill={rgb, 255:red, 0; green, 0; blue, 255 }  ,fill opacity=1 ] (109.67,195.88) -- (119.13,200.43) -- (110.4,206.16) -- cycle ;
            \draw [color={rgb, 255:red, 0; green, 0; blue, 255 }  ,draw opacity=1 ]   (278.01,208.06) .. controls (265.5,129.42) and (356.41,186.58) .. (278.17,207.83) ;
            \draw [color={rgb, 255:red, 0; green, 0; blue, 255 }  ,draw opacity=1 ][line width=3]    (278.01,208.06) .. controls (285.96,287.85) and (200.36,226.29) .. (277.86,208.29) ;
            \draw  [color={rgb, 255:red, 0; green, 0; blue, 255 }  ,draw opacity=1 ][fill={rgb, 255:red, 0; green, 0; blue, 255 }  ,fill opacity=1 ] (275.33,211.99) .. controls (272.99,210.46) and (272.36,207.36) .. (273.93,205.06) .. controls (275.5,202.77) and (278.67,202.15) .. (281.01,203.68) .. controls (283.35,205.21) and (283.98,208.31) .. (282.41,210.6) .. controls (280.84,212.89) and (277.67,213.51) .. (275.33,211.99) -- cycle ;
            \draw [color={rgb, 255:red, 0; green, 0; blue, 255 }  ,draw opacity=1 ]   (299.74,168.98) -- (301.08,177.23) ;
            \draw [color={rgb, 255:red, 0; green, 0; blue, 255 }  ,draw opacity=1 ]   (308.14,167.18) -- (299.74,168.98) ;
            \draw  [color={rgb, 255:red, 0; green, 0; blue, 255 }  ,draw opacity=1 ][fill={rgb, 255:red, 0; green, 0; blue, 255 }  ,fill opacity=1 ] (248.84,242.46) -- (259.36,243.03) -- (253.62,251.64) -- cycle ;
            \draw [color={rgb, 255:red, 0; green, 0; blue, 255 }  ,draw opacity=1 ]   (447.26,169.55) .. controls (366.29,160.93) and (446.4,89.99) .. (447.07,169.34) ;
            \draw [color={rgb, 255:red, 0; green, 0; blue, 255 }  ,draw opacity=1 ][line width=3]    (447.26,169.55) .. controls (528.17,182.78) and (445.09,247.56) .. (447.45,169.76) ;
            \draw  [color={rgb, 255:red, 0; green, 0; blue, 255 }  ,draw opacity=1 ][fill={rgb, 255:red, 0; green, 0; blue, 255 }  ,fill opacity=1 ] (450.44,173.1) .. controls (448.32,174.92) and (445.1,174.7) .. (443.24,172.62) .. controls (441.38,170.55) and (441.59,167.39) .. (443.71,165.58) .. controls (445.83,163.77) and (449.05,163.99) .. (450.91,166.06) .. controls (452.77,168.14) and (452.56,171.29) .. (450.44,173.1) -- cycle ;
            \draw [color={rgb, 255:red, 0; green, 0; blue, 255 }  ,draw opacity=1 ]   (410.02,134.3) -- (418.52,135.17) ;
            \draw [color={rgb, 255:red, 0; green, 0; blue, 255 }  ,draw opacity=1 ]   (418.52,135.17) -- (418.11,143.57) ;
            \draw  [color={rgb, 255:red, 0; green, 0; blue, 255 }  ,draw opacity=1 ][fill={rgb, 255:red, 0; green, 0; blue, 255 }  ,fill opacity=1 ] (481.45,199.14) -- (478.11,208.92) -- (471.13,201.23) -- cycle ;
            \draw [color={rgb, 255:red, 255; green, 0; blue, 0 }  ,draw opacity=1 ][line width=2.25]    (42.11,316.5) .. controls (85.33,267.7) and (86.71,363.36) .. (42.32,316.5) ;
            \draw  [color={rgb, 255:red, 255; green, 0; blue, 0 }  ,draw opacity=1 ][fill={rgb, 255:red, 255; green, 0; blue, 0 }  ,fill opacity=1 ] (38.47,316.5) .. controls (38.47,314.29) and (40.19,312.49) .. (42.32,312.49) .. controls (44.44,312.49) and (46.16,314.29) .. (46.16,316.5) .. controls (46.16,318.71) and (44.44,320.51) .. (42.32,320.51) .. controls (40.19,320.51) and (38.47,318.71) .. (38.47,316.5) -- cycle ;
            \draw  [color={rgb, 255:red, 255; green, 0; blue, 0 }  ,draw opacity=1 ][fill={rgb, 255:red, 255; green, 0; blue, 0 }  ,fill opacity=1 ] (74.89,311.51) -- (79.34,319.55) -- (71,319.83) -- cycle ;
            \draw [color={rgb, 255:red, 255; green, 0; blue, 0 }  ,draw opacity=1 ][line width=2.25]    (175.5,316.67) .. controls (218.72,267.88) and (220.09,363.53) .. (175.7,316.67) ;
            \draw  [color={rgb, 255:red, 255; green, 0; blue, 0 }  ,draw opacity=1 ][fill={rgb, 255:red, 255; green, 0; blue, 0 }  ,fill opacity=1 ] (171.85,316.67) .. controls (171.85,314.46) and (173.58,312.67) .. (175.7,312.67) .. controls (177.83,312.67) and (179.55,314.46) .. (179.55,316.67) .. controls (179.55,318.89) and (177.83,320.68) .. (175.7,320.68) .. controls (173.58,320.68) and (171.85,318.89) .. (171.85,316.67) -- cycle ;
            \draw  [color={rgb, 255:red, 255; green, 0; blue, 0 }  ,draw opacity=1 ][fill={rgb, 255:red, 255; green, 0; blue, 0 }  ,fill opacity=1 ] (208.27,311.68) -- (212.73,319.73) -- (204.38,320) -- cycle ;
            \draw [color={rgb, 255:red, 255; green, 0; blue, 0 }  ,draw opacity=1 ][line width=2.25]    (472.56,316.25) .. controls (515.77,267.45) and (517.15,363.11) .. (472.76,316.25) ;
            \draw  [color={rgb, 255:red, 255; green, 0; blue, 0 }  ,draw opacity=1 ][fill={rgb, 255:red, 255; green, 0; blue, 0 }  ,fill opacity=1 ] (468.91,316.25) .. controls (468.91,314.04) and (470.64,312.24) .. (472.76,312.24) .. controls (474.89,312.24) and (476.61,314.04) .. (476.61,316.25) .. controls (476.61,318.47) and (474.89,320.26) .. (472.76,320.26) .. controls (470.64,320.26) and (468.91,318.47) .. (468.91,316.25) -- cycle ;
            \draw  [color={rgb, 255:red, 255; green, 0; blue, 0 }  ,draw opacity=1 ][fill={rgb, 255:red, 255; green, 0; blue, 0 }  ,fill opacity=1 ] (505.33,311.26) -- (509.79,319.3) -- (501.44,319.58) -- cycle ;
            
            \draw (88,212) node [anchor=north west][inner sep=0.75pt]  [rotate=-317.37] [xscale = 1.7, yscale = 1.7] {$\textcolor[rgb]{0,0,1}{1}$};
            \draw (236,270) node [anchor=north west][inner sep=0.75pt]  [rotate=-294.34] [xscale = 1.7, yscale = 1.7] {$\textcolor[rgb]{0,0,1}{2}$};
            \draw (505,190) node [anchor=north west][inner sep=0.75pt]  [rotate=-42.76] [xscale = 2, yscale = 2] {\color{blue}$\frac{p-1}{2}$};
            \draw (346.91,235.47) node [anchor=north west][inner sep=0.75pt]  [xscale = 1.7, yscale = 1.7]  {$\dotsc $};
            \draw (88,306) node [anchor=north west][inner sep=0.75pt]  [color={rgb, 255:red, 255; green, 0; blue, 0 }  ,opacity=1 ,rotate=-1.55] [xscale = 1.7, yscale = 1.7] {$1$};
            \draw (220,306) node [anchor=north west][inner sep=0.75pt]  [color={rgb, 255:red, 255; green, 0; blue, 0 }  ,opacity=1 ,rotate=-1.55] [xscale = 1.7, yscale = 1.7] {$2$};
            \draw (510,300) node [anchor=north west][inner sep=0.75pt]  [color={rgb, 255:red, 255; green, 0; blue, 0 }  ,opacity=1 ,rotate=-1.55] [xscale = 2, yscale = 2] {$\frac{p-1}{2}$};
            \draw (275,55) node [anchor=north west][inner sep=0.75pt]  [xscale = 1.7, yscale = 1.7]  {\color{red}$1$};
            \draw (19,306) node [anchor=north west][inner sep=0.75pt] [xscale = 1.7, yscale = 1.7]   {\color{red}$1$};
            \draw (153,306) node [anchor=north west][inner sep=0.75pt] [xscale = 1.7, yscale = 1.7]   {\color{red}$1$};
            \draw (450,306) node [anchor=north west][inner sep=0.75pt] [xscale = 1.7, yscale = 1.7]  {\color{red}$1$};
            
            \end{tikzpicture}
           }
           \end{center}
    \caption{The intersection of $Y$ and $\Delta_{2, p}^{par}$.}
    \label{fig:intersectionOfLoci}
\end{figure}
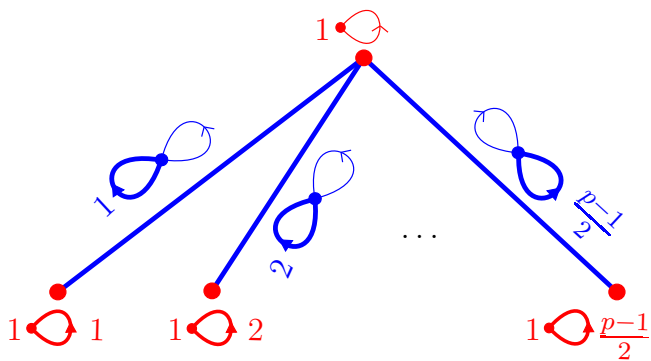

The cells in $X$ that are not in $Y \cup \Delta_{2,p}^{par}$ are the following:
\begin{enumerate}
    \item $D_{i}^{j}$ for $2 \leq i,j \leq (p-1)/2$ in entry $(1,3)$ of \cref{tab:pCoversGenusTwo}. There are $\frac{(p-3)^2}{4}$ of these.
    \item  $T^{i,j,k}$ in entry $(4,3)$ of \cref{tab:pCoversGenusTwo}. There are $\frac{(p-1)(p-5)}{12}$ of these cells.
    \item $T_{i,j,k}$ in entry $(5,3)$ of \cref{tab:pCoversGenusTwo}. There are  $\frac{(p-1)(p-5)}{12}$ of these cells as well.
\end{enumerate}
Note that the boundary of each of these remaining cells belongs to the complex $Y \cup \Delta_{2,p}^{par}$, so that in the quotient $X/(Y \cup \Delta_{2,p}^{par})$, each of these cells becomes a 2-dimensional sphere. We deduce that the quotient $X / (Y \cup \Delta_{2,p}^{par})$ is the wedge of $(p-1)(p-5)/6  + (p-3)^2 / 4$ spheres of dimension $2$. Since $Y \cup \Delta_{2,p}^{par}$ is contractible, we deduce that $X$ (and  consequently also $\Delta_{2,p}$) is homotopy equivalent to this quotient.
\end{proof}

    \addtocontents{toc}{\protect\setcounter{tocdepth}{0}}
    \section*{Acknowledgements}
    \addtocontents{toc}{\protect\setcounter{tocdepth}{1}}

    The authors thank Melody Chan and Martin Ulirsch and the anonymous referee for helpful comments. This work started during the authors' visit to Durham University, which was partially sponsored by Future Leaders Fellowship programme MR/S034463/2, being part of PH's Postdoctoral Fellowship. 
    The project was advanced further during FR's visit to RWTH Aachen University and YEM's visit to the University of T\"ubingen. We thank the aforementioned institutions for their generous hospitality.

    YEM, FR, and CY gratefully acknowledge support from the Deutsche Forschungsgemeinschaft (DFG, German Research Foundation) SFB-TRR 195 ``Symbolic Tools in Mathematics and their Application". PH was supported by UK Research and Innovation under the Future Leaders Fellowship programme MR/S034463/2 as a Postdoctoral Fellow at Durham University, and the JSPS Postdoctoral Fellowship 23769 and KAKENHI 23KF0187 as a Postdoctoral Fellow at the University of Tsukuba. FR and PS acknowledge support from the DFG through SFB-TRR 326 ``Geometry and Arithmetic of Uniformized Structures", project number 444845124, as well as Sachbeihilfe ``From Riemann surfaces to tropical curves (and back again)", project number 456557832. PS was also partially funded by the Deutscher Akademischer Austauschdienst (DAAD) under the scholarship program ``Forschungsstipendien - Promotionen in Deutschland'', 2020/21 - No. 57507871. There is no additional computational data associated to this paper.
    
    \newcommand{\etalchar}[1]{$^{#1}$}

\end{document}